\documentclass[12pt]{article}

\usepackage{epsf,epsfig,amsmath,amsfonts,upgreek}
\usepackage{graphicx}
\usepackage[multiple]{footmisc}

\newtheorem{thm}{Theorem}

\newtheorem{ex}{Example}

\newtheorem{lem}[thm]{Lemma}
\newtheorem{de}{Definition}
\newtheorem{rem}{Remark}
\newtheorem{ass}{Assumption}

\newcommand{\be}{\begin{eqnarray*}}
\newcommand{\ben}{\begin{eqnarray}}
\newcommand{\ee}{\end{eqnarray*}}
\newcommand{\een}{\end{eqnarray}}
\newcommand{\aln}{\begin{align}}
\newcommand{\ealn}{\end{align}}

\newcommand{\half}{\frac{1}{2}}

\newcommand{\F}{{\mathcal F}}

\newcommand{\R}{{\mathbb R}}

\renewcommand{\bar}[1]{\overline{#1\mkern-0mu}\mkern0mu }
\newcommand{\ubar}[1]{\underline{#1\mkern-0mu}\mkern0mu }

\renewcommand{\tilde}{\widetilde}

\newcommand{\bbe}{\mathbb{E}}
\newcommand{\bbp}{\mathbb{P}}

\newcommand{\bbr}{\mathbb{R}}
\newcommand{\bbs}{\mathbb{S}}

\newcommand{\bcal}{\mathcal{B}}

\newcommand{\ical}{\mathcal{I}}

\newcommand{\fcal}{\mathcal{F}}

\newcommand{\tcal}{\mathcal{T}}

\newcommand{\gc}{{\mathfrak C}}
\newcommand{\gzc}{{\mathfrak Z} _{\mathrm c}}

\newcommand{\gzb}{{\mathfrak Z} _\beta}
\newcommand{\gzo}{{\mathfrak Z} _0}
\newcommand{\gz}{{\mathfrak z}}

\newcommand{\liota}{{\ubar{\upiota}}}
\newcommand{\riota}{{\bar{\upiota}}}

\newcommand{\zminus}{z_\ominus}
\newcommand{\zplus}{z_\oplus}

\newcommand{\Cle}{C_{\mathrm l}}
\newcommand{\Cri}{C_{\mathrm r}}
\newcommand{\Dle}{D_{\mathrm l}}
\newcommand{\Dri}{D_{\mathrm r}}

\newcommand{\icalo}{\mathring{\ical}}

\usepackage{color}

\usepackage{tikz}
\usetikzlibrary{trees}

\begin{document}

\title{Discretionary stopping of stochastic differential
equations with generalised drift}

\author{
{\sc Mihail Zervos}\footnote{Department of Mathematics,
London School of Economics, Houghton Street, London
WC2A 2AE, UK, \texttt{mihalis.zervos@gmail.com}.}, \ 
{\sc Neofytos Rodosthenous}\footnote{School of
Mathematical Sciences, Queen Mary University
of London, London E1 4NS, UK, 
\texttt{n.rodosthenous@qmul.ac.uk}.
Research supported by EPSRC under Grant Number
EP/P017193/1.}, \\
{\sc Pui Chan Lon}
\ and
{\sc Thomas Bernhardt}\footnote{Department of
Mathematics, University of Michigan, 1851 East Hall,
530 Church Street, Ann Arbor, MI 48109,
\texttt{bernt@umich.edu}.}
}
\maketitle

\begin{abstract}
We consider the problem of optimally stopping a
general one-dimensional stochastic differential equation
(SDE) with generalised drift over an infinite time horizon.
First, we derive a complete characterisation of the solution
to this problem in terms of variational inequalities.
In particular, we prove that the problem's value function
is the difference of two convex functions and satisfies
an appropriate variational inequality in the sense of
distributions.
We also establish a verification theorem that is the
strongest one possible because it involves only the
optimal stopping problem's data.
Next, we derive the complete explicit solution to the
problem that arises when the state process is a skew
geometric Brownian motion and the reward function is
the one of a financial call option.
In this case, we show that the optimal stopping strategy
can take several qualitatively different forms, depending
on parameter values.
Furthermore, the explicit solution to this special case
shows that the so-called ``principle of smooth fit''
does not hold in general for optimal stopping problems
involving solutions to SDEs with generalised drift.
\end{abstract}

\newpage
\section{Introduction}

We consider the optimal stopping of the one-dimensional
SDE with generalised drift
\ben
X_t = x + \int _\liota^\riota L_t^z \, \nu (dz) +
\int _0^t \sigma (X_s) \, dW_s , \quad x \in \icalo ,
\label{X}
\een
where $L^z$ is the symmetric local time of $X$ at level $z$,
$W$ is a standard one-dimensional Brownian motion and
$\icalo = \mbox{} ]\liota, \riota[$ is the interior of a given
interval $\ical \subseteq [-\infty, \infty]$.
We assume that the signed Radon measure $\nu$ and the
Borel-measurable function $\sigma: \icalo \rightarrow \bbr
\setminus \{ 0 \}$ satisfy suitable conditions ensuring that
the SDE (\ref{X}) has a weak solution $\bigl( \Omega, \fcal,
(\fcal_t), \bbp_x, W, X \bigr)$ that is unique in the sense of
probability law up to a possible explosion time at which $X$
hits the boundary $\{ \liota, \riota \}$ of $\ical$ (see
Assumption~\ref{A1} in Section~\ref{sec:SDE}).
If the boundary point $\liota$ (resp., $\riota$) is inaccessible,
then the interval $\ical$ is open from the left (resp., open
from the right).
On the other hand, if the boundary point $\liota$ (resp.,
$\riota$) is not inaccessible, then we assume that it is
absorbing and the interval $\ical$ is closed from the left
(resp., closed from the right).
Comprehensive studies of these SDEs as well as
relevant literature surveys can be found in Engelbert
and Schmidt~\cite{ES91}, and Lejay~\cite{Lej06}.

In the special case when $\nu$ is absolutely continuous,
namely, when
\be
\nu (dx) = \frac{b(x)}{\sigma^2 (x)} \, dx ,
\ee 
for a Borel-measurable function $b: \icalo \rightarrow \bbr$
satisfying suitable integrability conditions,
an application of the occupation times formula shows
that the solution to (\ref{X}) admits the expression
\ben
X_t = x + \int _0^t b(X_s) \, ds + \int _0^t \sigma (X_s)
\, dW_s , \quad x \in \icalo . \label{X-usual}
\een
In view of this simple observation, we can see that
usual SDEs are special cases of SDEs with generalised
drift.
The skew Brownian motion, which is characterised by
the choices
\be
\nu (dx) = \beta \, \delta_0 (dx) , \quad \sigma \equiv 1
\quad \text{and} \quad \ical = \bbr ,
\ee 
where $\beta \in \mbox{} ]{-1}, 1[ \mbox{} \setminus \{ 0 \}$
and $\delta_0 (dx)$ is the Dirac probability measure that
assigns mass 1 at $\{ 0 \}$, and corresponds to the SDE
\ben
X_t = x + \beta L_t^0 + W_t , \quad x \in \bbr ,
\label{X-BM}
\een
is a fundamental example of an SDE with generalised
drift
(see It\^{o} and McKean~\cite[Problem~4.2.1]{IMcK74},
Harrison and Shepp~\cite{HS81},
Lejay~\cite{Lej06}, and several references therein).
A further important example is the skew geometric
Brownian motion, which is characterised by the
choices
\be
\nu (dx) = \frac{b}{\sigma^2 x} \, dx + \beta \, \delta_z
(dx) , \quad \sigma (x) = \sigma x \quad \text{and}
\quad \ical = \mbox{} ]0, \infty[ ,
\ee 
where $b \in \bbr$, $\beta \in \mbox{} ]{-1}, 1[
\mbox{} \setminus \{ 0 \}$, $\sigma \neq 0$ are
constants and $\delta_z (dx)$ is the Dirac probability
measure that assigns mass 1 at $\{ z \}$, for
some $z>0$, and corresponds to the SDE
\ben
dX_t = bX_t \, dt + \beta \, dL_t^z + \sigma X_t \, dW_t ,
\quad X_0 = x > 0 . \label{XGBM}
\een
Furthermore, an SDE with generalised drift whose
dynamics identify with the dynamics of a usual SDE
away from a finite number of points at which it
exhibits a skew behaviour is discussed
in Example~\ref{ex:levels} below.
At this point, we emphasise that the processes
$L^z$, $z \in \icalo$, in SDEs with generalised drift are
the local times of the solution $X$ to the corresponding
SDEs and not the local times of the driving Brownian
motion $W$ (see also Example~\ref{ex:skewBM} where
this observation is considered further in the context of
the skew Brownian motion given by (\ref{X-BM})).

The objective of the optimal stopping problem that we
study aims at maximising the performance criterion
\ben
\bbe _x \left[ \exp \left( - \int _0^\tau r(X_s) \, ds \right)
f(X_\tau) {\bf 1} _{\{ \tau < \infty \}}\right] \label{criterion}
\een
over all stopping times $\tau$, where the positive
Borel-measurable discounting rate function $r$ satisfies
Assumption~\ref{A2} in Section~\ref{sec:SDE}, while the
positive and possibly unbounded reward function $f$
satisfies Assumption~\ref{A3} in Section~\ref{sec:VI}.
To the best of our knowledge, there exist no results
in the literature addressing the solvability of such a
problem by means of variational inequalities when
$\nu$ is not absolutely continuous, even in special
cases.
We derive a complete characterisation of the solution
to this problem in terms of variational inequalities.
In particular, we prove that the value function $v$
of the optimal stopping problem associated with
(\ref{X}) and (\ref{criterion}) is the difference of two
convex functions and satisfies the variational 
inequality
\ben
\max \left\{ \frac{1}{2} \sigma^2 (x) p_-' (x) \left(
\frac{v_-'}{p_-'} \right) ' (dx) - r(x) v(x) \, dx ,
\ f(x) - v(x) \right\} = 0 , \label{v-VI}
\een
in the sense of distributions (see Definition~\ref{distr}
and Theorem~\ref{VThm}.(I)-(II) in Section~\ref{sec:VI}),
where $p$ is the scale function of the diffusion
associated with the SDE (\ref{X}).
We also establish a verification theorem that is the
strongest one possible because it involves only the
optimal stopping problem's data.
In particular, we derive a simple necessary and sufficient
condition for a solution to (\ref{v-VI}) to identify
with the problem's value function (see
Theorem~\ref{VThm}.(III)).

The second main contribution of the paper is to derive
the complete explicit solution to the special case that
arises if $f(x) = (x-K)^+$, for some constant $K>0$,
and $X$ is a skew geometric Brownian motion.
In this case, the SDE (\ref{XGBM}), has a unique
non-explosive strong solution.
Given such a solution $X$, which exists on any given
filtered probability space $\bigl( \Omega, \fcal, (\fcal_t),
\bbp \bigr)$ satisfying the usual conditions and
supporting a standard one-dimensional
$(\fcal_t)$-Brownian motion $W$, the value function
of the discretionary stopping problem that we solve is
defined by
\ben
v(x) = \sup _{\tau \in \tcal} \bbe \left[ e^{-r \tau}
(X_\tau - K)^+ {\bf 1} _{\{ \tau < \infty \}} \right] ,
\label{vEx}
\een
where $\tcal$ is the family of all $(\fcal_t)$-stopping
times and $r, K > 0$ are constants (we write $\bbe$
in place of $\bbe_x$ because we consider strong
rather than weak solutions here).
We prove that the optimal stopping strategy
can take several qualitatively different forms, depending
on parameter values (see
Theorems~\ref{propSol1}-\ref{prop:gz} and
Figures~4-13).
In contrast to the optimal stopping of an SDE with
absolutely continuous drift and reward function
such as the one of a financial call option, the optimal
stopping region may involve two distinct components,
one of which may be an isolated point.
Furthermore, the analysis of this problem shows that
the so-called ``principle of smooth fit'' does not hold
even in the case of a ``right-sided'' optimal stopping
strategy in the sense that none of the functions
\begin{gather}
v_-' (x) = \lim _{\varepsilon \downarrow 0}
\frac{v(x) - v(x-\varepsilon)}{\varepsilon} , \quad
\frac{v_-' (x)}{p_-' (x)} = \lim _{\varepsilon \downarrow 0}
\frac{v(x) - v (x-\varepsilon)} {p(x) - p(x-\varepsilon)}
\nonumber \\
\text{or} \quad
\frac{v_-' (x)}{\psi_-' (x)} = \lim _{\varepsilon \downarrow 0}
\frac{v(x) - v (x-\varepsilon)} {\psi (x) - \psi (x-\varepsilon)}
\label{SF}
\end{gather}
is continuous, where $p$ (resp., $\psi$) is the scale
function (resp., the increasing minimal excessive
function) of the diffusion associated with the SDE
(\ref{XGBM}) (see Remarks~\ref{rem:SF}
and~\ref{rem:SFex}).

To derive the solution to the optimal stopping problem
defined by (\ref{XGBM}) and (\ref{vEx}), we need
to consider a partition of the set $\bbr \times (\bbr
\setminus \{0\}) \times \mbox{} ]0, \infty[ \mbox{} \times
\mbox{} \bigl( ]{-1}, 1[ \setminus \{ 0 \} \bigr)$ in which
the parameter vector $(b, \sigma, z, \beta)$ takes
values into four sets (see Cases~(I)-(IV) of
Lemma~\ref{lem:psi}).
In Cases~(I) and~(II), $\beta \in \mbox{} ]{-1}, 0[$
and $\psi$ is convex.
In Cases~(III) and~(IV), $\beta \in \mbox{} ]{-1}, 0[$
and $\beta \in \mbox{} ]0, 1[$, respectively, and
$\psi$ fails to be convex.
The best part of our analysis' complexity (including
the whole of Section~\ref{sec:F}) is due to the facts
that (a) the optimal strategy takes several qualitatively
different forms, and (b) we establish necessary and
sufficient conditions on the problem's data that
differentiate between the different possible cases
without leaving any ``gap'' in the parameter space.
The solution to the problem in the easier Cases~(I),
(II) and~(IV) has been presented in the PhD thesis
Lon~\cite{Lon}.
Although the analysis of Case~(III) is not in itself
harder, linking it with Cases~(I) and~(II) with
necessary and sufficient conditions on the problem's
data requires rather tedious analysis, due to the
complex structure of $\psi$ (see Figure~2).
The solution to all cases was announced without
proofs in the conference paper Lon, Rodosthenous
and Zervos~\cite{LRZ}.

Variational inequalities take center stage in the continuous
time optimal stopping theory because they are efficient
for the investigation of specific problems.
In the context of this paper, they can be used to easily
identify critical parts of the state space $\ical$ that are
subsets of the so-called waiting region in a systematic
way that involves no guesswork.
For instance, the regularity of $v$ implies that all points
at which the reward function $f$ is discontinuous as well
as all ``minimal'' intervals in which $f$ cannot be
expressed as the difference of two convex functions
(e.g., intervals in which $f$ has the regularity of a
Brownian sample path) should be parts of the closure
of the waiting region.
For further analysis and discussion in this direction, see
Remark~\ref{rem:VIuse} at the end of Section~\ref{sec:VI}.
Beyond its usefulness in identifying optimal stopping
strategies, the variational inequality characterisation is
also very effective in verifying whether a candidate
function identifies with the value function of a specific
problem because, in the context of SDEs driven by
a Brownian motion, it has a local character in the sense
that it involves only derivatives.

The solution to optimal stopping problems using classical
solutions to variational inequalities has been extensively
studied.
Results in Friedman~\cite[Chapter~16]{F},
Bensoussan and Lions~\cite[Chapter~3]{BL},
{\O}ksendal~\cite[Chapter 10]{O} and Peskir and
Shiryaev~\cite{PS}, listed in chronological order, typically
make strong regularity assumptions on the problem data
(e.g., the problem data are assumed to be Lipschitz
continuous).
To relax such assumptions, {\O}ksendal and Reikvam~\cite{OR}
and Bassan and Ceci~\cite{BC} have considered
viscosity solutions to the variational inequalities associated
with the optimal stopping problems that they study.
Results with minimal assumptions on the problem data,
such as the ones that we derive here for the optimal
stopping problem associated with (\ref{X}) and (\ref{criterion}),
have been obtained by Lamberton~\cite{Lam} and
Lamberton and Zervos~\cite{LZ} who consider the
optimal stopping of the SDE (\ref{X-usual})
over a finite and an infinite time horizon, respectively.
Recently, the solution to suitable optimal stopping
problems by means of variational inequalities has
been used to characterise the boundary of Root's
solution to the Skorokhod embedding problem (see
Cox and Wang~\cite{CW},
Cox, Ob\l\'{o}j and Touzi~\cite{COT},
and references therein).
Furthermore, variational inequalities arise most
naturally in the study of optimal stopping problems
involving controlled stochastic processes (e.g., see
Bensoussan and Lions~\cite[Chapter~4]{BL},
Krylov~\cite[Chapters~3,~6]{K},
Bene\v{s}~\cite{Benes},
Davis and Zervos~\cite{DZ},
Karatzas and Sudderth~\cite{KS},
listed in chronological order,
as well as many more recent contributions).

There exist few references in the literature addressing
special cases of a general problem such as the one
we study.
Peskir~\cite{P} considered the validity of the ``principle
of smooth fit'' in terms of derivatives such as the first
two ones in (\ref{SF}).
In particular, failure of the ``principle of smooth fit'' was
exhibited by Examples~2.2 and~3.1 in this reference.
These examples involve the optimal stopping of the
processes $X = F(B)$, where $B$ is a standard Brownian
motion absorbed at the boundaries of $[{-1}, 1]$ and
\be
F(x) = \begin{cases} x^{1/3} , & \text{if } x \in [0,1] , \\
-|x|^{1/3} , & \text{if } x \in [{-1},0[ , \end{cases}
\quad \text{or} \quad
F(x) = \begin{cases} \sqrt{x} , & \text{if } x \in [0,1] , \\
-x^2 , & \text{if } x \in [{-1},0[ . \end{cases}
\ee
These processes are skew diffusions that cannot
be associated with solutions to SDEs with generalised
drift because It\^{o}-Tanaka's formula cannot be applied
due to the fact that $F_+'(0) = \infty$.
In the second paragraph of Section~3.2 in
Peskir~\cite{P}, a further similar example exhibiting
failure of the ``principle of smooth fit'' is briefly
discussed.
If we make a suitable choice for $\sigma$, such as
$\sigma \equiv 1$, then we can associate the diffusion
of this example with the solution to an SDE with
generalised drift by replicating the analysis in
Example~\ref{ex:skewBM} below.

Other references have considered the optimal
stopping of a skew Brownian motion, which is given
by the strong solution to the SDE (\ref{X-BM}),
with the objective to maximise the performance
index given by (\ref{criterion}) for $r>0$ being a
constant and for $f$ being an increasing function
associated with ``right-sided'' optimal stopping
strategies.
Crocce and Mordecki~\cite{CM} studied the
validity of the ``principle of smooth fit'' in terms of
derivatives such as the ones given by (\ref{SF}), and
presented two examples with optimal stopping
strategies such as the ones in Theorem~\ref{propSol1}
that are illustrated by Figures~5 and~10 below.
Alvarez and Salminen \cite{AS} derived
sufficient conditions on $f$ that are associated
with optimal stopping strategies of the same qualitative
nature as the ones in Theorems~\ref{propSol1}
and~\ref{prop:gz} that are illustrated by
Figures~10, 12 and~13 below.
The analysis in these references is based on Dynkin's
characterisation of an optimal stopping problem's
value function as the minimal excessive majorant
of its reward function and the Martin representation
theory of excessive functions.
In contrast to variational inequalities, which involve
only the problem's primary data, this approach has
a non-local character or it requires conditions involving
the elements of the set
\be
\left\{ x \in \ical \ \Big| \ \ \frac{f(x)}{\psi (x)} =
\sup _{u \in \ical} \frac{f(u)}{\psi (u)} \right\} .
\ee
As a consequence, its applicability has largely been
limited to problems with ``one-sided'' optimal stopping
strategies because, with notable exceptions such as
the ones associated with a Brownian motion or a
geometric Brownian motion, the minimal excessive
functions are not in general expressible in terms
of elementary functions.

Beyond its contributions to the optimal stopping theory,
the present paper has been motivated by applications
to the optimal timing of investment decisions involving
an underlying asset price or economic indicator.
In mathematical finance and the theory of real options,
such time series are typically modelled by SDEs driven
by a standard Brownian motion or, more generally,
a L\'{e}vy process.
A skew geometric Brownian motion or, more generally,
SDEs such as the ones considered in
Example~\ref{ex:levels} can be used to model
asset prices and economic indicators that
exhibit support and resistance levels\footnote{``Support
is a level or area on the chart under the market where
buying interest is sufficiently strong to overcome selling
pressure.
As a result, a decline is halted and prices turn back
again... Resistance is the opposite of support.''
(Murphy~\cite{Mu99})}
(see H\"{a}m\"{a}l\"{a}inen~\cite{Ha} for a recent
survey of studies focusing on such directional
predictability).
Indeed, the skew geometric Brownian motion (\ref{XGBM})
behaves like a standard geometric Brownian motion,
except that the sign of each excursion from $z$ is chosen
using an independent Bernoulli random variable of
parameter $p = \frac{1}{2} (\beta + 1)$, namely, $\bbp
(X_t > z \mid t > T_z) = p$, where $T_z$ is the first
hitting time of $z$.
Recently, several authors have studied financial models
involving SDEs with generalised drift
(e.g., see
Corns and Satchell~\cite{CS},
Decamps, Goovaerts and
Schoutens~\cite{DGS06a,DGS06b},
Rosello~\cite{Ro},
and references therein).
Alternatively, a skew geometric Brownian motion
can be used to capture phenomena of bounces and
sinks that are exhibited by financial firms in distress
(see Nilsen and Sayit~\cite{NS}).

The paper is organised as follows.
In Section~\ref{sec:SDE}, we derive an analytic
characterisation of the minimal excessive functions
of the one-dimensional diffusion associated with the
SDE (\ref{X}). 
In Section~\ref{sec:VI}, we establish a complete
characterisation of the solution to the general optimal
stopping problem given by (\ref{X}) and (\ref{criterion})
in terms of variational inequalities.
Section~\ref{sec:excessive} presents a study of a
skew geometric Brownian motion's minimal excessive
functions.
In Section~\ref{sec:F}, we prove a couple of technical
results that will facilitate the streamlining of the
presentation of the solution to the optimal stopping
problem defined by (\ref{XGBM}) and (\ref{vEx}).
We present the complete solution to this problem in
Section~\ref{sec:GBMsol}.
Finally, the proofs of all results stated in
Section~\ref{sec:GBMsol} are collected in
Section~\ref{sec:proofs}.

\section{The SDE (\ref{X}) and its associated minimal
excessive functions}
\label{sec:SDE}

We start with the following assumption.\footnote{We
denote by $\bcal (\icalo)$ the Borel $\sigma$-algebra
on $\icalo$.}

\begin{ass} \label{A1} {\rm
The function $\sigma: \icalo \rightarrow \bbr$ is
Borel-measurable,
\begin{align*}
\sigma (x) \neq 0 \text{ for all } x \in \icalo
\quad \text{and} \quad
\int _{\ubar{x}}^{\bar{x}} \sigma^{-2} (u) \, du < \infty
\text{ for all } \ubar{x} < \bar{x} \text{ in } \icalo .
\end{align*}
Also, $\nu$ is a signed Radon measure on
$\bigl( \icalo , \bcal (\icalo) \bigr)$ such that
$\nu \bigl( \{ z \} \bigr) \in \mbox{} ]{-1},1[$.
} \mbox{}\hfill$\Box$ \end{ass}

\noindent
In the presence of this assumption, the SDE (\ref{X})
has a weak solution that is unique in the sense of
probability law (see Engelbert and
Schmidt~\cite[Theorems~4.35 and~4.37]{ES91}).
In particular, given any initial point $x \in \icalo$,
there is a collection $\bbs_x = \bigl( \Omega, \fcal,
(\fcal_t), \bbp_x, W, X \bigr)$ such that $(\Omega, \fcal,
(\fcal_t), \bbp_x)$ is a filtered probability space
satisfying the usual conditions,  $W$ is a standard
$(\fcal_t)$-Brownian motion and $X$ is a continuous
$(\fcal_t)$-adapted stochastic process such that
(\ref{X}) holds true in the stochastic interval
$[0, T_\liota \wedge T_\riota[$, where
\ben
T_y = \inf \{ t \geq 0 \mid \ X_t = y \} ,
\quad \text{for } y \in [\liota, \riota] , \label{1sthit}
\een
with the usual assumption that $\inf \emptyset
= \infty$.
We assume that either of the endpoints $\liota$,
$\riota$ is either inaccessible or absorbing.
Accordingly, if $\liota$ (resp., $\riota$) is absorbing,
then $X_t = \liota$ (resp., $X_t = \riota$) for all
$t \geq T_\liota$ (resp., $t \geq T_\riota$).

The scale function of the diffusion associated with
the SDE (\ref{X}) is the unique, up to a strictly increasing
affine transformation, continuous strictly increasing
function $p: \ical \rightarrow \bbr$ that satisfies
\be
\bbp_x (T_{\bar{x}} < T_{\ubar{x}}) = 1 - \bbp_x
(T_{\ubar{x}} < T_{\bar{x}}) = \frac{p(x) - p({\ubar{x}})}
{p({\bar{x}}) - p({\ubar{x}})} ,
\ee
for all points $\ubar{x} < x < \bar{x}$ in $\ical$. 
The restriction of $p$ to $\icalo$ is the difference of
two convex functions\footnote{A function
$g: \icalo \rightarrow \bbr$ is the difference of two
convex functions if and only if it is absolutely continuous
with left-hand derivative that is a function of
finite variation.
Given such a function, we denote by
$g_\pm'$ its right-hand and left-hand side first
derivatives, which are defined by
\be
g_+' (x) = \lim _{\varepsilon \downarrow 0}
\frac{g(x+\varepsilon) - g(x)}{\varepsilon}
\quad \text{and} \quad
g_-' (x) = \lim _{\varepsilon \downarrow 0}
\frac{g(x) - g(x-\varepsilon)}{\varepsilon} ,
\ee
and by $g'' (dx)$ the measure that identifies with
its second distributional derivative.}
and satisfies the ordinary differential equation (ODE)
\ben
p'' (dx) = - \bigl[ p_+' (x) + p_-' (x) \bigr] \, \nu (dx)
\label{p'ODE}
\een
in the sense that
\ben
p_+' (\bar{x}) - p_+' (\ubar{x}) =
- \int_{]\ubar{x}, \bar{x}]} \bigl[ p_+' (z) + p_-' (z) \bigr]
\, \nu (dz) \quad \text{for all } \ubar{x} < \bar{x}
\text{ in } \icalo . \label{p''-sense}
\een
In particular, it is given by
\begin{align}
p_+' (x) & = e^{- 2 \nu ([x_1, x])} \prod _{z \in [x_1, x]}
\frac{1 - \nu \bigl( \{ z \} \bigr)}{1 + \nu \bigl( \{ z \} \bigr)}
e^{2 \nu (\{ z \})} , \quad \text{if } x \geq x_1 ,
\\
\text{and} \quad
p_+' (x) & = e^{2 \nu (]x, x_1[)} \prod
_{z \in \mbox{} ]x, x_1[} \frac{1 + \nu \bigl( \{ z \} \bigr)}
{1 - \nu \bigl( \{ z \} \bigr)} e^{- 2 \nu (\{ z \})} , \quad
\text{if } x < x_1 ,
\end{align}
where $x_1 \in \icalo$ is an arbitrary fixed point.
All these claims can be found, e.g., in Engelbert
and Schmidt~\cite[Section~4.3]{ES91}.
For future reference, we also note that these
expressions imply that
\ben
p_+' (x) = p_-' (x) \frac{1 - \nu \bigl( \{ x \} \bigr)}
{1 + \nu \bigl( \{ x \} \bigr)} \quad \text{for all }
x \in \icalo . \label{p-jumps}
\een

We will need the following real analysis
result.\footnote{We denote by $p(\icalo)$
the interval $]p(\liota), p(\riota)[$ and
by $\bcal \bigl( p(\icalo) \bigr)$ the Borel
$\sigma$-algebra on $p(\icalo)$.}

\begin{lem} \label{u-utilde}
Let $\tilde{u} : p(\icalo) \rightarrow \bbr$ be a
difference of two convex functions and define
$u(x) = \tilde{u} \bigl( p(x) \bigr)$, for $x \in \icalo$.
The following statements hold true:
\smallskip

\noindent {\rm (I)}
$u$ is the difference of two convex functions.
\smallskip

\noindent {\rm (II)}
The function $u_-' / p_-'$ is of finite variation.
Furthermore, the measure on $\bigl( \icalo, \bcal
(\icalo) \bigr)$ that identifies with the first distributional
derivative of the function $u_-' / p_-'$ is the image
of the measure $\tilde{u}''$ under the function $p^{-1}$,
namely, $(u_-' / p_-') ' (dx) = (\tilde{u}'' \circ p) (dx)$.
In particular,
\begin{gather}
\left( \frac{u_-'}{p_-'} \right) ' \bigl( [\ubar{x} ,
\bar{x}[ \bigr) \equiv \frac{u_-'}{p_-'} (\bar{x})
- \frac{u_-'}{p_-'} (\ubar{x}) = \tilde{u}'' \bigl(
[p(\ubar{x}), p(\bar{x})[ \bigr) \quad \text{for all }
\ubar{x} < \bar{x} \text{ in } \icalo \nonumber \\
\text{and} \quad
\tilde{u}'' \bigl( [\ubar{q}, \bar{q}[ \bigr) \equiv
\tilde{u}_-' (\bar{q}) - \tilde{u}_-' (\ubar{q}) =
\left( \frac{u_-'}{p_-'} \right) ' \bigl( [p^{-1} (\ubar{q}) ,
p^{-1} (\bar{q})[ \bigr) \quad \text{for all } \ubar{q}
< \bar{q} \text{ in } p(\icalo) . \nonumber
\end{gather}

\noindent {\rm (III)}
If $\tilde{u}$ has absolutely continuous first
derivative $\tilde{u}'$ $(= \tilde{u}_-' = \tilde{u}_+')$,
then $u_-' / p_-'$ is absolutely continuous
and\,\footnote{\label{foot1}
In this part of the lemma, we use the same notation
for the signed Radon measure that identifies with the second
distributional derivative of $\tilde{u}$ as well as for
the Radon-Nikodym derivative of this measure
with respect to the Lebesgue measure, namely,
we write $\tilde{u}'' (dq) = \tilde{u}'' (q) \, dq$.
We refer to this footnote whenever we make
such an abuse of notation; confusion is unlikely
to occur.}
\be
\left( \frac{u_-'}{p_-'} \right) ' (x) = p_-' (x) \tilde{u}''
\bigl( p(x) \bigr) , \quad \text{for } x \in \icalo .
\ee
\end{lem}
{\bf Proof.}
We first note that $u =  \tilde{u} \circ p$ is
absolutely continuous because it is the composition
of absolutely continuous functions and $p$ is
increasing.
Given any $x \in \icalo$,
\be
\lim _{\varepsilon \downarrow 0}
\frac{u (x) - u (x-\varepsilon)} {p(x) - p(x-\varepsilon)}
= \lim _{\varepsilon \downarrow 0}
\frac{\tilde{u} \bigl( p(x) \bigr) - \tilde{u} \bigl(
p(x-\varepsilon) \bigr)} {p(x) - p(x-\varepsilon)}
= \tilde{u}_-' \bigl( p(x) \bigr) .
\ee
Combining this observation with the fact that the
limit $p_-' (x) = \lim _{\varepsilon \downarrow 0}
\frac{1}{\varepsilon} \bigl[ p(x) - p(x-\varepsilon) \bigr]$
exists, we can see that the limit $u_-' (x) = \lim
_{\varepsilon \downarrow 0} \frac{1}{\varepsilon}
\bigl[ u(x) - u(x-\varepsilon) \bigr]$ exists for all
$x \in \icalo$.
Given any points $\ubar{x} < \bar{x}$ in $\icalo$,
we use the change of variables formula (e.g., see
Revuz and Yor~\cite[Proposition~0.4.10]{RY}) to
calculate
\begin{align}
\frac{u_-'}{p_-'} (\bar{x}) - \frac{u_-'}{p_-'}
(\ubar{x}) & = \tilde{u}_-' \bigl( p(\bar{x})
\bigr) - \tilde{u}_-' \bigl( p(\ubar{x}) \bigr)
\nonumber \\
& = \int _{[p(\ubar{x}), p(\bar{x})[} \tilde{u}''
(dq) = \int _{[\ubar{x}, \bar{x}[} (\tilde{u}'' \circ p)
(dx) , \nonumber
\end{align}
and (II) follows.
Furthermore, (I) follows from the absolute
continuity of $u$ and the fact that $u_-'$
is the product of the finite variation functions
$p_-'$ and $u_-' / p_-'$.

Finally, if $\tilde{u}'$ is absolutely continuous
(see also footnote~\ref{foot1}), then
\be
\frac{u_-'}{p_-'} (\bar{x}) - \frac{u_-'}{p_-'} (\ubar{x})
= \int _{p(\ubar{x})}^{p(\bar{x})} \tilde{u}'' (q) \,
dq = \int _{\ubar{x}}^{\bar{x}} \tilde{u}'' \bigl( p
(x) \bigr) p_-' (x) \, dx \quad \text{for all } \ubar{x}
< \bar{x} \text{ in } \icalo ,
\ee
and (III) follows. 
\mbox{}\hfill$\Box$
\bigskip

Given a weak solution $\bbs_x = \bigl( \Omega, \fcal,
(\fcal_t), \bbp_x, W, X \bigr)$ to the SDE (\ref{X}),
the collection $\bigl( \Omega, \fcal, (\fcal_t), \bbp_x,
W, p(X) \bigr)$ is a weak solution to the SDE
\ben
d\tilde{X}_t  = \bigl( \sigma \circ p^{-1} \bigr) (\tilde{X}_t)
\, \bigl( p_-' \circ p^{-1} \bigr) (\tilde{X}_t) \ d\tilde{W}_t ,
\quad \tilde{X}_0 = p(x) \in p(\icalo) , \label{X-tilde}
\een
for $\tilde{W} = W$.
Conversely, given a weak solution $\tilde{\bbs}_x =
\bigl( \tilde{\Omega}, \tilde{\fcal}, (\tilde{\fcal}_t),
\tilde{\bbp}_x, \tilde{W}, \tilde{X} \bigr)$ to the SDE
(\ref{X-tilde}), the collection $\bigl( \tilde{\Omega},
\tilde{\fcal}, (\tilde{\fcal}_t), \tilde{\bbp}_x, \tilde{W},
p^{-1} (\tilde{X}) \bigr)$ is a weak solution to the SDE
(\ref{X}) for $W = \tilde{W}$.
These results, which are established in Engelbert
and Schmidt~\cite[Proposition~4.29]{ES91}),
will play a fundamental role in our analysis.

To proceed further, we consider the discounting
rate function $r$ appearing in (\ref{criterion}) and we
make the following assumption.

\begin{ass} \label{A2} {\rm
The function $r: \icalo \rightarrow \bbr_+$ is
Borel-measurable, uniformly bounded away from
$0$, namely, $r(x) \geq r_0$ for all $x \in \icalo$,
for some $r_0 > 0$, and such that
\be
\int _{\ubar{x}}^{\bar{x}} \frac{r(u)}{\sigma^{2} (u)} \, du < \infty
\quad \text{for all } \ubar{x} < \bar{x} \text{ in } \icalo .
\ee
} \vspace{-13.3mm}

\mbox{}\hfill$\Box$ \end{ass}
\smallskip

\noindent
The minimal $r$-excessive functions
$\varphi, \psi: \icalo \rightarrow \bbr_+$
of the diffusion associated with the SDE (\ref{X}) are the
unique, modulo multiplicative constants, functions
that satisfy
\begin{align}
\varphi (\bar{x}) &= \varphi (\ubar{x}) \, \bbe
_{\bar{x}} \left[ \exp \left( - \int _0^{T_{\ubar{x}}}
r(X_s) ds \right) \right] \text{ for all } \ubar{x}
<  \bar{x} \text{ in } \icalo \label{rphi} \\
\text{and} \quad
\psi (\ubar{x}) &= \psi (\bar{x}) \, \bbe
_{\ubar{x}} \left[ \exp \left( - \int _0^{T_{\bar{x}}}
r(X_s) ds \right) \right] \text{ for all } \ubar{x}
< \bar{x} \text{ in } \icalo , \label{rpsi}
\end{align}
where $T_{\ubar{x}}$, $T_{\bar{x}}$ are as in
(\ref{1sthit}) (see Borodin and
Salminen~\cite[Section~II.1]{BS}).

\begin{lem} \label{lem:exess}
The functions $\varphi, \psi: \icalo \rightarrow \bbr_+$
given by (\ref{rphi})--(\ref{rpsi}) are such
that $\varphi$ (resp., $\psi$) is strictly decreasing
(resp., increasing),
\begin{gather}
\text{if } \liota \text{ is absorbing, then } \varphi
(\liota) := \lim _{x \downarrow \liota} \varphi (x)
< \infty \text{ and }  \psi (\liota) := \lim
_{x \downarrow \liota} \psi (x) = 0 ,
\nonumber \\
\text{if } \riota \text{ is absorbing, then } \varphi
(\riota) := \lim _{x \uparrow \riota} \varphi (x) = 0
\text{ and } \psi (\riota) := \lim _{x \uparrow \riota}
\psi (x) < \infty , \nonumber \\
\text{and, if } \liota \text{ (resp., } \riota \text{)
is inaccessible, then } \lim _{x \downarrow \liota}
\varphi (x) = \infty \text{ (resp., } \lim
_{x   \uparrow \riota} \psi (x) = \infty \text{)} .
\nonumber
\end{gather}
Both of the functions $\varphi$ and $\psi$ are
absolutely continuous.
Furthermore, the functions $\varphi _-' / p_-'$
and $\psi _-' / p_-'$ are absolutely continuous
and the homogeneous ODE\,\footnote{
\label{foot2}
As in Lemma~\ref{u-utilde}.(III), we use here
$(\varphi _-' / p_-')'$ and $(\psi _-' / p_-')'$
to denote the Radon-Nikodym derivatives of
the measures that identify with the first distributional
derivatives of  the functions $\varphi _-' / p_-'$ and
$\psi _-' / p_-'$ with respect to the Lebesgue measure
(see also footnote~\ref{foot1}).}
\ben
\frac{1}{2} \sigma^2 (x) p_-' (x) \left(
\frac{g_-'}{p_-'} \right) ' (x) - r(x) g(x)
= 0 \nonumber
\een
is satisfied Lebesgue-a.e.\ in $\icalo$ for $g$ standing
for either $\varphi$ or $\psi$.
\end{lem}
{\bf Proof.}
In view of the results connecting the solvability
of (\ref{X}) with the solvability of (\ref{X-tilde})
that we have discussed above, we can see
that, given any $\ubar{x} < \bar{x}$ in $\icalo$,
\be
\frac{\varphi (\bar{x})}{\varphi (\ubar{x})} =
\bbe _{\bar{x}} \left[ \exp \left( - \int
_0^{\tilde{T}_{p(\ubar{x})}} (r \circ p^{-1})
(\tilde{X}_s) ds \right) \right] =
\frac{\tilde{\varphi} \bigl( p(\bar{x}) \bigr)}
{\tilde{\varphi} \bigl( p(\ubar{x}) \bigr)}
\quad \text{and} \quad
\frac{\psi (\bar{x})}{\psi (\ubar{x})} =
\frac{\tilde{\psi} \bigl( p(\bar{x}) \bigr)}
{\tilde{\psi} \bigl( p(\ubar{x}) \bigr)} ,
\ee
where
\be
\tilde{T}_y = \inf \{ t \geq 0 \mid \ \tilde{X}_t
\equiv p(X_t) = y \} , \quad \text{for } y \in \bigl[
p(\liota), p(\riota) \bigr] ,
\ee
and $\tilde{\varphi}, \tilde{\psi} : p(\icalo)
\rightarrow \bbr_+$ are the minimal
$(r \circ p^{-1})$-excessive functions of the
diffusion associated with the SDE (\ref{X-tilde}),
given by
\begin{align}
\tilde{\varphi} (\bar{q}) &= \tilde{\varphi}
(\ubar{q}) \, \tilde{\bbe} _{\bar{q}} \left[ \exp
\left( - \int _0^{\tilde{T}_{\ubar{q}}} (r \circ p^{-1})
(\tilde{X}_s) \, ds \right) \right] \text{ for all }
\ubar{q} <  \bar{q} \text{ in } p(\icalo)
\label{phi-til} \\
\text{and} \quad
\tilde{\psi} (\ubar{q}) &= \tilde{\psi} (\bar{q})
\, \tilde{\bbe} _{\ubar{q}} \left[ \exp \left( - \int
_0^{\tilde{T}_{\bar{q}}} (r \circ p^{-1})
(\tilde{X}_s) \, ds \right) \right] \text{ for all }
\ubar{q} < \bar{q} \text{ in } p(\icalo) .
\label{psi-til}
\end{align}
It follows that
\ben
\varphi = \tilde{\varphi} \circ p \quad \text{and}
\quad \psi = \tilde{\psi} \circ p . \label{phi-psi-til}
\een

In view of the general theory reviewed, e.g., in
Borodin and Salminen~\cite[Section~II.1]{BS},
the functions $\tilde{\varphi}$, $\tilde{\psi}$
are unique modulo multiplicative constants,
$C^1$ with absolutely continuous first derivatives,
and such that $\tilde{\varphi}$ (resp., $\tilde{\psi}$)
is strictly decreasing (resp., increasing).
Also, since $\liota$ (resp., $\riota$) is an absorbing
(resp., inaccessible) boundary point for $X$ if and
only if $p(\liota)$ (resp., $p(\riota)$) is an absorbing
(resp., inaccessible) boundary point for $\tilde{X}
\equiv p(X)$,
\begin{gather}
\text{if } \liota \text{ is absorbing for $X$, then }
\tilde{\varphi} \bigl( p(\liota) \bigr) := \lim
_{y \downarrow p(\liota)} \tilde{\varphi} (y) < \infty
\text{ and } \tilde{\psi} \bigl( p(\liota) \bigr) := \lim
_{y \downarrow p(\liota)} \tilde{\psi} (y) = 0 ,
\nonumber \\
\text{if } \riota \text{ is absorbing for $X$, then }
\tilde{\varphi} \bigl( p(\riota) \bigr) := \lim
_{y \uparrow p(\riota)} \tilde{\varphi} (y) = 0
\text{ and } \tilde{\psi} \bigl( p(\riota) \bigr) :=
\lim _{y \uparrow p(\riota)} \tilde{\psi} (y) < \infty
, \nonumber \\
\text{and, if } \liota \text{ (resp., } \riota \text{) is
inaccessible for $X$, then } \lim _{y \downarrow
p(\liota)} \tilde{\varphi} (y) = \infty \text{ (resp., }
\lim _{y \uparrow p(\riota)} \tilde{\psi} (y) = \infty
\text{)} , \nonumber
\end{gather}
Furthermore, $\tilde{\varphi}$ and $\tilde{\psi}$
satisfy the homogeneous ODE in $\tilde{g}$
\be
\frac{1}{2} \bigl( \sigma \circ p^{-1} \bigr)^2 (y) 
\, \bigl( p_-' \circ p^{-1} \bigr)^2 (y) \, \tilde{g}'' (y)
- \bigl( r \circ p^{-1} \bigr) (y) \, \tilde{g} (y) = 0 ,
\ee
Lebesgue-a.e.\ in $p(\icalo)$.
This fact and the absolute continuity of $p^{-1}$
imply that the ODE
\be
\frac{1}{2} \sigma^2 (x) (p_-')^2 (x) \tilde{g}''
\bigl( p(x) \bigr) - r(x) \tilde{g} \bigl( p(x) \bigr)
= 0
\ee
is satisfied Lebesgue-a.e.\ in $\icalo$ for
$\tilde{g}$ standing for either $\tilde{\varphi}$
or $\tilde{\psi}$.
Combining these observations with (\ref{phi-psi-til})
and Lemma~\ref{u-utilde}.(III), we obtain all of
the required results.
\mbox{}\hfill$\Box$

\begin{ex} \label{ex:levels} {\rm
Suppose that the measure $\nu$ is of the form
\ben
\nu (dz) = \frac{b(z)}{\sigma^2 (z)} \, dz + \sum _{j=1}^k
\beta_j \, \delta _{z_j} (dz) , \nonumber
\een
for some function $b: \icalo \rightarrow \bbr$ such that
\ben
\int _{\ubar{x}}^{\bar{x}} \frac{\bigl| b(u) \bigr|}{\sigma^2 (u)}
\, du < \infty \quad \text{for all } \liota < \ubar{x} < \bar{x}
< \riota , \nonumber
\een
some constants $\beta_1 , \ldots , \beta_k \in \mbox{}
]{-1}, 1[$ and some distinct points $z_1, \ldots , z_k \in
\icalo$, where $\delta _{z_j} (dz)$ is the Dirac probability
measure that assigns unit mass on $\{ z_j \}$.
Using the occupation times formula, we can see that,
in this case, $X$ satisfies the SDE
\ben
X_t = x + \int _0^t b(X_s) \, ds + \sum _{j=1}^k
\beta_j L_t^{z_j} (X) + \int _0^t \sigma (X_s) \, dW_s
, \quad x \in \icalo . \nonumber
\een
In view of (\ref{p'ODE}), the restriction of the scale
function $p$ to $\icalo \setminus \{ z_1 , \ldots, z_k \}$
has absolutely continuous first derivative $p'$
$(= p_-' = p_+')$ that satisfies the ODE
\ben
\frac{1}{2} \sigma^2 (x) p'' (x) + b(x) p' (x) = 0 ,
\label{p'ODEex}
\een
Lebesgue-a.e.\ in $\icalo \setminus \{ z_1 ,
\ldots, z_k \}$ (see also footnote~\ref{foot1} about
$p''$).
Furthermore, (\ref{p-jumps}) implies that
\ben
(1 + \beta_j) p_+' (z_j) = (1 - \beta_j) p_-' (z_j) ,
\quad \text{for } j = 1 , \ldots , k . \label{p'ODEBCex}
\een
Using these observations, we derive the
expression
\ben
p_+' (x) = \exp \left( - \int _{x_1}^x \frac{2 b(u)}
{\sigma^2 (u)} \, du \right) \prod _{j=1}^k
\left( \frac{1 - \beta_j} {1 + \beta_j} \right) ^{{\bf 1}
_{\{ x_1 \leq z_j \leq x \}}} , \quad \text{for }
x \geq x_1 , \label{p'ODEsolex}
\een
as well as a similar one for $x<x_1$,
where $x_1 \in \icalo$ is an arbitrary fixed point.
If we denote by $g$ either of the excessive
functions $\varphi$ or $\psi$ given by
(\ref{rphi}) and (\ref{rpsi}), then (\ref{p'ODEBCex})
and the (absolute) continuity of $g_-' / p_-'$
(see Lemma~\ref{lem:exess}) imply that
\begin{align}
\frac{g_-'}{p_-'} (z_j) = \lim _{y >z_j ,\; y \downarrow z_j}
& \frac{g_-'}{p_-'} (y) = \frac{g_+'}{p_+'} (z_j)
\nonumber \\
& \Rightarrow \quad
(1 + \beta_j) g_+' (z_j) = (1 - \beta_j) g_-' (z_j) ,
\quad \text{for } j = 1 , \ldots , k . \label{ODEexBC}
\end{align}
Furthermore, Lemma~\ref{lem:exess} and
(\ref{p'ODEex}) imply that $g$ satisfies the ODE
\ben
\frac{1}{2} \sigma^2 (x) g'' (x) + b(x) g' (x) -
r(x) g(x) = 0 , \label{ODEex}
\een
Lebesgue-a.e.\ in $\icalo \setminus \{ z_1 ,
\ldots, z_k \}$ (see also footnotes~\ref{foot1},
\ref{foot2} about $g''$).
} \mbox{}\hfill$\Box$ \end{ex}

\begin{ex} \label{ex:skewBM} {\rm
In the case of the skew Brownian motion, which is the
unique strong solution to the SDE (\ref{X-BM}), we can
see that (\ref{p'ODEsolex}) yields the expressions
\be
p(x) = \begin{cases} x , & \text{if } x<0 , \\
\frac{1-\beta}{1+\beta} x , & \text{if } x \geq 0 ,
\end{cases} \quad \text{and} \quad
p^{-1} (\tilde{x}) = \begin{cases} \tilde{x} , & \text{if }
\tilde{x} < 0 , \\ \frac{1+\beta}{1-\beta} \tilde{x} , &
\text{if } \tilde{x} \geq 0 , \end{cases} 
\ee
for the scale function $p$ and its inverse
$p^{-1}$.
Accordingly, the SDE (\ref{X-tilde}) takes the form
\be
d\tilde{X}_t  = \bigl( p_-' \circ p^{-1} \bigr) (\tilde{X}_t)
\, dW_t = \left( {\bf 1} _{]{-\infty}, 0]} (\tilde{X}_t)
+ \frac{1-\beta}{1+\beta} \, {\bf 1} _{]0, \infty[}
(\tilde{X}_t) \right) dW_t , \quad \tilde{X}_0
= p(x) \in \bbr ,
\ee
where $W$ is a standard Brownian motion.
This SDE has a unique strong solution and
the skew Brownian motion is the process
$X = p^{-1} (\tilde{X})$.
To verify that this process indeed satisfies the SDE
(\ref{X-BM}), we first use It\^{o}-Tanaka's formula
(e.g., see Assing and
Schmidt~\cite[Proposition~1.14]{AS98}) to obtain
\begin{align}
dX_t & = \frac{1}{2} \left( \frac{1+\beta}{1-\beta} - 1
\right) dL_t^{\tilde{X}, 0} + \bigl( p^{-1} \bigr)' (\tilde{X}_t)
\bigl( p_-' \circ p^{-1} \bigr) (\tilde{X}_t) \, dW_t
\nonumber \\
& = \frac{\beta}{1-\beta} \, dL_t^{\tilde{X}, 0} + dW_t ,
\nonumber
\end{align}
where $L^{\tilde{X}, 0}$ is the local time process
of $\tilde{X}$ at level 0.
In view of this result and the connection
\be
L^0 = \frac{1}{2} \bigl[ (p^{-1})_+' (0) + (p^{-1})_-' (0)
\bigr] L^{\tilde{X}, 0} = \frac{1}{1-\beta} L^{\tilde{X}, 0} .
\ee
of the local time $L^0$ of $X$ with the local time
$L^{\tilde{X}, 0}$ of $\tilde{X}$ (e.g., see Assing
and Schmidt~\cite[Lemma~1.18]{AS98} or Engelbert
and Schmidt~\cite[Proposition~4.29.iii]{ES91}),
we can see that $X$ satisfies the SDE (\ref{X-BM}).
} \mbox{}\hfill$\Box$ \end{ex}

\section{The solution to the general optimal stopping problem}
\label{sec:VI}

The value function of the optimal stopping problem
that aims at maximising the performance criterion
appearing in (\ref{criterion}) is defined by
\ben
v(x) = \sup _{(\bbs_x, \tau) \in \tcal_x} \bbe_x \left[ 
\exp \left( - \int _0^\tau r(X_s) \, ds \right) f(X_\tau)
{\bf 1} _{\{ \tau < \infty \}} \right] , \quad \text{for }
x \in \ical , \label{v}
\een
where the set of all stopping strategies $\tcal_x$ consists
of all pairs $(\bbs_x, \tau)$ such that $\bbs_x = \bigl(
\Omega, \fcal, (\fcal_t), \bbp_x, W, X \bigr)$ is a weak
solution to (\ref{X}) and $\tau$ is an associated
$(\fcal_t)$-stopping time.
We assume that the discounting rate function $r$
satisfies Assumption~\ref{A2}, while the reward
function $f$ satisfies the following assumption.

\begin{ass} \label{A3} {\rm
The positive function $f: \ical \rightarrow \bbr_+$
is Borel-measurable and its restriction to $\icalo$
is upper semicontinuous, namely,
\be
f(x) = \limsup _{y \rightarrow x} f(y) \quad \text{for all }
x \in \icalo .
\ee
} \vspace{-13.3mm}

\mbox{}\hfill$\Box$ \end{ass}
\smallskip

Our main result in this section establishes a complete
characterisation of the general optimal stopping problem
defined by (\ref{X}), (\ref{v}) in terms of solutions to the
variational inequality (\ref{v-VI}) in the sense of distributions,
which are introduced by the following definition.

\begin{de} \label{distr}
A function $v: \ical \rightarrow \R_+$ is a solution to the
variational inequality (\ref{v-VI}) if it satisfies the following
conditions:
\smallskip

\noindent {\rm (I)}
The restriction of $v$ to $\icalo$ is the difference of two
convex functions.
\smallskip

\noindent {\rm (II)}
The signed Radon measure on $\bigl( \icalo, \bcal (\icalo)
\bigr)$ defined by\,\footnote{See Lemma~\ref{u-utilde}
about the measure $(v_-' / p_-')'$.}
\be
\mu_v (dx) = - \frac{1}{2} \sigma^2 (x) p_-' (x) \left(
\frac{v_-'}{p_-'} \right) ' (dx) + r(x) v(x) \, dx
\ee
is positive and such that $\mu_v \bigl( \bigl\{ x \in \icalo
\mid \ v(x) > f(x) \bigr\} \bigr) = 0$.
\smallskip

\noindent {\rm (III)}
The inequality $v(x) \geq f(x)$ holds true for all
$x \in \icalo$.
\end{de}

\begin{thm} \label{VThm}
Consider the optimal stopping problem defined by
(\ref{X}) and (\ref{v}), and recall the minimal excessive
functions $\varphi$ and $\psi$  given by
(\ref{rphi}) and (\ref{rpsi}).
The following statements hold true:
\smallskip

\noindent {\rm (I)} 
If the problem data is such that
\ben
f(x) < \infty \text{ for all } x \in \ical , \quad 
\limsup _{y \downarrow \liota} \frac{f(y)} {\varphi (y)}
< \infty  \quad \text{and} \quad 
\limsup _{y \uparrow \riota} \frac{f(y)} {\psi (y)} < \infty 
, \label{opt0}
\een
then $v(x) < \infty$ for all $x \in \ical$. 
If any of the inequalities in (\ref{opt0}) fails, then 
$v(x) = \infty$ for all $x \in \ical$.
\smallskip

\noindent {\rm (II)} 
If the problem data is such that the inequalities in
(\ref{opt0}) all hold true, then the value function $v$
satisfies the variational inequality (\ref{v-VI}) 
in the sense of Definition~\ref{distr},
\begin{gather}
\lim_{y \in \icalo ,\; y \downarrow \liota} \frac{v(y)}
{\varphi (y)}  = \limsup _{y \downarrow \liota} \frac{f(y)}
{\varphi (y)} , \quad \lim _{y \in \icalo ,\; y \uparrow \riota}
\frac{v(y)} {\psi (y)}  = \limsup _{y \uparrow \riota}
\frac{f(y)} {\psi (y)} \label{BCv} \\
\text{and} \quad v(\liota) = f(\liota) \ \bigl(
\text{resp., } v(\riota) = f(\riota) \bigr) \quad  \text{if }
\liota \ (\text{resp., } \riota) \text{ is absorbing} .
\nonumber
\end{gather}

\noindent {\rm (III)} 
In the presence of (\ref{opt0}), if a positive function
$w: \ical \rightarrow \R_+$ satisfies the variational
inequality (\ref{v-VI}) in the sense of Definition~\ref{distr} 
as well as the growth conditions 
\begin{gather}
\lim _{y \in \icalo ,\; y \downarrow \liota} \frac{w(y)}
{\varphi (y)} = \limsup _{y \downarrow \liota} \frac{f(y)}
{\varphi (y)} , \quad \lim _{y \in \icalo ,\; y \uparrow \riota}
\frac{w(y)}{\psi(y)} = \limsup _{y \uparrow \riota}
\frac{f(y)} {\psi(y)} \label{BCw} \\
\quad \text{and} \quad  
w(\liota) = f(\liota) \ \bigl(\text{resp., } w(\riota) =
f(\riota) \bigr) \quad  \text{if } \liota \ (\text{resp., }
\riota) \text{ is absorbing} , \nonumber
\end{gather}
then $w(x) = v(x)$ for all $x \in \ical$. 
\smallskip

\noindent {\rm (IV)}
If
\begin{gather}
\limsup _{y \downarrow \liota} \frac{f(y)} {\varphi (y)}
= 0 \text{ if } \liota \text{ is inaccessible} , \quad
\limsup _{y \uparrow \riota} \frac{f(y)} {\psi (y)} = 0
\text{ if } \riota \text{ is inaccessible}  , \nonumber \\
f(\liota) = \limsup_{y \downarrow \liota} {f(y)}
\text{ if } \liota \text{ is absorbing}  \quad \text{and}
\quad f(\riota) = \limsup_{y \uparrow \riota} {f(y)}
\text{ if } \riota \text{ is absorbing} , \nonumber
\end{gather}
then the stopping strategy $(\bbs_x, \tau^\star)$,
where $\bbs_x$ is
a weak solution to (\ref{X}) and
\be
\tau_\star = \inf \bigl\{ t \geq 0 \mid \ v(X_t) = f(X_t)
\bigr\} ,
\ee
is optimal.
\end{thm}
{\bf Proof.}
In view of the results connecting the solvability
of (\ref{X}) with the solvability of (\ref{X-tilde})
that we have already used in the proof of
Lemma~\ref{lem:exess}, we can see that
\ben
v(x) = \sup _{(\bbs_x, \tau) \in \tcal_x} \bbe_x \left[
e^{- \int _0^\tau (r \circ p^{-1}) (\tilde{X}_s) \, ds} \tilde{f}
(\tilde{X}_\tau) \right] =: \tilde{v} \bigl( p(x) \bigr)
\quad \text{for all } x \in \ical , \label{vs-VT}
\een
where $\tilde{f} = f \circ p^{-1}$.
Also, we note that
\ben
\limsup _{y \downarrow \liota} \frac{f(y)}{\varphi (y)}
= \limsup _{y \downarrow \liota} \frac{\tilde{f} \bigl(
p(y) \bigr)}{\tilde{\varphi} \bigl( p(y) \bigr)}
\quad \text{and} \quad 
\limsup _{y \uparrow \riota} \frac{f(y)} {\psi (y)}
= \limsup _{y \uparrow \riota} \frac{\tilde{f} \bigl( p(y)
\bigr)}{\tilde{\psi} \bigl( p(y) \bigr)} ,
\label{phi-psi-lims-VT}
\een
where $\tilde{\varphi} = \varphi \circ p^{-1}$ and
$\tilde{\psi} = \psi \circ p^{-1}$ are the minimal
$(r \circ p^{-1})$-excessive functions of the
diffusion associated with the SDE (\ref{X-tilde}),
given by (\ref{phi-til}) and (\ref{psi-til}).

The identities in (\ref{phi-psi-lims-VT}) and
Theorem~6.3.(I) in Lamberton and
Zervos~\cite{LZ} imply (I).
If the inequalities in (\ref{opt0}) all hold true, then
Theorem~6.3 in Lamberton and Zervos~\cite{LZ}
asserts that the restriction of the value function
$\tilde{v}$ to $p(\icalo)$ is the difference of two convex
functions and satisfies the variational inequality
\be
\max \left\{ \frac{1}{2} \bigl( \sigma \circ p^{-1} \bigr)^2
(q) \bigl( p_-' \circ p^{-1} \bigr)^2 (q) \, \tilde{v}'' (dq) -
\bigl( r \circ p^{-1} \bigr) (q) \tilde{v} (q) \, dq , \
\tilde{f} (q) - \tilde{v} (q) \right\} = 0
\ee
in the sense that the Radon measure on $\bigl(
p(\icalo), \bcal (p(\icalo)) \bigr)$ defined by
\be
\mu _{\tilde{v}} (dq) = - \frac{1}{2} \bigl( \sigma \circ
p^{-1} \bigr)^2 (q) \bigl( p_-' \circ p^{-1} \bigr)^2 (q)
\, \tilde{v}'' (dq) + \bigl( r \circ p^{-1} \bigr) (q) \tilde{v}
(q) \, dq
\ee
is positive and such that $\mu _{\tilde{v}} \bigl(
\bigl\{ q \in p(\icalo) \mid \ \tilde{v} (q) > \tilde{f} (q)
\bigr\} \bigr) = 0$, while $\tilde{v} (q) \geq \tilde{f}
(q)$ for all $q \in p(\icalo)$.

In view of Lemma~\ref{u-utilde}.(I)-(II), $v =
\tilde{v} \circ p$ is the difference of two convex
functions and $\mu_v = \mu_{\tilde{v}} \circ p$
because $(v_-' / p_-') ' = \tilde{v}'' \circ p$.
Combining these observations with
(\ref{vs-VT})--(\ref{phi-psi-lims-VT}) and
Theorems~6.3,~6.4 in Lamberton and
Zervos~\cite{LZ}, we obtain all of the required
results in (II)-(IV).
\mbox{}\hfill$\Box$

\begin{rem} {\rm
It is worth stressing the precise nature of the
boundary conditions appearing in (\ref{BCv}) and
(\ref{BCw}).
The existence of the limits on the left-hand side of
(\ref{BCv}) is a result, while, the existence of the
limits on the left-hand side of (\ref{BCw}) is an
assumption.
Also, the limits on the left-hand sides of (\ref{BCv}),
(\ref{BCw}) are taken from inside the interior $\icalo$
of $\ical$.
On the other hand, the limsups on the right-hand
sides of (\ref{BCv}), (\ref{BCw}) are taken from
inside $\ical$ itself.
In view of these observations, we can see that,
e.g., if $\liota$ is absorbing, then we are faced
in (\ref{BCv}) with either the possibility that
\be
v(\liota) = f(\liota) = \lim _{y \in \icalo , \, y \downarrow \liota}
v(y) = \limsup _{y \downarrow \liota}  f(y) , \quad \text{if }
f(\liota) = \limsup _{y \downarrow \liota}  f(y) \geq \limsup
_{y \in \icalo , \, y \downarrow \liota} f(y) ,
\ee
or the possibility that
\be
v(\liota) = f(\liota) < \lim _{y \in \icalo , \, y \downarrow \liota}
v(y) = \limsup _{y \downarrow \liota} f(y) , \quad \text{if }
f(\liota) < \limsup _{y \downarrow \liota}  f(y) = \limsup
_{y \in \icalo , \, y \downarrow \liota} f(y) ,
\ee
where we have used the fact that, in this case,
$\varphi (\liota) := \lim _{x \downarrow \liota}
\varphi (x) < \infty$ (see Lemma~\ref{lem:exess}).
} \mbox{}\hfill$\Box$ \end{rem}

\begin{ex} \label{ex:VI} {\rm
Suppose that the measure $\nu$ is as in
Example~\ref{ex:levels}.
Given $\ubar{x} < \bar{x}$ such that $[\ubar{x}, \bar{x}[
\mbox{} \subseteq \icalo \setminus \{ z_1 , \ldots, z_k \}$,
we use the integration by parts formula and the fact
that the scale function $p$ has absolutely continuous
derivative satisfying (\ref{p'ODEex}) to calculate
\be
\frac{v_-'}{p'} (\bar{x}) - \frac{v_-'}{p'} (\ubar{x})
= \int _{[\ubar{x}, \bar{x}[} \left[ \frac{1}{p' (y)} v'' (dy)
+ \frac{2b(y)}{\sigma^2 (y) p' (y)} v_-' (y) \, dy \right] .
\ee
Also, we note that the measure $\mu_v$ defined in
Definition~\ref{distr}.(II) is such that
\begin{align}
\mu_v \bigl( \{ z_j \} \bigr) & = - \frac{1}{2} \sigma^2
(z_j) p_-' (z_j) \left[ \frac{v_+'}{p_+'} (z_j) -
\frac{v_-'}{p_-'} (z_j) \right] \nonumber \\
& \stackrel{(\ref{p-jumps})}{=} - \frac{\sigma^2 (z_j)}
{2 (1 - \beta_j)} \bigl[ (1+\beta_j)
v_+' (z_j) - (1-\beta_j) v_-' (z_j) \bigr] . \nonumber
\end{align}
It follows that, in this case, the variational inequality
(\ref{v-VI}) takes the form
\ben
\max \left\{ \frac{1}{2} \sigma^2 (x) v'' (dx) + b(x)
v_-' (x) \, dx - r(x) v(x) \, dx , \ f(x) - v(x) \right\} = 0
\label{VI-sc1}
\een
inside $\bigl( \icalo \setminus \{ z_1 , \ldots, z_k \} \, ,
\bcal (\icalo \setminus \{ z_1 , \ldots, z_k \}) \bigr)$,
coupled with the conditions
\ben
\max \Big\{ (1+\beta_j) v_+' (z_j) - (1-\beta_j) v_-'
(z_j) , \ f(z_j) - v(z_j) \Big\} = 0 , \quad \text{for }
j = 1, \ldots, k . \label{VI-sc-cond}
\een
Furthermore, if $v$ has absolutely continuous
first derivative, namely, if $v''(dx)$ is equal to
$v'' (x) \, dx$ (see also footnote~\ref{foot1}),
then $v$ should satisfy
\ben
\max \left\{ \frac{1}{2} \sigma^2 (x) v'' (x) + b(x)
v' (x) - r(x) v(x) , \ f(x) - v(x) \right\} = 0 ,
\label{VI-sc2}
\een
Lebesgue-a.e.\ in $\icalo \setminus \{ z_1 , \ldots,
z_k \}$ as well as the conditions (\ref{VI-sc-cond}).
} \mbox{}\hfill$\Box$ \end{ex}

\begin{rem} \label{rem:VIuse} {\rm
To appreciate how variational inequalities can be
used to systematically identify critical parts of the
state space $\ical$ that belong to the waiting region,
consider the previous example.
In this context, we can make the following observations:
\smallskip

\noindent (a)
Given $x \in \icalo$, if one of the limits
\be
f_+' (x) = \lim _{\varepsilon \downarrow 0}
\frac{f(x+\varepsilon) - f(x)}{\varepsilon}
\quad \text{and} \quad
f_-' (x) = \lim _{\varepsilon \downarrow 0}
\frac{f(x) - f(x-\varepsilon)}{\varepsilon} ,
\ee
does not exist, then $x$ belongs to the closure of
the waiting region.
If both $f_\pm' (x)$ fail to exist, then $x$ belongs
to the waiting region because it cannot belong to
the stopping region.
\smallskip

\noindent (b)
Given $j = 1, \ldots, k$, if both of the derivatives
$f_\pm' (z_j)$ exist and
\be
(1+\beta_j) f_+' (z_j) - (1-\beta_j) f_-' (z_j) > 0 ,
\ee
then (\ref{VI-sc-cond}) implies that $z_j$ belongs
to the waiting region.
In particular, if $f$ is $C^1$ at $z_j$, then $z_j$
belongs to the waiting region if $\beta_j \in
\mbox{} ]0, 1[$.
\smallskip

\noindent (c)
Suppose that the restriction of $f$ to an interval
$]\ubar{\iota} , \bar{\iota}[ \mbox{} \subseteq \ical$
is $C^2$.
The validity of (\ref{VI-sc1}) implies that
\be
\left\{ x \in \mbox{} ]\ubar{\iota} , \bar{\iota}[ \mbox{}
\ \Big| \ \ \frac{1}{2} \sigma^2 (x) f'' (x) + b(x)
f' (x) - r(x) f(x) > 0 \right\}
\ee
is a subset of the waiting region.
On the other hand, the intersection of the stopping
region with $]\ubar{\iota} , \bar{\iota}[$ is a
(usually strict) subset of the complement of this set.
} \mbox{}\hfill$\Box$ \end{rem}

\begin{rem} \label{rem:SF} {\rm
In the context of Example~\ref{ex:VI}, we can
make the following observations relative to the
so-called ``principle of smooth fit'':
\smallskip

\noindent (a)
If $f$ is $C^1$ then (\ref{VI-sc1}) implies that
the restriction of the value function $v$ to
$\icalo \setminus \{ z_1 , \ldots, z_k \}$ is $C^1$
(see Lamberton and Zervos~\cite[Corollary 7.5]{LZ}).
\smallskip

\noindent (b)
Given $j = 1, \ldots, k$, if $z_j$ belongs to the
stopping region, namely, $f(z_j) = v(z_j)$, then
(\ref{p'ODEBCex}) and (\ref{VI-sc-cond})
imply that
\ben
\frac{v_+' (z_j)}{p_+' (z_j)} - \frac{v_-' (z_j)}{p_-' (z_j)}
= \frac{1}{(1-\beta_j) p_-' (z_j)} \bigl[ (1+\beta_j)
v_+' (z_j) - (1-\beta_j) v_-' (z_j) \bigr] \leq 0 .
\label{SF-calc}
\een
In general, this inequality can be strict: see
Remark~\ref{rem:SFex} in Section~\ref{sec:GBMsol}.

\noindent (c)
Given $j = 1, \ldots, k$, if $z_j$ belongs to the
waiting region, namely, $f(z_j) < v(z_j)$,
then (\ref{p'ODEBCex}) and (\ref{VI-sc-cond})
imply that
\be
\frac{v_+' (z_j)}{p_+' (z_j)} - \frac{v_-' (z_j)}{p_-' (z_j)}
= 0 .
\ee
} \mbox{}\hfill$\Box$ \end{rem}


\section{The minimal excessive functions
of a skew geometric Brownian motion}
\label{sec:excessive}

From this point onwards, we fix a filtered probability space
$\bigl( \Omega, \fcal, (\fcal_t), \bbp \bigr)$ satisfying the
usual conditions and supporting a standard one-dimensional
$(\fcal_t)$-Brownian motion $W$.
In such a setting, we denote by $X$ the unique non-explosive
strong solution to the SDE (\ref{XGBM}).

The conditions (\ref{ODEexBC}) in Example~\ref{ex:levels}
reduce to
\ben
(1+\beta) g_+' (z) = (1-\beta) g_-' (z) , \label{g'z}
\een
while, given a constant $r>0$, the ODE (\ref{ODEex}) in
Example~\ref{ex:levels} reduces to the Euler ODE
\ben
\half \sigma^2 x^2 g''(x) + bx g' (x) - rg(x) = 0 .
\label{ODE}
\een
It is well-known that every solution to (\ref{ODE})
is given by
\be
g(x) = A x^n + B x^m ,
\ee
for some constants $A, B \in \bbr$, where $m<0<n$
are the solutions to the quadratic equation
\be
\half \sigma^2 k^2 + \left( b - \half \sigma ^2 \right) k - r
= 0 ,
\ee
given by
\be
m, n = \frac{- \left( b - \half \sigma^2 \right) \mp
\sqrt{\left( b - \half \sigma^2 \right) ^2 + 2 \sigma^2 r}}
{\sigma^2} .
\ee
It is straightforward to verify that
\begin{gather}
r > b \ \Leftrightarrow \ n > 1 , \quad
n+m-1 = - \frac{2b}{\sigma^2} , \quad
nm = - \frac{2r}{\sigma^2} , \label{nm} \\
r-bm = \half \sigma^2 m (m-1) > 0
\quad \text{and} \quad
r-bn = \half \sigma^2 n (n-1) > 0 . \label{rbnm}
\end{gather}

The excessive functions $\psi = \psi (\cdot ; z)$ and
$\varphi = \varphi (\cdot ; z)$ that satisfy the ODE
(\ref{ODE}) inside $]0,z[ \mbox{} \cup \mbox{} ]z,\infty[$
as well as the condition (\ref{g'z}) are given by
\ben
\psi (x;z) = \begin{cases} x^n , & \text{if } x<z , \\
Ax^n + B(z) x^m , & \text{if } x \geq z , \end{cases}
\quad \text{and} \quad 
\varphi (x;z) = \begin{cases} C(z) x^n + D x^m , &
\text{if } x<z , \\ x^m , & \text{if } x \geq z , \end{cases}
\label{psi}
\een
for
\begin{gather}
A = \frac{n (1 - \beta) - m (1 + \beta)}{(n-m)(1 + \beta)}
\begin{cases} >1, & \text{if } \beta \in \mbox{} ]{-1}, 0[
, \\  \in \mbox{} ]0,1[ , & \text{if } \beta \in \mbox{} ]0, 1[
, \end{cases} \label{A} \\
B(z) = \frac{2n\beta z^{n-m}}{(n-m)(1 + \beta)} 
\begin{cases} < 0, & \text{if } \beta \in \mbox{}
]{-1}, 0[ , \\ > 0 , & \text{if } \beta \in \mbox{} ]0, 1[ ,
\end{cases} \label{B} \\
C(z) = \frac{2n\beta z^{n-m}}{(n-m)(1 - \beta)} 
\begin{cases} < 0 , & \text{if } \beta \in \mbox{}
]{-1}, 0[ , \\ > 0 , & \text{if } \beta \in \mbox{} ]0, 1[ ,
\end{cases} \nonumber \\
\text{and} \quad
D = \frac{n (1 - \beta) - m (1 + \beta)}{(n-m)(1 - \beta)}
\begin{cases} \in \mbox{} ]0,1[ , & \text{if } \beta
\in \mbox{} ]{-1}, 0[ , \\ > 1, & \text{if } \beta \in \mbox{}
]0, 1[ . \end{cases} \nonumber
\end{gather}
It is straightforward to verify that
\ben
\psi' (z-;z) \equiv nz^{n-1} < nAz^{n-1} + mB(z) z^{m-1}
\equiv \psi' (z+;z) \quad \Leftrightarrow \quad \beta < 0 .
\label{psi-conv1}
\een
Here, as well as in the rest of the paper, we adopt
the notation
$\psi' (x;z) = \frac{\partial \psi}{\partial x} (x;z)$ and
$\psi'' (x;z) = \frac{\partial^2 \psi}{\partial x^2} (x;z)$.

In the rest of the paper, we make the following assumption,
which is sufficient for the value function of the optimal
stopping problem defined by (\ref{XGBM}) and (\ref{vEx})
to be real-valued.

\begin{ass} \label{assm}
$r>b \stackrel{(\ref{nm})}{\Leftrightarrow} n>1$.
\end{ass}
Indeed, if $r<b \Leftrightarrow n<1$, then
Theorem~\ref{VThm}.(I) implies that the value
function given by (\ref{vEx}) is identically equal
to $\infty$.

In the presence of Assumption~\ref{assm}, we can verify
that
\begin{gather}
\frac{r}{r-b} = \frac{nm}{(n-1)(m-1)} < \frac{n}{n-1}
\label{rb-nm} \\
\text{and} \quad
\beta_{\mathrm c} := \frac{n-1}{n+2m-1} \in
\begin{cases}
]1, \infty[ , & \text{if } n+2m-1 > 0 , \\
]{-\infty}, {-1}] , & \text{if } n+2m-1 < 0 \text{ and } b \leq 0 , \\
]{-1}, 0[ , & \text{if } n+2m-1 < 0
\text{ and } b > 0 . \end{cases} \label{beta-c}
\end{gather}
Here, deriving the possible values of $\beta_{\mathrm c}$
involves the observation that, if $n+2m-1 < 0$, then
\ben
\beta_{\mathrm c} \in \mbox{} ]{-1}, 0[ \mbox{}
\quad \Leftrightarrow \quad n+m-1< 0
\quad \stackrel{(\ref{nm})}{\Leftrightarrow} \quad
b > 0 . \label{n+m-1/b}
\een
Combining the range of values of the point
$\beta_{\mathrm c}$ given by (\ref{beta-c}), with the
observation that
\be
(n-1) (1-\beta) - 2m \beta < 0
\quad \Leftrightarrow \quad
\beta \begin{cases} > \beta_{\mathrm c} ,
& \text{if } n+2m-1 > 0 , \\ < \beta_{\mathrm c} ,
& \text{if } n+2m-1 < 0 , \end{cases}
\ee
we can see that
\ben
\text{given any $\beta \in \mbox{} ]{-1}, 1[$,} \quad
(n-1) (1-\beta) - 2m \beta < 0 \ \Leftrightarrow \
\Bigl( b > 0 \text{ and } \beta \in
\mbox{} ]{-1}, \beta_{\mathrm c}[ \Bigr) .
\label{convexity}
\een
Furthermore, 
\ben
\text{given any $\beta \in \mbox{} ]{-1}, 0[$,} \quad
(n-1)(1-\beta) - 2m\beta < 0 \ \Leftrightarrow \
\frac{n}{n - \frac{1+\beta}{1-\beta}} < \frac{r}{r-b} .
\label{convexity*}
\een

In the following result, we concentrate on the increasing
function $\psi$ because only this is involved in the
solution to the optimal stopping problem we consider
in this main section.
In particular, the critical points defined by
\ben
\gzc = \frac{rK}{r-b} , \quad
\gzb = \frac{nK}{n - \frac{1+\beta}{1-\beta}} ,
\text{ if } n \neq \frac{1+\beta}{1-\beta} ,
\quad \text{and} \quad
\gzo = \frac{nK}{n-1} \label{cp}
\een
play a critical role in differentiating the different
qualitative forms of the optimal strategy.

\begin{lem} \label{lem:psi}
Suppose that Assumption~\ref{assm} holds true.
The function $\psi (\cdot ; z)$ defined by (\ref{psi})--(\ref{B})
is such that the following statements hold true:
\smallskip

\noindent
{\rm (I)}
If $b \leq 0$ and $\beta \in \mbox{} ]{-1}, 0[$, then
$\psi (\cdot ; z)$ is convex.
\smallskip

\noindent
{\rm (II)}
If $b > 0$ and $\beta \in [\beta_{\mathrm c}, 0[$,
where $\beta_{\mathrm c} \in \mbox{} ]{-1}, 0[$ is defined
by (\ref{beta-c}), then $\psi (\cdot ; z)$ is convex.
\smallskip

\noindent
{\rm (III)}
If $b > 0$ and $\beta \in \mbox{} ]{-1}, \beta_{\mathrm c}[$,
where $\beta_{\mathrm c} \in \mbox{} ]{-1}, 0[$ is defined
by (\ref{beta-c}), then the restrictions of $\psi (\cdot ; z)$
to $[0,z]$ as well as to $[\gc^{-1} z , \infty[$ are convex,
while the restriction of $\psi (\cdot ; z)$ to $[z, \gc^{-1} z]$
is concave, where
\ben
\gc = \left( - \frac{(n-1) \bigl[ n(1 - \beta) - m(1 + \beta) \bigr]}
{2m(m-1)\beta} \right)^{\frac{1}{n-m}} \in \mbox{} ]0,1[ .
\label{gc}
\een

\noindent
{\rm (IV)}
If $\beta \in \mbox{} ]0,1[$, then  the restrictions of
$\psi (\cdot ; z)$ to $[0,z]$ as well as to $[z,\infty[$
are convex but $\psi (\cdot ; z)$ is not convex in its
entire domain.
\smallskip

\noindent
Cases~(I) and~(II) are illustrated by Figure~1,
while Cases~(III) and~(IV) are illustrated by Figures~2
and~3, respectively.
Furthermore, the critical points $\gzc$, $\gzb$, $\gzo$
defined by (\ref{cp}) are such that,
\begin{align}
\text{in Cases~(I), (II)} , & \quad (n-1)(1-\beta) - 2m
\beta \geq 0 \text{ \,and \ } 0 < \gzc \leq \gzb < \gzo ,
\label{cp-I-II} \\
\text{in Case~(III)} , & \quad (n-1)(1-\beta) - 2m
\beta < 0 \text{ \,and \ } 0 < \gzb < \gzc < \gzo ,
\label{cp-III} \\
\text{and in Case~(IV)} , & \quad (n-1)(1-\beta) - 2m
\beta < 0 , \quad 0 < \gzc < \gzo , \nonumber \\
& \quad \frac{1+\beta}{1-\beta} < n \ \Rightarrow \
0 < \gzo < \gzb \text{ \,and \,} n < \frac{1+\beta}{1-\beta}
\ \Rightarrow \ \gzb < 0 < \gzo \label{cp-IV} ,
\end{align}
with equalities (resp., strict inequalities) in place of
weak inequalities in (\ref{cp-I-II}) if $b > 0$ and $\beta
= \beta_{\mathrm c}$ (resp., otherwise).
\end{lem}

\begin{picture}(160,110)
\put(20,15){\begin{picture}(120,95)

\put(0,0){\line(1,0){120}}
\put(120,0){\line(0,1){95}}
\put(0,0){\line(0,1){95}}
\put(0,95){\line(1,0){120}}

\put(10,10){\vector(1,0){100}}
\put(10,10){\vector(0,1){80}}
\put(12,87){$\psi (x)$}
\put(108,6){$x$}

\color{blue}
\put(0,-0.2){\qbezier(10,10)(15,10.5)(25,12)}
\put(0,-0.1){\qbezier(10,10)(15,10.5)(25,12)}
\put(0,0){\qbezier(10,10)(15,10.5)(25,12)}
\put(0,0.1){\qbezier(10,10)(15,10.5)(25,12)}
\put(0,0.1){\qbezier(10,10)(15,10.5)(25,12)}

\put(0,-0.2){\qbezier(25,12)(40,17.5)(50,25)}
\put(0,-0.1){\qbezier(25,12)(40,17.5)(50,25)}
\put(0,0){\qbezier(25,12)(40,17.5)(50,25)}
\put(0,0.1){\qbezier(25,12)(40,17.5)(50,25)}
\put(0,0.2){\qbezier(25,12)(40,17.5)(50,25)}

\put(0,-0.2){\qbezier(50,25)(70,37)(100,85)}
\put(0,-0.1){\qbezier(50,25)(70,37)(100,85)}
\put(0,0){\qbezier(50,25)(70,37)(100,85)}
\put(0,0.1){\qbezier(50,25)(70,37)(100,85)}
\put(0,0.2){\qbezier(50,25)(70,37)(100,85)}

\color{black}
\put(25,9){\line(0,1){2}}
\put(0,0){\qbezier[4](25,10)(25,11)(25,12)}
\put(23.5,5.5){$z$}

\end{picture}}

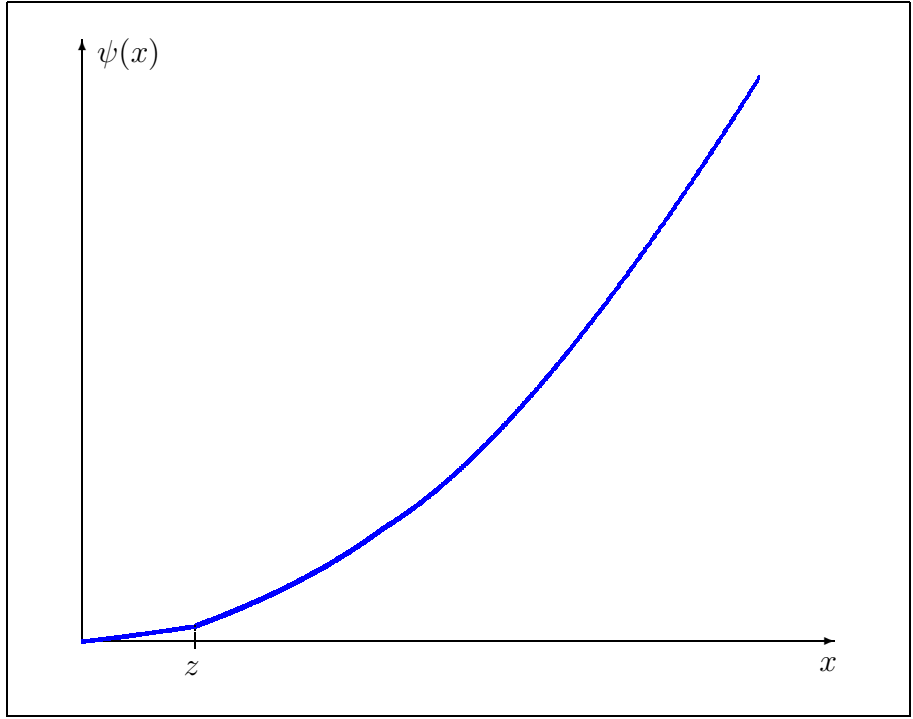
\put(20,10){\small{{\bf Figure 1.} Graph of the
function $\psi$ in Cases~(I) and~(II) of Lemma~\ref{lem:psi}.}}
\end{picture}

\begin{picture}(160,110)
\put(20,15){\begin{picture}(120,95)

\put(0,0){\line(1,0){120}}
\put(120,0){\line(0,1){95}}
\put(0,0){\line(0,1){95}}
\put(0,95){\line(1,0){120}}

\put(10,10){\vector(1,0){100}}
\put(10,10){\vector(0,1){80}}
\put(12,87){$\psi (x)$}
\put(108,6){$x$}

\color{blue}
\put(0,-0.2){\qbezier(10,10)(30,15)(50,25)}
\put(0,-0.1){\qbezier(10,10)(30,15)(50,25)}
\put(0,0){\qbezier(10,10)(30,15)(50,25)}
\put(0,0.1){\qbezier(10,10)(30,15)(50,25)}
\put(0,0.1){\qbezier(10,10)(30,15)(50,25)}

\put(0,-0.2){\qbezier(50,25)(55,35)(80,47.5)}
\put(0,-0.1){\qbezier(50,25)(55,35)(80,47.5)}
\put(0,0){\qbezier(50,25)(55,35)(80,47.5)}
\put(0,0.1){\qbezier(50,25)(55,35)(80,47.5)}
\put(0,0.2){\qbezier(50,25)(55,35)(80,47.5)}

\put(0,-0.2){\qbezier(80,47.5)(92,55)(100,85)}
\put(0,-0.1){\qbezier(80,47.5)(92,55)(100,85)}
\put(0,0){\qbezier(80,47.5)(92,55)(100,85)}
\put(0,0.1){\qbezier(80,47.5)(92,55)(100,85)}
\put(0,0.2){\qbezier(80,47.5)(92,55)(100,85)}

\color{black}
\put(50,9){\line(0,1){2}}
\put(0,0){\qbezier[20](50,10)(50,17.5)(50,25)}
\put(49,5.5){$z$}

\put(80,9){\line(0,1){2}}
\put(0,0){\qbezier[40](80,10)(80,28.5)(80,47.5)}
\put(75,4.8){$\gc^{-1} z$}

\end{picture}}

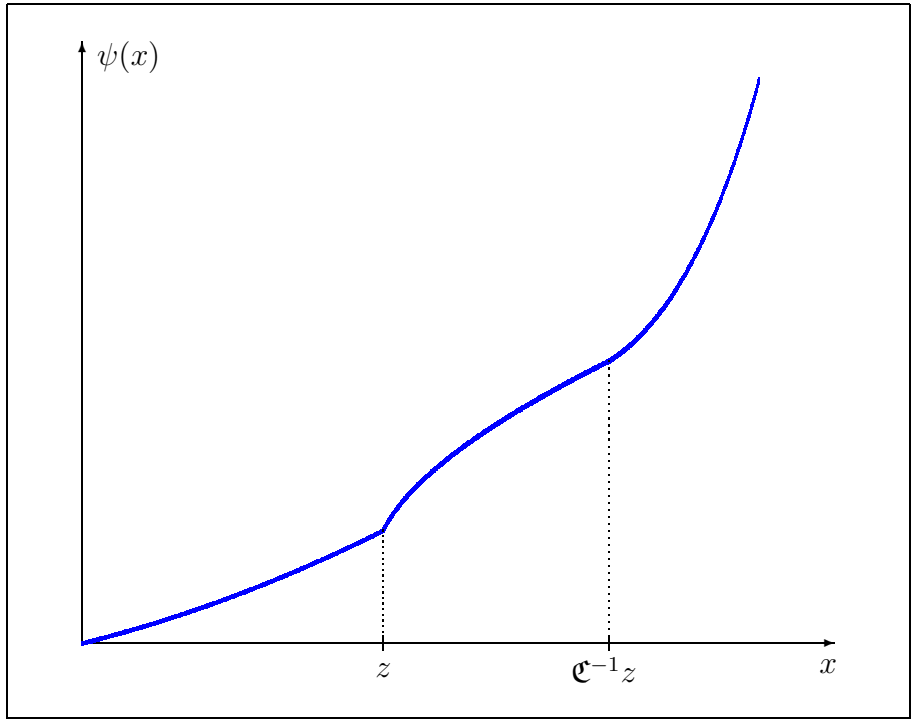
\put(20,10){\small{{\bf Figure 2.} Graph of the
function $\psi$ in Case~(III) of Lemma~\ref{lem:psi}.}}
\end{picture}

\normalsize
\begin{picture}(160,110)
\put(20,15){\begin{picture}(120,95) 

\put(0,0){\line(1,0){120}}
\put(120,0){\line(0,1){95}}
\put(0,0){\line(0,1){95}}
\put(0,95){\line(1,0){120}}

\put(10,10){\vector(1,0){100}}
\put(10,10){\vector(0,1){80}}
\put(12,87){$\psi (x)$}
\put(108,6){$x$}

\color{blue}
\put(0,-0.2){\qbezier(10,10)(30,11)(40,25)}
\put(0,-0.1){\qbezier(10,10)(30,11)(40,25)}
\put(0,0){\qbezier(10,10)(30,11)(40,25)}
\put(0,0.1){\qbezier(10,10)(30,11)(40,25)}
\put(0,0.1){\qbezier(10,10)(30,11)(40,25)}

\put(0,-0.2){\qbezier(40,25)(53,27)(60,32.5)}
\put(0,-0.1){\qbezier(40,25)(53,27)(60,32.5)}
\put(0,0){\qbezier(40,25)(53,27)(60,32.5)}
\put(0,0.1){\qbezier(40,25)(53,27)(60,32.5)}
\put(0,0.2){\qbezier(40,25)(53,27)(60,32.5)}

\put(0,-0.2){\qbezier(60,32.5)(80,48)(100,85)}
\put(0,-0.1){\qbezier(60,32.5)(80,48)(100,85)}
\put(0,0){\qbezier(60,32.5)(80,48)(100,85)}
\put(0,0.1){\qbezier(60,32.5)(80,48)(100,85)}
\put(0,0.2){\qbezier(60,32.5)(80,48)(100,85)}

\color{black}
\put(40,9){\line(0,1){2}}
\put(0,0){\qbezier[22](40,10)(40,17.5)(40,25)}
\put(39,6){$z$}

\end{picture}}

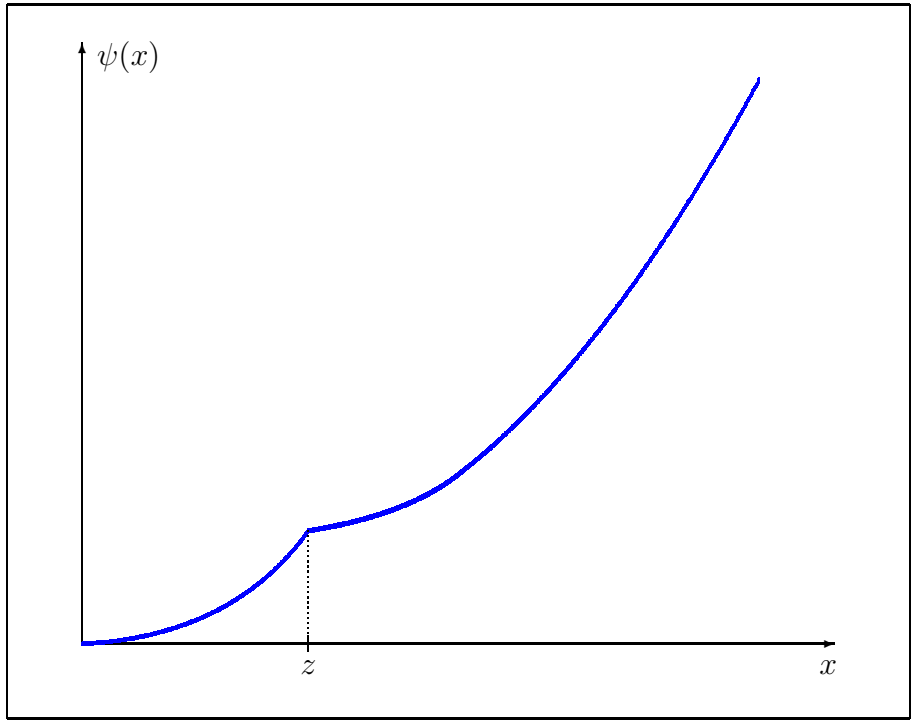
\put(20,10){\small{{\bf Figure 3.} Graph of the
function $\psi$ in Case~(IV) of Lemma~\ref{lem:psi}.}}
\end{picture}

\noindent
{\bf Proof.}
We first note that $\psi (\cdot ; z)$ is always convex
in $[0,z]$.
On the other hand, the inequality $n>1$ and the calculation
\be
\psi '' (x;z) & = \bigl[ n(n-1) Ax^{n-m} + m(m-1) B(z) \bigr]
x^{m-2} , \quad \text{for } x>z ,
\ee
imply the equivalences
\begin{align}
\psi '' (x;z) > 0 \text{ for all } x>z \quad 
& \Leftrightarrow \quad \psi '' (z+;z) \equiv
\frac{n \bigl[ (n-1) (1-\beta) - 2m\beta \bigr] z^{n-2}}
{1+\beta} \geq 0 \nonumber \\
& \Leftrightarrow \quad (n-1) (1-\beta) - 2m \beta
\geq 0 . \nonumber
\end{align}
Combining these observations with
(\ref{psi-conv1})--(\ref{cp}), we obtain the required results.
$\mbox{}$\hfill $\Box$

\section{Preliminary analytic results for the solution to
the optimal stopping problem defined by (\ref{XGBM})
and (\ref{vEx})}
\label{sec:F}

We now establish a pair of technical analytic results that
we will need for the solution of the optimal stopping problem
defined by (\ref{XGBM}) and (\ref{vEx}), which we
derive in the next section.
(This section can easily be skipped at a first reading.)
Given any $z>0$ fixed, we consider the equation
\ben
F(x;z) = 0 , \label{Feqn}
\een
for $x>z$, where $F$ is the function defined by
\ben
F(x;z) = \bigl[ (n-1)x - nK \bigr] A x^{n-m} + \bigl[ (m-1)x
- mK \bigr] B(z) , \quad \text{for } x>0 , \label{F}
\een
which admits the expression
\ben
F(x;z) = \bigl[ (x-K) \psi' (x;z) - \psi (x;z) \bigr] x^{-m+1}
, \quad \text{for } x>z . \label{F-psi}
\een

The following result involves the critical points $\gzc$,
$\gzb$, $\gzo$ defined by (\ref{cp}) and is structured
based on the four cases of Lemma~\ref{lem:psi}.

\begin{lem} \label{lem:a(z)}
In the presence of Assumption~\ref{assm}, the
following statements hold true:
\smallskip

\noindent {\rm (i)}
If the problem's parameters are as in Cases~(I)
or~(II) of Lemma~\ref{lem:psi}, then equation (\ref{Feqn})
defines uniquely a strictly decreasing $C^1$ function
$\upalpha : \mbox{} ]0, \gzb[ \mbox{} \rightarrow \mbox{}
]\gzb, \gzo[$ such that
\begin{gather}
\lim _{z \rightarrow 0} \upalpha (z) = \gzo , \quad
\lim _{z \rightarrow \gzb} \upalpha (z) = \gzb ,
\label{gotha-i1} \\
F(x; z) \begin{cases} 
< 0 & \text{for all } x \in [z \vee K , \upalpha (z)[ , \\ 
> 0 & \text{for all } x > \upalpha (z) , \end{cases}
\quad \text{and} \quad
\frac{\partial F}{\partial x} (x; z) > 0
\text{ for all } x \geq \upalpha (z) . \label{gotha-i2}
\end{gather}

\noindent {\rm (ii)}
If the problem's parameters are as in Case~(III)
of Lemma~\ref{lem:psi}, then $K < \gc \gzc$, where
$\gc \in \mbox{} ]0,1[$ is given by (\ref{gc}), and
equation (\ref{Feqn}) defines uniquely a strictly
decreasing $C^1$ function $\upalpha : \mbox{}
]0, \gc \gzc[ \mbox{} \rightarrow \mbox{} ]\gzc, \gzo[$
such that
\begin{gather}
\lim _{z \rightarrow 0} \upalpha (z) = \gzo , \quad
\lim _{z \rightarrow \gc \gzc} \upalpha (z) = \gzc ,
\label{gotha-ii1} \\
F(x; z) \begin{cases} 
< 0 & \text{for all } x \in [\gzc , \upalpha (z)[ , \\ 
> 0 & \text{for all } x > \upalpha (z) , \end{cases}
\quad \text{and} \quad
\frac{\partial F}{\partial x} (x; z) > 0
\text{ for all } x \geq \upalpha (z) . \label{gotha-ii2}
\end{gather}
Furthermore,
\ben
F(x;z) > 0 \quad \text{for all } z \in [\gc \gzc, \gzc[
\mbox{} \text{ and } x > \gzc . \label{F>0-ii}
\een

\noindent {\rm (iii)}
If $\beta \in \mbox{} ]0,1[$ (Case~(IV) of
Lemma~\ref{lem:psi}), then equation (\ref{Feqn})
defines uniquely a $C^1$ function
$\upalpha : \mbox{} ]0, \infty[ \mbox{} \rightarrow \mbox{}
]\gzo, \infty[$ such that
\begin{gather}
\upalpha (z) \in \begin{cases} 
]z \vee \gzo , \infty[ \text{ \ for all } z \in \mbox{}
]0, \infty[ , & \text{if } n \leq \frac{1+\beta}{1-\beta} , \\
]z \vee \gzo , \gzb[ \text{ \ for all } z \in \mbox{}
]0, \gzb[ , & \text{if } n > \frac{1+\beta}{1-\beta} , \\
]\gzb , z[ \text{ \ for all } z \in \mbox{}
]\gzb, \infty[ , & \text{if } n > \frac{1+\beta}{1-\beta} ,
\end{cases} \label{gotha-iii1} \\
F(x; z) \begin{cases} 
< 0 & \text{for all } x \in \mbox{} ]z \wedge \gzo,
\upalpha (z)[ , \\  > 0 & \text{for all } x > \upalpha (z) ,
\end{cases}
\quad \text{and} \quad
\frac{\partial F}{\partial x} (x; z) > 0
\text{ for all } x \geq \upalpha (z) . \label{gotha-iii2}
\end{gather}
\end{lem}
{\bf Proof.}
Throughout the proof, we use repeatedly the expressions
and signs of $A$, $B$, given by (\ref{A}), (\ref{B}), as well
as the results in (\ref{nm}), (\ref{rbnm}) and (\ref{rb-nm})
without special mention.
We first note that
\begin{gather}
\lim _{x \rightarrow 0} F(x;z) = -mK B(z)
\begin{cases} < 0 , & \text{if } \beta < 0 , \\
> 0 , & \text{if } \beta > 0 , \end{cases}
\label{Flims+1} \\
F (\gzo ; z) = - \frac{(n-m)K}{n-1} B(z) 
\begin{cases} > 0 , & \text{if } \beta < 0 , \\
< 0 , & \text{if } \beta > 0 , \end{cases}
\quad \text{and} \quad
\lim _{x \rightarrow \infty} F(x;z) = \infty .
\label{Flims+2}
\end{gather}
We also calculate
\begin{align}
\frac{\partial F}{\partial x} (x;z) = \mbox{} & (n-1)
\Bigr[ (n-m+1) x - (n-m) \gzo \Bigr] A x^{n-m-1}
+ (m-1) B(z) , \label{Fa()} \\
\frac{\partial^2 F}{\partial x^2} (x;z) = \mbox{} &
(n-1) (n-m) \Bigl[ (n-m+1) x - (n-m-1)\gzo \Bigr]
A x^{n-m-2} , \label{Faa()} \\
\frac{\partial F}{\partial z} (x; z) =
\mbox{} & \frac{2n \bigl[ (m-1) x - mK \bigr] \beta}
{1+\beta} z^{n-m-1} , \label{Fz()} \\
\text{and} \quad 
F(z;z) = \mbox{} & \frac{1-\beta}{1+\beta} \left[ \left(
n - \frac{1+\beta} {1-\beta} \right) z - nK \right] z^{n-m} .
\label{F(z,z)}
\end{align}
The calculation (\ref{Faa()}) implies that
\ben
\frac{\partial F}{\partial x} (\cdot; z) \text{ is strictly}
\begin{cases}
\text{decreasing in } ]0, x_\dagger[ , \\
\text{increasing in } ]x_\dagger, \infty[ ,
\end{cases}  \text{ where } x_\dagger =
\frac{(n-m-1) \gzo}{n-m+1} \in \mbox{} ]0, \gzo[ .
\label{xdag}
\een
Combining this observation with the limits
\be
\lim _{x \rightarrow 0} \frac{\partial F}{\partial x} (x;z)
= (m-1) B(z) \begin{cases} > 0 , & \text{if } \beta < 0, \\
< 0 , & \text{if } \beta > 0 , \end{cases}
\quad \text{and} \quad
\lim _{x \rightarrow \infty} \frac{\partial F}{\partial x} (x;z)
= \infty ,
\ee
which follow from (\ref{Fa()}), we can see that
\be
\text{if } \beta < 0 \text{ and } \frac{\partial F}{\partial x}
(x_\dagger; z) \geq 0 , \quad \text{then }
\frac{\partial F}{\partial x} (x; z) \geq 0 \text{ for all }
x > 0 ,
\ee
or there exist strictly positive constants
$\underline{x}_\dagger (z) < x_\dagger <
\overline{x}_\dagger (z)$ such that
\begin{align}
\text{if } \beta < 0 \text{ and } \frac{\partial F}{\partial x}
(x_\dagger; z) < 0 , & \nonumber \\
\text{then } \frac{\partial F}{\partial x} (x; z) & \begin{cases}
> 0 , & \text{for all } x \in \mbox{} ]0, \underline{x}_\dagger
(z)[ \mbox{} \cup \mbox{} ]\overline{x}_\dagger (z) , \infty[ ,
\\ < 0 , & \text{for all } x \in \mbox{} ]\underline{x}_\dagger
(z) , \overline{x}_\dagger (z)[ , \end{cases}  \label{Fsh2}
\end{align}
or there exists a constant $x_\ddagger (z) >
x_\dagger$ such that
\ben
\text{if } \beta > 0 , \quad \text{then }
\frac{\partial F}{\partial x} (x; z)
\begin{cases}
< 0 , & \text{for all } x \in \mbox{} ]0, x_\ddagger (z)[ , \\  
> 0 , & \text{for all } x \in \mbox{} ]x_\ddagger (z), \infty[ .
\end{cases}  \label{Fsh3}
\een

Keeping in mind that $n>1$, and $\gzb > 0$ if and
only if $n > \frac{1+\beta}{1-\beta}$, we can see that
\begin{align}
F (\gzb ; z) & = \frac{K}{\left( n - \frac{1+\beta}{1-\beta}
\right) (1-\beta)} \Bigl[ 2n\beta A \gzb ^{n-m}
- \bigl[ n(1-\beta) - m(1+\beta) \bigr] B(z) \Bigr]
\nonumber \\
& = \frac{2\bigl[ n(1-\beta) - m(1+\beta) \bigr] \beta}
{(n-m) (1-\beta) (1+\beta)} \gzb \bigl[ \gzb ^{n-m}
- z^{n-m} \bigr] \nonumber \\
& \begin{cases}
< 0 , & \text{if } \left( \beta < 0 \text{ and } z < \gzb \right)
\text{ or } \left( \beta > 0 \text{ and } z > \gzb \right)
, \\ > 0 , & \text{if } \left( \beta < 0 \text{ and } z > \gzb \right)
\text{ or } \left( \beta > 0 \text{ and } z \in \mbox{} ]0, \gzb[
\right) . \end{cases} \label{F-gzb}
\end{align}
Furthermore, we can use the definition (\ref{F}) of $F$
and (\ref{Fa()}) to see that, given any $\beta \in \mbox{}
]{-1}, 0[$,
\begin{align}
\frac{\partial F}{\partial x} (K; z) & = - (m-1) K^{-1} F(K; z)
\nonumber \\
& =  \frac{(m-1) K^{n-m}}{(n-m) (1+\beta)} \left[
n(1-\beta) - m(1+\beta) + 2n \beta \left( \frac{z}{K}
\right) ^{n-m} \right] \nonumber \\
& < 0 \qquad \text{for all } z < \tilde{z} , \label{F_x(K,z)}
\end{align}
where
\be
\tilde{z} = \left( - \frac{n(1-\beta) - m(1+\beta)}
{2n \beta} \right) ^{\frac{1}{n-m}} K > K .
\ee

\underline{\em Proof of (i)\/.}
Given any $z \in \mbox{} ]0 , \tilde{z} \wedge \gzb[$,
the calculations in (\ref{F_x(K,z)}) imply that (\ref{Fsh2})
is true with $\underline{x}_\dagger (z) < K <
\overline{x}_\dagger (z)$ as well as that $F(K; z) < 0$.
Also, (\ref{F(z,z)}) implies that $F(z,z) < 0$.
Combining these observations with (\ref{cp-I-II}),
(\ref{Flims+2}) and the relevant inequality in (\ref{F-gzb}),
we can see that there exists a unique $\upalpha (z)
\in \mbox{} ]\gzb, \gzo[$ such that (\ref{gotha-i2}) holds
true.

If $ \tilde{z} < \gzb$ and $z \in [\tilde{z}, \gzb[$,
then (\ref{F_x(K,z)}) implies that $F(K; z) \geq 0$
and $\frac{\partial F}{\partial x} (K; z) \geq 0$.
In this case, the inequality $F(z,z) < 0$, which follows
from (\ref{F(z,z)}), implies that (\ref{Fsh2}) is true with
$\underline{x}_\dagger (z) < z$.
This observation, (\ref{Flims+2}) and the relevant
inequality in (\ref{F-gzb}) imply that there exists
a unique $\upalpha (z) \in \mbox{} ]\gzb, \gzo[$ such that
(\ref{gotha-i2}) holds true.

Differentiating the identity $F \bigl( \upalpha (z) ; z \bigr)
= 0$ with respect to $z$, and using (\ref{Fz()}),
the inequalities $\frac{mK}{m-1} < \gzc \leq
\gzb < \upalpha (z)$ (see also (\ref{cp-I-II})) and
(\ref{gotha-i2}), we obtain
\be
\upalpha ' (z) = - \frac{\frac{\partial F}{\partial z} \bigl(
\upalpha (z); z \bigr)} {\frac{\partial F}{\partial x} \bigl(
\upalpha (z); z \bigr)} < 0 \quad \text{for all } z \in \mbox{}
]0, \gzb[ ,
\ee
which proves that $\upalpha$ is strictly decreasing.
Furthermore, the first limit in (\ref{gotha-i1}) follows
from the calculation
\ben
0 = \lim _{z \rightarrow 0} F \bigl( \upalpha (z); z \bigr)
= \lim _{z \rightarrow 0} \bigl[ (n-1) \upalpha (z) - nK
\bigr] A \upalpha ^{n-m} (z) , \label{a(0)-calc}
\een
while, the second limit in (\ref{gotha-i1}) follows
from (\ref{gotha-i2}) and (\ref{F-gzb}).

\underline{\em Proof of (ii)\/.}
We first note that, in this case,
\ben
\frac{mK}{m-1} < K < x_\dagger < \gzc , \label{K<xdag}
\een
where $x_\dagger$ is given by (\ref{xdag})
(see (\ref{rb-nm}), (\ref{n+m-1/b}) and the statement
of Lemma~\ref{lem:psi}.(III)).
Combining this observation with (\ref{Fsh2}) and the
identities $F(K; \tilde{z}) = \frac{\partial F}{\partial x}
(K; \tilde{z}) = 0$, which follow from (\ref{F_x(K,z)}),
we can see that
\ben
K = \underline{x} _\dagger (\tilde{z}) < x_\dagger
\quad \text{and} \quad F(x_\dagger ; \tilde{z}) < 0
. \label{III-K-daggs}
\een

Using the definition (\ref{F}) of $F$ and (\ref{Fa()}),
we calculate
\begin{align}
F (\gzc ; z) & = - \frac{rK}{r-bm} \left[ A \gzc ^{n-m}
+ \frac{m(m-1)}{n(n-1)} B(z) \right] , \nonumber \\
\text{and} \quad
\frac{\partial F}{\partial x} (\gzc ; z) & = -
\frac{n (r-bn)}{r} \left[ A \gzc ^{n-m} + \frac{m(m-1)}{n(n-1)}
B(z) \right] \nonumber \\
& = \frac{n (r-bn) (r-bm)}{r^2 K} F (\gzc ; z) ,
\nonumber
\end{align}
It follows that
\ben
F (\gzc ; z) , \ \frac{\partial F}{\partial x} (\gzc ; z)
\begin{cases}
< 0 & \text{for all } z \in \mbox{} ]0, \gc \gzc[ , \\
> 0 & \text{for all } z \in \mbox{} ]\gc \gzc, \gzc[ ,
\end{cases} \label{III-gzc}
\een
where $\gc$ is defined by (\ref{gc}).
These inequalities and (\ref{K<xdag}) imply that
\be
x_\dagger < \overline{x} _\dagger (\gc \gzc) = \gzc
\quad \text{and} \quad F(x_\dagger ; \gc \gzc) > 0 .
\ee
Combining this result with (\ref{III-K-daggs}) and the
fact that $\frac{\partial F}{\partial z} (x_\dagger; z)
> 0$, which follows from (\ref{Fz()}) and (\ref{K<xdag}),
we can see that $\tilde{z} < \gc \gzc$, which implies
that $K < \gc \gzc$.

Given any $z \in \mbox{} ]0, \gc \gzc[$, the inequalities
in (\ref{III-gzc}) imply that (\ref{Fsh2}) is true with
$\underline{x}_\dagger (z) < \gzc < \overline{x}_\dagger
(z)$.
It follows that, given any $z \in \mbox{} ]0, \gc \gzc[$,
there exists a unique $\upalpha (z) \in \mbox{} ]\gzc, \gzo[$
such that (\ref{gotha-ii2}) holds true.
Furthermore, (\ref{gotha-ii2}), (\ref{a(0)-calc}) and
(\ref{III-gzc}) imply the limits in (\ref{gotha-ii1}).
Taking into account the inequalities $\frac{mK}{m-1}
< \gzc < \upalpha (z)$, we can show that $\upalpha$ is
strictly decreasing in the same way as in Part~(i).
The inequality (\ref{F>0-ii}) follows from (\ref{Fz()}), the
inequality $\frac{mK}{m-1} < \gzc$ and the fact that
\be
F(x; \gc \gzc) > 0 \quad \text{for all } x > \gzc
\ee
(see also (\ref{gotha-ii1}) and (\ref{gotha-ii2})).

\underline{\em Proof of (iii)\/.}
If $n \leq \frac{1+\beta}{1-\beta}$, then (\ref{Flims+1}),
(\ref{Flims+2}) and (\ref{F(z,z)}) imply that
\be
\lim _{x \rightarrow 0} F(x;z) > 0 , \quad
F(z;z) \leq 0 , \quad F(\gzo; z) < 0 \quad \text{and}
\quad \lim _{x \rightarrow \infty} F(x;z) = \infty .
\ee
If $n > \frac{1+\beta}{1-\beta}$, which implies that 
$\gzo < \gzb$ (see (\ref{cp-IV})), then (\ref{Flims+1}),
(\ref{Flims+2}), (\ref{F(z,z)}) and (\ref{F-gzb})
imply that
\begin{gather}
\lim _{x \rightarrow 0} F(x;z) > 0 , \quad
F(\gzo; z) < 0 , \quad \lim _{x \rightarrow \infty}
F(x;z) = \infty , \nonumber \\
F(z;z) \begin{cases} < 0 , & \text{if } z \in \mbox{}
]0, \gzb[ , \\ > 0 , & \text{if } z \in \mbox{} ]\gzb,
\infty[ , \end{cases} \quad \text{and} \quad
F(\gzb; z) \begin{cases} > 0 , & \text{if } z \in \mbox{}
]0, \gzb[ , \\ < 0 , & \text{if } z \in \mbox{} ]\gzb,
\infty[ . \end{cases} \nonumber
\end{gather}
Combining these observations with (\ref{Fsh3}),
we can see that, given any $z > 0$, there exists
a unique $\upalpha (z) > z \vee \gzo$ such that
(\ref{gotha-iii1}) and (\ref{gotha-iii2}) both hold
true.
\mbox{}\hfill$\Box$
\bigskip

The next result addresses the inequality
\ben
g(x,z) := \frac{\upalpha (z) - K}{\psi \bigl( \upalpha
(z) ; z \bigr)} - \frac{x-K}{\psi (x;z)} \geq 0 ,
\quad \text{for } x \in \mbox{} ]0, \upalpha (z)] ,
\label{stop-ineq}
\een
which will play an important role in our analysis
in the next section.

\begin{lem} \label{lem:HJB-ineq}
Suppose that Assumption~\ref{assm} holds true,
and let the function $\upalpha$ be as in
Lemma~\ref{lem:a(z)}.
The following statements hold true:
\smallskip

\noindent {\rm (i)}
Given any $z \in \mbox{} ]0, \gzb[$, the inequality
(\ref{stop-ineq}) holds true for all $x \in \mbox{}
]0, \upalpha (z)]$ if the problem's parameters are
as in Cases~(I) or~(II) of Lemma~\ref{lem:psi}.
\smallskip

\noindent {\rm (ii)}
If the problem's parameters are as in Case~(III)
of Lemma~\ref{lem:psi}, then there exists a unique
point $\zminus \in \mbox{} ]K, \gc \gzc[$ such that
\ben
g(z,z) = \frac{\upalpha (z) - K} {\psi \bigl( \upalpha
(z) ; z \bigr)} - \frac{z - K} {z^n} \begin{cases} > 0 ,
& \text{if } z \in \mbox{} ]0, \zminus[ , \\ < 0 , & \text{if }
z \in \mbox{} ]\zminus, \gc \gzc[ . \end{cases}
\label{zo-minus}
\een
Given any $z \in \mbox{} ]0, \zminus]$,
the inequality (\ref{stop-ineq}) holds true for all
$x \in \mbox{} ]0, \upalpha (z)]$.
Furthermore, there exists a function $\gz : [\zminus,
\gc \gzc[ \mbox{} \rightarrow [\zminus, \gzc[$ such that
\ben
\gz (\zminus) = \zminus, \quad z < \gz (z) \quad
\text{and} \quad g \bigl( \gz (z) , z \bigr) = 0
\quad \text{for all } z \in \mbox{} ]\zminus, \gc \gzc[ .
\label{gz-props}
\een

\noindent {\rm (iii)}
If $\beta \in \mbox{} ]0,1[$ (Case~(IV) of
Lemma~\ref{lem:psi}), then there exists a unique
point $\zplus \in \mbox{} ]\gzo , \infty[$ such that
\ben
g(\gzo, z) = \frac{\upalpha (z) - K} {\psi \bigl( \upalpha
(z) ; z \bigr)} - \frac{\gzo - K} {\psi(\gzo; z)} \begin{cases}
> 0 , & \text{if } z \in \mbox{} ]0, \zplus[ , \\ < 0 , &
\text{if } z \in \mbox{} ]\zplus, \infty[ . \end{cases}
\label{zo-plus}
\een
Furthermore, given any $z \in \mbox{} ]0, \zplus]$,
the inequality (\ref{stop-ineq}) holds true for all
$x \in \mbox{} ]0, \upalpha (z)]$.

\end{lem}
{\bf Proof.}
We first calculate
\ben
g \bigl( \upalpha (z) , z) = 0 \quad \text{and} \quad
\frac{\partial g}{\partial x} (x,z) = \begin{cases}
(n-1) (x - \gzo) x^{-1-n} , & \text{if } x \in \mbox{}
]0, z[ , \\ x^{m-1} \psi^{-2} (x;z) F(x;z) , & \text{if }
x \in \mbox{} ]z, \upalpha (z)[ ,  \end{cases}
\label{(x-K)/psi-der}
\een
for all $z$ in the domain of $\upalpha$.
(Note that $\frac{\partial g}{\partial x} (z,z)$
does not exist, the corresponding left and
right partial derivatives of $g (\cdot, z)$
are discontinuous at $z$.)

\underline{\em Proof of (i)\/.}
This case follows immediately from (\ref{cp-I-II}),
(\ref{gotha-i2}), (\ref{(x-K)/psi-der}) and the fact that
$\upalpha (z) > K$.

\underline{\em Proof of (ii)\/.}
We first recall that, in this case, $K < \gc \gzc$ (see
Lemma~\ref{lem:a(z)}.(ii)).
In view of (\ref{F-psi}) and Lemma~\ref{lem:psi}.(III),
we can see that
\ben
\frac{\partial}{\partial x} \bigl( x^{m-1} F(x;z) \bigr)
= (x-K) \psi '' (x;z) \begin{cases}
< 0 , & \text{if } z < K \text{ and } x \in \mbox{}
]K, K \vee \gc^{-1} z[ , \\
> 0 , & \text{if } z < K \text{ and } x \in \mbox{}
]K \vee \gc^{-1} z , \infty[ , \\
< 0 , & \text{if } z \geq K \text{ and } x \in \mbox{}
]z , \gc^{-1} z[ , \\
> 0 , & \text{if } z \geq K \text{ and } x \in \mbox{}
]\gc^{-1} z, \infty[ \end{cases} \label{F''-signs}
\een
(in the inequalities here, we list only the cases
we will use).
Also, given any $z \in \mbox{} ]0, \gc \gzc[$, we use
the identities in (\ref{(x-K)/psi-der}) to calculate
\begin{align}
g(x,z) & = g \bigl( \upalpha (z) , z \bigr) - \int
_x^{\upalpha (z)} \frac{\partial g}{\partial x} (y,z)
\, dy \nonumber \\
& = - \int _x^{\upalpha (z)} \frac{y^{m-1} F(y;z)}
{\psi^2 (y;z)} \, dy , \quad \text{for } x \in [z ,
\upalpha (z)] . \label{ii-F(x,z)}
\end{align}

Defining
\be
\upalpha (\gc \gzc) := \lim _{z \rightarrow \gc \gzc}
\upalpha (z) \stackrel{(\ref{gotha-ii1})}{=} \gzc ,
\ee
we note that (\ref{gotha-ii1}), (\ref{gotha-ii2}) imply
that $F(\gzc; \gc \gzc) = 0$.
This observation, the fact that $K < \gc \gzc$ and
the second pair of inequalities in (\ref{F''-signs})
imply that
\be
F(x; \gc \gzc) > 0 \quad \text{for all } x \in \mbox{}
]\gc \gzc , \gzc[ .
\ee
This inequality and (\ref{ii-F(x,z)}) imply that
\ben
g(\gc \gzc , \gc \gzc) = - \int _{\gc \gzc}^{\gzc}
\frac{y^{m-1} F(y;z)}{\psi^2 (y;z)} \, dy < 0 .
\label{g(gc.gzc)<0}
\een
Combining this result with the observation that
$g(z,z) > 0$ for all $z \leq K$, which follows from
the definition (\ref{stop-ineq}) of $g$ and the fact that
$\upalpha (z) > \gzc > K$ for all $z < \gc \gzc$, we
can see that
\be
\zminus = \inf \bigl\{ z \in \mbox{} ]0, \gc \gzc] \mid
\ g(z,z) \leq 0 \bigr\} \in \mbox{} ]K, \gc \gzc[ .
\ee

We will establish (\ref{zo-minus}) if we show that
$g(z,z) < 0$ for all $z \in \mbox{} ]\zminus, \gc \gzc[$.
To this end, we differentiate the expression of the
function $]0, \gc \gzc[ \mbox{} \ni z \rightarrow \bar{g}
(z) := g(z,z)$ given by (\ref{zo-minus}) and we use the
identities
\be
F \bigl( \upalpha (z) ; z \bigr) = 0 \quad
\stackrel{(\ref{F-psi})}{\Rightarrow} \quad
\bigl( \upalpha (z) - K \bigr) \psi' \bigl( \upalpha (z);
z \bigr) = \psi \bigl( \upalpha (z); z \bigr)
\ee
to calculate
\begin{align}
\bar{g}' (z) & =
- \frac{\upalpha^{m-1} (z) F \bigl(\upalpha (z) ; z \bigr)}
{\psi^2 \bigl( \upalpha (z) ; z \bigr)} \upalpha' (z)
- \frac{\bigl( \upalpha (z) - K \bigr) \frac{\partial \psi}
{\partial z} \bigl( \upalpha (z); z \bigr)}
{\psi^2 \bigl( \upalpha (z); z \bigr)}
- (n-1) (\gzo - z) z^{-n-1} \nonumber \\
& = - \frac{2n \beta \bigl( \upalpha (z) - K \bigr) z^{n-m-1}
\upalpha^m (z)} {(1+\beta) \psi^2 \bigl( \upalpha (z); z \bigr)}
- (n-1) (\gzo - z) z^{-n-1} \nonumber
\end{align}
and
\begin{align}
\bar{g}'' (z) = \mbox{} &
- \frac{2 n (m-1) \beta z^{n-m-1} \upalpha^{m-1} (z)
\upalpha' (z)}{(1+\beta) \psi^2 \bigl( \upalpha (z); z \bigr)}
\left[ \upalpha(z) - \frac{mK}{m-1} \right] \nonumber\\
& - \frac{2 n \beta z^{n-m-2} \bigl( \upalpha(z) - K \bigr)
\upalpha^m (z)}{(1+\beta) \psi^3 \bigl( \upalpha (z); z \bigr)}
\left[ (n-m-1) \psi \bigl( \upalpha(z); z \bigr) 
- \frac{4 n \beta z^{n-m} \upalpha^m (z)} {1+\beta} \right]
\nonumber \\
& + n(n-1) \left[ \frac{n+1}{n} \gzo - z \right] z^{-n-2} .
\nonumber
\end{align}
In view of the fact that $\upalpha : \mbox{} ]0,
\gc \gzc[ \mbox{} \rightarrow \mbox{} ]\gzc, \gzo[$
is strictly decreasing and the inequalities
$\frac{mK}{m-1} < K < \zminus < \gc \gzc < \gzc < \gzo$
(see (\ref{cp-III}) in Lemma~\ref{lem:psi} and
Lemma~\ref{lem:a(z)}.(ii)), the latter expression implies
that the function $]0, \gc \gzc[ \mbox{} \ni z \rightarrow
\bar{g} (z)$ is strictly convex.
Combining this observation with the inequalities
$\bar{g} (\zminus) = 0$ and $\bar{g} (\gc \gzc) < 0$
(see (\ref{g(gc.gzc)<0}) for the last one), we can see
that $\bar{g} (z) \equiv g(z,z) < 0$ for all $z \in \mbox{}
]\zminus, \gc \gzc]$, as required.

To proceed further, we note that, if $z < K$, then the
expression (\ref{F-psi}) of $F$ implies that $F(x,z)
< 0$ for all $x \in [z,K]$.
Combining this observation with the identity $F \bigl(
\upalpha (z) ; z \bigr) = 0$ and the first pair of inequalities
in (\ref{F''-signs}), we can see that
\be
\text{given any } z \in \mbox{} ]0, K[ , \quad
F(x;z) < 0 \quad \text{for all } x \in [z, \upalpha (z)[ .
\ee
On the other hand, the second pair of inequalities in
(\ref{F''-signs}) and the identity $F \bigl( \upalpha
(z) ; z \bigr) = 0$ imply that
\begin{gather*}
\text{given any } z \in [K, \gc \gzc[ , \text{ either }
F(x;z) < 0 \quad \text{for all } x \in \mbox{}
]z, \upalpha (z)[ , \nonumber \\
\text{or there exists } x_\star (z)  \in \mbox{} ]z, \gzc[
\text{ such that } F(x;z) \begin{cases}
> 0 , & \text{if } x \in \mbox{} ]z, x_\star (z)[ , \\
< 0 , & \text{if } x \in \mbox{} ]x_\star (z), \upalpha (z)[ .
\end{cases}
\end{gather*}
These observations and (\ref{(x-K)/psi-der}) imply
that either
\ben
\frac{\partial g}{\partial x} (x,z) < 0 \text{ for all }
x \in \mbox{} ]0, \upalpha (z)[ , \quad \text{or} \quad
\frac{\partial g}{\partial x} (x,z) \begin{cases}
< 0  , & \text{if } x \in \mbox{} ]0, z[ \mbox{}
\cup \mbox{} ]x_\star (z), \upalpha (z)[ , \\
> 0 , & \text{if } x \in \mbox{} ]z, x_\star (z)[ .
\end{cases} \label{gmono}
\een
Given any $z \in \mbox{} ]0, \zminus]$, the inequality
$g(z,z) > 0$ (see (\ref{zo-minus})), the identity
$g \bigl( \upalpha (z), z \bigr) = 0$ and (\ref{gmono})
imply that (\ref{stop-ineq}) holds true for all
$x \in [0, \upalpha (z)]$.
On the other hand, given any $z \in \mbox{}
[\zminus, \gc \gzc[$, the inequality $g(z,z) < 0$
(see (\ref{zo-minus})), the identity $g \bigl(
\upalpha (z), z \bigr) = 0$ and (\ref{gmono})
imply that there exists a unique $\gz (z) \in
[\zminus, \gzc[$ such that (\ref{gz-props})
holds true (note that $g$ is as in the second
case of (\ref{gmono}) here).

\underline{\em Proof of (iii)\/.}
In this case, (\ref{cp-IV}) and (\ref{gotha-iii1}) imply
that
\be
\lim _{z \rightarrow \infty} \upalpha (z) = \infty ,
\text{ if } n \leq \frac{1+\beta}{1-\beta} ,
\quad \text{and} \quad
\upalpha (\gzb) = \gzb > \gzo , \text{ if } n >
\frac{1+\beta}{1-\beta} .
\ee
It follows that
\begin{gather}
\lim _{z \rightarrow \infty} g(\gzo, z) = - \frac{\gzo - K}
{\gzo^n} < 0 , \quad \text{if } n \leq \frac{1+\beta}{1-\beta}
, \nonumber \\
\text{and} \quad
g(\gzo, \gzb) = \frac{\gzb - K} {\gzb^n} - \frac{\gzo - K}
{\gzo^n} < 0 , \quad \text{if } n > \frac{1+\beta}{1-\beta}
. \nonumber
\end{gather}
On the other hand, (\ref{gotha-iii2}) and (\ref{(x-K)/psi-der})
imply that
\be
g(\gzo, \gzo) = - \int _{\gzo}^{\upalpha (\gzo)}
\frac{\partial g}{\partial x} (y,z) \, dy = - \int
_{\gzo}^{\upalpha (\gzo)} \frac{y^{m-1} F(y;z)}
{\psi^2 (y;z)} \, dy > 0 .
\ee
Combining these observations with the calculation
\begin{align}
\frac{\partial}{\partial z} g(\gzo, z) & = -
\frac{\upalpha^{m-1} (z) F \bigl( \upalpha (z) ; z \bigr)}
{\psi^2 \bigl( \upalpha (z) ; z \bigr)} \upalpha' (z)
- \frac{\bigl( \upalpha (z) - K \bigr) \frac{\partial \psi}
{\partial z} \bigl( \upalpha (z); z \bigr)}
{\psi^2 \bigl( \upalpha (z); z \bigr)}
+ \frac{(\gzo - K) \frac{\partial \psi}{\partial z}
(\gzo; z)}{\psi^2 (\gzo; z)} \nonumber \\
& = - \frac{2n \beta \bigl( \upalpha (z) - K \bigr)
z^{n-m-1} \upalpha^m (z)} {(1+\beta) \psi^2
\bigl( \upalpha (z); z \bigr)} < 0 , \quad \text{for }
z \in \mbox{} ]\gzo , \infty[
, \nonumber
\end{align}
we can see that there exists a unique $\zplus
\in \mbox{} ]\gzo , \infty[$ such that (\ref{zo-plus})
holds true.

Finally, we fix any $z \in \mbox{} ]0, \zplus]$.
In view of the inequality $g(\gzo, z) \geq 0$,
the identity $g \bigl( \upalpha (z) , z) = 0$ and
the observation that
\be
\frac{\partial g}{\partial x} (x,z) \begin{cases}
> 0 , & \text{if } x \in \mbox{} ]\gzo, z[ , \\ < 0 , &
\text{if } x \in \mbox{} ]0, \gzo[ \mbox{} \cup \mbox{}
]z, \upalpha (z)[ ,  \end{cases}
\ee
which follows from (\ref{gotha-iii2}) and
(\ref{(x-K)/psi-der}), we can see that the inequality
(\ref{stop-ineq}) holds true for all $x \in \mbox{}
]0, \upalpha (z)[$.
\mbox{}\hfill$\Box$

\begin{rem} \label{rem:zminus...} {\rm
Our analysis in the next sections will make use of the
following observation.
Suppose that the problem's parameters are as in Case~(III)
of Lemma~\ref{lem:psi} and fix any $z \in [\zminus, \gc
\gzc[$, where $\zminus$ is as in Lemma~\ref{lem:HJB-ineq}.(ii).
The function $u (\cdot ; z) : [z, \infty[ \mbox{} \rightarrow \bbr$
defined by
\be
u(x ; z) = \Gamma (z) \psi (x; z) - (x-K)
= A \Gamma (z) x^n + B (z) \Gamma (z) x^m - (x-K) ,
\ee
where $\Gamma (z)  = \bigl( \upalpha (z) - K \bigr)
/ \psi \bigl( \upalpha (z) ; z \bigr)$, is such that
\be
u \bigl( \gz (z) ; z \bigr) = 0 \quad \text{and} \quad
u \bigl( \upalpha (z) ; z \bigr) = \frac{\partial u}{\partial x}
\bigl( \upalpha (z) ; z \bigr) = 0 .
\ee
The first of these identities follows immediately from
(\ref{gz-props}) and the fact that $u(x;z) = g(x,z)
\psi (x; z)$ for all $x \geq z$.
On the other hand, the identities for $x = \upalpha (z)$
hold true because they are equivalent to the identity
$F \bigl( \upalpha (z) ; z \bigr) = 0$.
} \mbox{}\hfill$\Box$ \end{rem}


\section{The solution to the optimal stopping problem
defined by (\ref{XGBM}) and (\ref{vEx})}
\label{sec:GBMsol}

We expect that the value function $v$ of the discretionary
problem defined by (\ref{XGBM}) and (\ref{vEx}) should
be strictly positive.
Combining this observation with the fact that the restriction
of the function $x \mapsto (x-K)^+$ to $\bbr_+ \setminus
\{ K \}$ is $C^\infty$ and the so-called ``principle of
smooth fit'' (see also Remark~\ref{rem:SF}), we
expect that the restriction of $v$ to $]0, \infty[ \mbox{}
\setminus \{ z \}$ should be $C^1$ with absolutely
continuous first derivative.
In view of (\ref{VI-sc-cond}), (\ref{VI-sc2}) in
Example~\ref{ex:VI}, we therefore expect that $v$
should identify with a function $w$ satisfying
\begin{gather}
\max \left\{ \half \sigma^2 x^2 w'' (x) + bx w' (x) - r w(x)
, \ (x-K)^+ - w(x) \right\} = 0 , \text{ inside } \mbox{}
]0,z[ \mbox{} \cup \mbox{} ]z,\infty[ , \label{VI1} \\
\text{and} \quad
\max \Big\{ (1+\beta) w_+' (z) - (1-\beta) w_-' (z)
, \ (z-K)^+ - w(z) \Big\} = 0 . \label{VI2}
\end{gather}
Furthermore, the strict positivity of $v$ and
Remark~\ref{rem:VIuse} imply that the waiting
region includes the interval $\bigl] 0, K
\vee rK/(r-b) \bigr[$ and, if $\beta \in \mbox{} ]0,1[$,
then $z$ also belongs to the waiting region.

We now solve the optimal stopping problem we consider
in this section by constructing an appropriate solution
to the variational inequality (\ref{VI1})--(\ref{VI2}).
In its simplest form, we expect that the required solution
has the same qualitative form as the solution to the optimal
stopping problem associated with the usual perpetual 
American call option, which involves a standard geometric
Brownian motion ($\beta = 0$).
Accordingly, we expect that the value function $v$ should
identify with the function
\ben
w(x) = w(x;z) = \begin{cases} \Gamma (z) \psi (x; z)
, & \text{if } x \leq a , \\ x-K , & \text{if } x > a , \end{cases}
\label{w-}
\een
for some constants $a = a(z) > 0$ and $\Gamma (z) > 0$,
while
\ben
\tau _\star = \inf \bigl\{ t \geq 0 \mid \ X_t \geq a \bigr\}
\label{optau}
\een
should identify an optimal stopping time.
It turns out that this is indeed the case for a wide range
of parameter values (see Figures~4-10).
To determine the constant $\Gamma (z)$ and the
free-boundary point $a$, we first appeal to the continuity
of the value function, which yields the expression
\ben
\Gamma (z) = (a-K) \psi^{-1} (a; z) . \label{Gamma}
\een
With the exception of the possibilities depicted by
Figures~5 and~8, we expect that the value function
should be $C^1$ at $a$, which gives rise to the equation
$\Gamma (z) \psi' (a;z) = 1$.\footnote{
Recall that we have adopted the notation $\psi' (x;z)
= \frac{\partial \psi}{\partial x} (x;z)$.}
This equation and (\ref{Gamma}) imply that $a$
should satisfy equation (\ref{Feqn}) if $z<a$ (see
Figures~4, 7 and~10) and should be given by
\be
a = \frac{nK}{n-1} \stackrel{(\ref{cp})}{=:} \gzo > 0
\ee
if $a<z$ (see Figures~6 and~9).

The following result, which we prove in
Section~\ref{sec:proofs},
involves the parameters $\gzc$, $\gzb$, $\gzo$ and
$\zminus$, $\zplus$ that are as in (\ref{cp}) and
Lemma~\ref{lem:HJB-ineq}.(ii)-(iii), respectively.

\begin{thm} \label{propSol1}
Consider the optimal stopping problem defined by
(\ref{XGBM}), (\ref{vEx}) and suppose that
Assumption~\ref{assm} holds true. 
If the problem parameters are as in Cases~(I) or~(II)
of Lemma~\ref{lem:psi}, define
\ben
a = \begin{cases} \upalpha (z) , & \text{if } z \in
\mbox{} ]0, \gzb[ , \\ z \wedge \gzo , & \text{if } z
\in \mbox{} [\gzb, \infty[ , \end{cases} \label{a-1}
\een
where the function $\upalpha$ is as in
Lemma~\ref{lem:a(z)}.(i).
If the problem parameters are as in Case~(III) of 
Lemma~\ref{lem:psi}, suppose that $z \in \mbox{}
]0, \zminus] \cup [\gzc, \infty[$ and define
\ben
a = \begin{cases} \upalpha (z) , & \text{if } z \in \mbox{}
]0, \zminus] , \\ z \wedge \gzo , & \text{if } z \in \mbox{}
[\gzc, \infty[ , \end{cases} \label{a-2}
\een
where the function $\upalpha$ is as in
Lemma~\ref{lem:a(z)}.(ii) and $\zminus$ is as in
Lemma~\ref{lem:HJB-ineq}.(ii).
If $\beta \in \mbox{} ]0, 1[$ (Case~(IV) of
Lemma~\ref{lem:psi}), suppose that $z \in \mbox{}
]0, \zplus]$ and define
\ben
a = \upalpha (z) , \quad \text{for } z \in \mbox{}
]0, \zplus] , \label{a-3}
\een
where the function $\upalpha$ is as in
Lemma~\ref{lem:a(z)}.(iii) and $\zplus$ is as in
Lemma~\ref{lem:HJB-ineq}.(iii).
For such choices of $a$ and for $\Gamma (z) > 0$
given by (\ref{Gamma}), the function $w$ defined
by (\ref{w-}) identifies with the value  function $v$
of the discretionary stopping problem and the stopping
time given by (\ref{optau}) is optimal.
\end{thm}

In the context of (\ref{a-1}), we can see that
Figure~4 transforms ``continuously'' into
Figure~5 and then into Figure~6
as $z$ increases from 0 to $\infty$, thanks to
the second limit in (\ref{gotha-i1}).
\bigskip

\begin{picture}(160,110)
\put(20,15){\begin{picture}(120,95) 

\put(0,0){\line(1,0){120}}
\put(120,0){\line(0,1){95}}
\put(0,0){\line(0,1){95}}
\put(0,95){\line(1,0){120}}

\put(10,10){\vector(1,0){100}}
\put(10,10){\vector(0,1){80}}
\put(6,88){$y$}
\put(108,6){$x$}

\put(0,0){\qbezier(30,10)(70,40)(110,70)}

\color{blue}
\put(0,-0.2){\qbezier(10,10)(15,10.5)(25,12)}
\put(0,-0.1){\qbezier(10,10)(15,10.5)(25,12)}
\put(0,0){\qbezier(10,10)(15,10.5)(25,12)}
\put(0,0.1){\qbezier(10,10)(15,10.5)(25,12)}
\put(0,0.1){\qbezier(10,10)(15,10.5)(25,12)}

\put(0,-0.2){\qbezier(25,12)(40,17.5)(50,25)}
\put(0,-0.1){\qbezier(25,12)(40,17.5)(50,25)}
\put(0,0){\qbezier(25,12)(40,17.5)(50,25)}
\put(0,0.1){\qbezier(25,12)(40,17.5)(50,25)}
\put(0,0.2){\qbezier(25,12)(40,17.5)(50,25)}

\color{red}
\put(0,-0.2){\qbezier(50,25)(70,40)(110,70)}
\put(0,-0.1){\qbezier(50,25)(70,40)(110,70)}
\put(0,0.1){\qbezier(50,25)(70,40)(110,70)}
\put(0,0.2){\qbezier(50,25)(70,40)(110,70)}
\color{black}

\put(0,0){\qbezier(50,25)(67,37)(80,72)}

\put(30,9){\line(0,1){2}}
\put(28.5,5.5){$K$}

\put(50,9){\line(0,1){2}}
\put(0,0){\qbezier[20](50,10)(50,17.5)(50,25)}
\put(49,5.5){$a$}

\put(25,9){\line(0,1){2}}
\put(0,0){\qbezier[4](25,10)(25,11)(25,12)}
\put(23.5,5.5){$z$}

\put(37,35.5){\vector(1,-3){5.1}}
\put(39.8,39){\vector(1,0){28.5}}
\put(30,38){$v(x)$}

\put(65,60){\vector(2,-1){8}}
\put(56,60){$\psi(x)$}

\put(92,41){\vector(-1,3){4}}
\put(87,38){$x-K$}

\end{picture}}

\put(20,10){\small{{\bf Figure 4.} The value function
$v$ if the problem parameters are as in}}
\put(40,5){\small{Cases~(I) or~(II)
of Lemma~\ref{lem:psi} and $z \in \mbox{} ]0, \gzb[$.}}
\end{picture}

\begin{picture}(160,110)
\put(20,15){\begin{picture}(120,95) 

\put(0,0){\line(1,0){120}}
\put(120,0){\line(0,1){95}}
\put(0,0){\line(0,1){95}}
\put(0,95){\line(1,0){120}}

\put(10,10){\vector(1,0){100}}
\put(10,10){\vector(0,1){80}}
\put(6,88){$y$}
\put(108,6){$x$}

\qbezier(30,10)(70,40)(110,70)

\color{blue}
\put(0,-0.2){\qbezier(10,10)(30,15)(50,25)}
\put(0,-0.1){\qbezier(10,10)(30,15)(50,25)}
\put(0,0){\qbezier(10,10)(30,15)(50,25)}
\put(0,0.1){\qbezier(10,10)(30,15)(50,25)}
\put(0,0.1){\qbezier(10,10)(30,15)(50,25)}

\color{red}
\put(0,-0.2){\qbezier(50,25)(70,40)(110,70)}
\put(0,-0.1){\qbezier(50,25)(70,40)(110,70)}
\put(0,0.1){\qbezier(50,25)(70,40)(110,70)}
\put(0,0.2){\qbezier(50,25)(70,40)(110,70)}
\color{black}

\put(0,0){\qbezier(50,25)(67,40)(80,75)}

\put(30,9){\line(0,1){2}}
\put(28,5.5){$K$}

\put(50,9){\line(0,1){2}}
\put(0,0){\qbezier[20](50,10)(50,17.5)(50,25)}
\put(45,5.5){$a=z$}

\put(37,35.5){\vector(1,-3){4.7}}
\put(39.8,39){\vector(1,0){28.5}}
\put(30,38){$v(x)$}

\put(63,60){\vector(2,-1){8.5}}
\put(54,60){$\psi(x)$}

\put(92,41){\vector(-1,3){4}}
\put(87,38){$x-K$}

\end{picture}}

\put(20,10){\small{{\bf Figure 5.} The value function
$v$ if the problem parameters are as in}}
\put(40,5){\small{Cases~(I) or~(II)
of Lemma~\ref{lem:psi} and $z \in \mbox{} [\gzb, \gzo]$.}}
\end{picture}

\begin{picture}(160,110)
\put(20,15){\begin{picture}(120,95) 

\put(0,0){\line(1,0){120}}
\put(120,0){\line(0,1){95}}
\put(0,0){\line(0,1){95}}
\put(0,95){\line(1,0){120}}

\put(10,10){\vector(1,0){100}}
\put(10,10){\vector(0,1){80}}
\put(6,88){$y$}
\put(108,6){$x$}

\qbezier(30,10)(70,40)(110,70)

\color{blue}
\put(0,-0.2){\qbezier(10,10)(30,11)(50,25)}
\put(0,-0.1){\qbezier(10,10)(30,11)(50,25)}
\put(0,0){\qbezier(10,10)(30,11)(50,25)}
\put(0,0.1){\qbezier(10,10)(30,11)(50,25)}
\put(0,0.1){\qbezier(10,10)(30,11)(50,25)}

\color{red}
\put(0,-0.2){\qbezier(50,25)(70,40)(110,70)}
\put(0,-0.1){\qbezier(50,25)(70,40)(110,70)}
\put(0,0.1){\qbezier(50,25)(70,40)(110,70)}
\put(0,0.2){\qbezier(50,25)(70,40)(110,70)}
\color{black}

\put(0,0){\qbezier(50,25)(55,28.5)(60,35)}
\put(0,0){\qbezier(60,35)(65,45)(70,70)}

\put(30,9){\line(0,1){2}}
\put(28,5.5){$K$}

\put(50,9){\line(0,1){2}}
\put(0,0){\qbezier[20](50,10)(50,17.5)(50,25)}
\put(45,5.5){$a = \gzo$}

\put(60,9){\line(0,1){2}}
\put(0,0){\qbezier[32](60,10)(60,22.5)(60,35)}
\put(59,5.5){$z$}

\put(35,41.5){\vector(0,-1){24.7}}
\put(39.8,45){\vector(1,0){36.5}}
\put(30,44){$v (x)$}

\put(59,64){\vector(2,-1){8.5}}
\put(50,64){$\psi(x)$}

\put(92,41){\vector(-1,3){4}}
\put(87,38){$x-K$}

\end{picture}}

\put(20,10){\small{{\bf Figure 6.} The value function
$v$ if the problem parameters are as in}}
\put(40,5){\small{Cases~(I) or~(II)
of Lemma~\ref{lem:psi} and $z \in \mbox{} ]\gzo, \infty[$.}}
\end{picture}

\begin{picture}(160,110)
\put(20,15){\begin{picture}(120,95)

\put(0,0){\line(1,0){120}}
\put(120,0){\line(0,1){95}}
\put(0,0){\line(0,1){95}}
\put(0,95){\line(1,0){120}}

\put(10,10){\vector(1,0){100}}
\put(10,10){\vector(0,1){80}}
\put(6,88){$y$}
\put(108,6){$x$}

\put(0,0){\qbezier(30,10)(70,40)(110,70)}

\color{blue}
\put(0,-0.2){\qbezier(10,10)(15,10.5)(25,12)}
\put(0,-0.1){\qbezier(10,10)(15,10.5)(25,12)}
\put(0,0){\qbezier(10,10)(15,10.5)(25,12)}
\put(0,0.1){\qbezier(10,10)(15,10.5)(25,12)}
\put(0,0.1){\qbezier(10,10)(15,10.5)(25,12)}

\put(0,-0.2){\qbezier(25,12)(30,18.5)(38,20)}
\put(0,-0.1){\qbezier(25,12)(30,18.5)(38,20)}
\put(0,0){\qbezier(25,12)(30,18.5)(38,20)}
\put(0,0.1){\qbezier(25,12)(30,18.5)(38,20)}
\put(0,0.2){\qbezier(25,12)(30,18.5)(38,20)}

\put(0,-0.2){\qbezier(38,20)(46,21.9)(50,25)}
\put(0,-0.1){\qbezier(38,20)(46,21.9)(50,25)}
\put(0,0){\qbezier(38,20)(46,21.9)(50,25)}
\put(0,0.1){\qbezier(38,20)(46,21.9)(50,25)}
\put(0,0.2){\qbezier(38,20)(46,21.9)(50,25)}

\color{red}
\put(0,-0.2){\qbezier(50,25)(70,40)(110,70)}
\put(0,-0.1){\qbezier(50,25)(70,40)(110,70)}
\put(0,0.1){\qbezier(50,25)(70,40)(110,70)}
\put(0,0.2){\qbezier(50,25)(70,40)(110,70)}

\color{black}
\put(0,0){\qbezier(50,25)(67,37)(80,72)}

\put(30,9){\line(0,1){2}}
\put(28.5,5.5){$K$}

\put(50,9){\line(0,1){2}}
\put(0,0){\qbezier[20](50,10)(50,17.5)(50,25)}
\put(49,5.5){$a$}

\put(25,9){\line(0,1){2}}
\put(0,0){\qbezier[4](25,10)(25,11)(25,12)}
\put(23.5,5.5){$z$}

\put(37,35.5){\vector(0,-1){15.1}}
\put(39.8,39){\vector(1,0){28.5}}
\put(30,38){$v(x)$}

\put(65,60){\vector(2,-1){8}}
\put(56,60){$\psi(x)$}

\put(92,41){\vector(-1,3){4}}
\put(87,38){$x-K$}

\end{picture}}

\put(20,10){\small{{\bf Figure 7.} The value function
$v$ if the problem parameters are as in}}
\put(40,5){\small{Case~(III) of Lemma~\ref{lem:psi}
and $z \in \mbox{} ]0, \zminus[$.}}
\end{picture}

\begin{picture}(160,110)
\put(20,15){\begin{picture}(120,95)

\put(0,0){\line(1,0){120}}
\put(120,0){\line(0,1){95}}
\put(0,0){\line(0,1){95}}
\put(0,95){\line(1,0){120}}

\put(10,10){\vector(1,0){100}}
\put(10,10){\vector(0,1){80}}
\put(6,88){$y$}
\put(108,6){$x$}

\qbezier(30,10)(70,40)(110,70)

\color{blue}
\put(0,-0.2){\qbezier(10,10)(30,15)(50,25)}
\put(0,-0.1){\qbezier(10,10)(30,15)(50,25)}
\put(0,0){\qbezier(10,10)(30,15)(50,25)}
\put(0,0.1){\qbezier(10,10)(30,15)(50,25)}
\put(0,0.1){\qbezier(10,10)(30,15)(50,25)}

\color{red}
\put(0,-0.2){\qbezier(50,25)(70,40)(110,70)}
\put(0,-0.1){\qbezier(50,25)(70,40)(110,70)}
\put(0,0.1){\qbezier(50,25)(70,40)(110,70)}
\put(0,0.2){\qbezier(50,25)(70,40)(110,70)}

\color{black}
\put(0,0){\qbezier(50,25)(52,34)(58,37)}

\put(0,0){\qbezier(58,37)(77,46)(90,75)}

\put(30,9){\line(0,1){2}}
\put(28,5.5){$K$}

\put(50,9){\line(0,1){2}}
\put(0,0){\qbezier[20](50,10)(50,17.5)(50,25)}
\put(45,5.5){$a=z$}

\put(37,35.5){\vector(1,-3){4.7}}
\put(39.8,39){\vector(1,0){28.5}}
\put(30,38){$v(x)$}

\put(63,60){\vector(2,-1){13.5}}
\put(54,60){$\psi(x)$}

\put(92,41){\vector(-1,3){4}}
\put(87,38){$x-K$}

\end{picture}}

\put(20,10){\small{{\bf Figure 8.} The value function
$v$ if the problem parameters are as in}}
\put(40,5){\small{Case~(III) of Lemma~\ref{lem:psi}
and $z \in \mbox{} [\gzc, \gzo]$.}}
\end{picture}

\begin{picture}(160,110)
\put(20,15){\begin{picture}(120,95) 

\put(0,0){\line(1,0){120}}
\put(120,0){\line(0,1){95}}
\put(0,0){\line(0,1){95}}
\put(0,95){\line(1,0){120}}

\put(10,10){\vector(1,0){100}}
\put(10,10){\vector(0,1){80}}
\put(6,88){$y$}
\put(108,6){$x$}

\qbezier(30,10)(70,40)(110,70)

\color{blue}
\put(0,-0.2){\qbezier(10,10)(30,11)(50,25)}
\put(0,-0.1){\qbezier(10,10)(30,11)(50,25)}
\put(0,0){\qbezier(10,10)(30,11)(50,25)}
\put(0,0.1){\qbezier(10,10)(30,11)(50,25)}
\put(0,0.1){\qbezier(10,10)(30,11)(50,25)}

\color{red}
\put(0,-0.2){\qbezier(50,25)(70,40)(110,70)}
\put(0,-0.1){\qbezier(50,25)(70,40)(110,70)}
\put(0,0.1){\qbezier(50,25)(70,40)(110,70)}
\put(0,0.2){\qbezier(50,25)(70,40)(110,70)}

\color{black}
\put(0,0){\qbezier(50,25)(55,28.5)(60,35)}
\put(0,0){\qbezier(60,35)(61,45)(68,49)}
\put(0,0){\qbezier(68,49)(81,56)(92,80)}

\put(30,9){\line(0,1){2}}
\put(28,5.5){$K$}

\put(50,9){\line(0,1){2}}
\put(0,0){\qbezier[20](50,10)(50,17.5)(50,25)}
\put(45,5.5){$a = \gzo$}

\put(60,9){\line(0,1){2}}
\put(0,0){\qbezier[32](60,10)(60,22.5)(60,35)}
\put(59,5.5){$z$}

\put(35,41.5){\vector(0,-1){24.7}}
\put(39.8,45){\vector(1,0){36.5}}
\put(30,44){$v(x)$}

\put(59,64){\vector(2,-1){16.5}}
\put(50,64){$\psi (x)$}

\put(92,41){\vector(-1,3){4}}
\put(87,38){$x-K$}

\end{picture}}

\put(20,10){\small{{\bf Figure 9.} The value function
$v$ if the problem parameters are as in}}
\put(40,5){\small{Case~(III) of Lemma~\ref{lem:psi}
and $z \in \mbox{} ]\gzo, \infty[$.}}
\end{picture}

\begin{picture}(160,110)
\put(20,15){\begin{picture}(120,95)

\put(0,0){\line(1,0){120}}
\put(120,0){\line(0,1){95}}
\put(0,0){\line(0,1){95}}
\put(0,95){\line(1,0){120}}

\put(10,10){\vector(1,0){100}}
\put(10,10){\vector(0,1){80}}
\put(6,88){$y$}
\put(108,6){$x$}

\qbezier(30,10)(70,40)(110,70)

\color{blue}
\put(0,-0.2){\qbezier(10,10)(30,11)(40,25)}
\put(0,-0.1){\qbezier(10,10)(30,11)(40,25)}
\put(0,0){\qbezier(10,10)(30,11)(40,25)}
\put(0,0.1){\qbezier(10,10)(30,11)(40,25)}
\put(0,0.1){\qbezier(10,10)(30,11)(40,25)}

\put(0,-0.2){\qbezier(40,25)(53,27)(60,32.5)}
\put(0,-0.1){\qbezier(40,25)(53,27)(60,32.5)}
\put(0,0){\qbezier(40,25)(53,27)(60,32.5)}
\put(0,0.1){\qbezier(40,25)(53,27)(60,32.5)}
\put(0,0.2){\qbezier(40,25)(53,27)(60,32.5)}

\color{red}
\put(0,-0.2){\qbezier(60,32.5)(80,47.5)(110,70)}
\put(0,-0.1){\qbezier(60,32.5)(80,47.5)(110,70)}
\put(0,0){\qbezier(60,32.5)(80,47.5)(110,70)}
\put(0,0.1){\qbezier(60,32.5)(80,47.5)(110,70)}
\put(0,0.2){\qbezier(60,32.5)(80,47.5)(110,70)}

\color{black}
\put(0,0){\qbezier(60,32.5)(75,45)(80,70)}

\put(30,9){\line(0,1){2}}
\put(28,5.5){$K$}

\put(60,9){\line(0,1){2}}
\put(0,0){\qbezier[30](60,10)(60,21.25)(60,32.5)}
\put(59,6){$a$}

\put(40,9){\line(0,1){2}}
\put(0,0){\qbezier[22](40,10)(40,17.5)(40,25)}
\put(39,6){$z$}

\put(35,41.5){\vector(0,-1){21.7}}
\put(39.8,45){\vector(1,0){36.5}}
\put(30,44){$v(x)$}

\put(69,64){\vector(2,-1){8.1}}
\put(60,64){$\psi(x)$}

\end{picture}}

\put(20,10){\small{{\bf Figure 10.} The value function
$v$ if the problem parameters are as in}}
\put(42,5){\small{Case~(IV) of Lemma~\ref{lem:psi}
and $z \in \mbox{} ]0, \zplus[$.}}
\end{picture}

\begin{rem} \label{rem:SFex} {\rm
Suppose that the problem parameters are as in
Cases~(I) or~(II) of Lemma~\ref{lem:psi}.
In view of the identity in (\ref{SF-calc}) and
Theorem~\ref{propSol1}, we can see that, given
any $z \in [\gzb, \gzo]$,
\begin{align}
\frac{v_+' (z)}{p_+' (z)} - \frac{v_-' (z)}{p_-' (z)}
& = \frac{1}{(1-\beta) p_-' (z)} \bigl[ 1+\beta -
(1-\beta) v_-' (z) \bigr] \nonumber \\
& = - \frac{n - \frac{1+\beta}{1-\beta}}
{z p_-' (z)} (z - \gzb) \in \left[ \frac{2 n \beta K}
{(n-1) (1-\beta) \gzo p_-' (\gzo)} , \ 0 \right] ,
\nonumber
\end{align}
while
\be
\frac{v_+' (z)}{\psi' (z+; z)} - \frac{v_-' (z)}{\psi' (z-; z)}
= - \frac{n - \frac{1+\beta}{1-\beta}}{n z^n} (z - \gzb)
\in \left[ \frac{2 \beta K}{(n-1) (1-\beta) \gzo^n}
, \ 0 \right] .
\ee
We are thus faced with an example of
``right-sided'' optimal stopping of a skew
geometric Brownian motion in which the
``principle of smooth fit'' does not hold in the
sense that none of $v_-'$, $v_-' / p_-'$ or
$v_-' / \psi_-'$ is continuous.
} \mbox{}\hfill$\Box$ \end{rem}

If the problem parameters are as in Case~(III) of 
Lemma~\ref{lem:psi}, then the function $w = w (\cdot; z)$
given by (\ref{w-}), (\ref{Gamma}) is such that
\begin{align}
w(\zminus; \zminus) & \equiv \Gamma (\zminus)
\psi (\zminus; \zminus) = \zminus - K \nonumber \\
\text{and} \quad
w \bigl( \upalpha (\zminus) ; \zminus \bigr) & \equiv
\Gamma (\zminus) \psi \bigl( \upalpha (\zminus) ;
\zminus \bigr) = \upalpha (\zminus) - K
\label{w-zminus}
\end{align}
(see Lemma~\ref{lem:HJB-ineq}.(ii)).
This observation and the ``singularity'' associated
with $z$ give rise to the following possibility.
For $z \geq \zminus$, the stopping
time
\ben
\underline{\tau}^\star = \inf \bigl\{ t \geq 0 \, \mid \
X_t \in \mbox{}  \{ z \} \mbox{} \cup \mbox{} [\xi, \infty[
\bigr\} , \label{optau-}
\een
where $\xi = \xi (z) > z$ is a constant, may be optimal.
In such a context, we expect that the value function
$v$ should identify with the function
\ben
\underline{w} (x) = \underline{w} (x;z) = \begin{cases}
(z-K) z^{-n} x^n , & \text{if } x \leq z , \\ C (z) x^n + D (z)
x^m , & \text{if } x \in \mbox{} ]z, \xi[ , \\ x-K , & \text{if }
x \geq \xi , \end{cases} \label{wJ}
\een 
for some $C(z), D(z) \in \R$ (see Figure~11). 
To determine the constants $C(z)$, $D(z)$ and the
free-boundary point $\xi = \xi (z)$, we require that
$\underline{w}$ should be $C^1$ at $\xi$, which is
suggested by the ``principle of smooth fit'', as well
as continuous at $z$. 
The system of equations arising from these requirements
is equivalent to the expressions
\begin{gather}
C(z) = - \frac{1}{n-m} \bigl[ (m-1) \xi - mK \bigr] \xi^{-n} ,
\quad D(z) = \frac{1}{n-m} \bigl[ (n-1) \xi - nK \bigr]
\xi^{-m} , \label{CD}
\end{gather}
and the algebraic equation
\begin{equation}
J(\xi; z) = 0 , \label{Jeqn}
\end{equation}
where
\begin{align}
J(x; z) 
& = \bigl[ (n-1) x - nK \bigr] A x^{-m} - \bigl[ (m-1) x - mK \bigr]
A z^{n-m} x^{-n} - (n-m) A (z-K) z^{-m} \nonumber \\
& = x^{-n} F(x; z) - \bigl[ (m-1) x - mK \bigr] z^{n-m} x^{-n}
- (n-m) A (z-K) z^{-m} . \label{J}
\end{align}
To establish the second identity here, we have used
the definitions (\ref{A}), (\ref{B}) of $A$, $B$, as well as
the definition (\ref{F}) of $F$.

We prove the following result in Section~\ref{sec:proofs}.

\begin{thm} \label{prop:xi}
Consider the optimal stopping problem defined by
(\ref{XGBM}), (\ref{vEx}) and suppose that
Assumption~\ref{assm} holds true.
Also, suppose that the problem parameters are as
in Case~(III) of  Lemma~\ref{lem:psi}.
Equation~(\ref{Jeqn}) defines uniquely a strictly
decreasing function $\upxi : \mbox{} ]0, \gzc[ \mbox{}
\rightarrow \mbox{} ]\gzc, \gzo[$ such that
\ben
\lim _{z \rightarrow 0} \upxi (z) = \gzo , \quad
\upxi (\zminus) = \upalpha (\zminus)
\quad \text{and} \quad
\lim _{z \rightarrow \gzc} \upxi (z) = \gzc ,
\label{xi-lims}
\een
where $\upalpha$ is as in Lemma~\ref{lem:a(z)}.(ii)
and $\zminus$ is as in Lemma~\ref{lem:HJB-ineq}.(ii).
Given any $z \in \mbox{} ]\zminus, \gzc[$, the function
$\underline{w}$ defined by (\ref{wJ})--(\ref{CD}) 
for $\xi = \upxi (z)$ identifies with the value function
$v$ of the discretionary stopping problem and the
$(\F_t)$-stopping time $\underline{\tau}^\star$ defined
by (\ref{optau-}) is optimal.
\end{thm}

In the context of the previous result and the part
of Theorem~\ref{propSol1} addressing the case when
the problem parameters are as in Case~(III) of 
Lemma~\ref{lem:psi} (see (\ref{a-2}) in particular),
we can see that Figure~7 transforms ``continuously''
into Figure~11, then into Figure~8 and then into
Figure~9 as $z$ increases from 0 to $\infty$, thanks
to the identities
\ben
\underline{w} (\cdot; \zminus) = w(\cdot; \zminus) ,
\quad \upxi (\zminus) = \upalpha (\zminus) \quad
\text{and} \quad \lim _{z \rightarrow \gzc} \upxi (z)
= \gzc < \gzo
\een
(see (\ref{w-zminus}) and (\ref{xi-lims})).
\bigskip

\begin{picture}(160,110)
\put(20,15){\begin{picture}(120,95)

\put(0,0){\line(1,0){120}}
\put(120,0){\line(0,1){95}}
\put(0,0){\line(0,1){95}}
\put(0,95){\line(1,0){120}}

\put(10,10){\vector(1,0){100}}
\put(10,10){\vector(0,1){80}}
\put(6,88){$y$}
\put(108,6){$x$}

\qbezier(30,10)(70,40)(110,70)

\color{blue}
\put(0,-0.2){\qbezier(10,10)(30,15)(50,25)}
\put(0,-0.1){\qbezier(10,10)(30,15)(50,25)}
\put(0,0){\qbezier(10,10)(30,15)(50,25)}
\put(0,0.1){\qbezier(10,10)(30,15)(50,25)}
\put(0,0.1){\qbezier(10,10)(30,15)(50,25)}

\put(0,-0.2){\qbezier(50,25)(52,33)(60,36)}
\put(0,-0.1){\qbezier(50,25)(52,33)(60,36)}
\put(0,0){\qbezier(50,25)(52,33)(60,36)}
\put(0,0.1){\qbezier(50,25)(52,33)(60,36)}
\put(0,0.2){\qbezier(50,25)(52,33)(60,36)}

\put(0,-0.2){\qbezier(60,36)(73,42)(80,47.5)}
\put(0,-0.1){\qbezier(60,36)(73,42)(80,47.5)}
\put(0,0){\qbezier(60,36)(73,42)(80,47.5)}
\put(0,0.1){\qbezier(60,36)(73,42)(80,47.5)}
\put(0,0.2){\qbezier(60,36)(73,42)(80,47.5)}

\color{red}
\put(0,-0.2){\qbezier(80,47.5)(90,55)(110,70)}
\put(0,-0.1){\qbezier(80,47.5)(90,55)(110,70)}
\put(0,0){\qbezier(80,47.5)(90,55)(110,70)}
\put(0,0.1){\qbezier(80,47.5)(90,55)(110,70)}
\put(0,0.2){\qbezier(80,47.5)(90,55)(110,70)}

\put(50,25){\circle*{1.2}}

\color{black}
\put(30,9){\line(0,1){2}}
\put(28,5.5){$K$}

\put(50,9){\line(0,1){2}}
\put(0,0){\qbezier[20](50,10)(50,17.5)(50,25)}
\put(49,5.5){$z$}

\put(60,9){\line(0,1){2}}
\put(59,5.5){$\gzc$}

\put(80,9){\line(0,1){2}}
\put(0,0){\qbezier[45](80,10)(80,28.75)(80,47.5)}
\put(79,5.5){$\xi$}

\put(37,35.5){\vector(1,-3){4.7}}
\put(39.8,39){\vector(3,-1){15}}
\put(30,38){$v(x)$}

\put(92,41){\vector(-1,3){4}}
\put(87,38){$x-K$}

\end{picture}}

\put(20,10){\small{{\bf Figure 11.}  The value function
$v$ if the problem parameters are as in}}
\put(42,5){\small{Case~(III) of Lemma~\ref{lem:psi}
and $z \in \mbox{} ]\zminus, \gzc[$.}}
\end{picture}

If the problem parameters are as in Case~(IV) of 
Lemma~\ref{lem:psi}, then the function $w = w (\cdot; z)$
given by
(\ref{w-}), (\ref{Gamma}) is such that
\begin{align}
w(\gzo ; \zplus) & \equiv \Gamma (\zplus) \psi
(\gzo; \zplus) = \gzo - K \nonumber \\
\text{and} \quad
w \bigl( \upalpha (\zplus) ; \zplus \bigr) & \equiv \Gamma
(\zplus) \psi \bigl( \upalpha (\zplus) ; \zplus \bigr) =
\upalpha (\zplus) - K \label{w-zplus}
\end{align}
(see Lemma~\ref{lem:HJB-ineq}.(iii) and Figure~12).
This observation suggests the possibility for the stopping
time
\ben
\overline{\tau}^\star = \inf \bigl\{ t \geq 0 \mid \ X_t \in
[\gzo, \gamma] \mbox{} \cup \mbox{} [\zeta, \infty[ \bigr\}
, \label{optautil}
\een
where $\gamma = \gamma (z) < \zeta = \zeta (z)$
are constants, to be optimal.
In such a context, we expect that the value function
$v$ should identify with the function
\ben
\overline{w} (x) = \overline{w} (x;z) = \begin{cases}
\frac{1}{n} \gzo ^{-n+1} x^n , & \text{if } x \in \mbox{} ]0,
\gzo[ , \\ \Cle x^n + \Dle x^m , & \text{if } x \in \mbox{}
]\gamma, z] , \\ \Cri x^n + \Dri x^m , & \text{if } x \in
\mbox{} ]z, \zeta[ , \\ x-K , & \text{if } x \in[\gzo, \gamma]
\cup [\zeta, \infty[ , \end{cases} \label{w-complex1}
\een
for some $\Cle, \Dle, \Cri, \Dri \in \bbr$ (see
Figure~13).
We suppress the actual dependence of $\Cle$, $\Dle$,
$\Cri$, $\Dri$ on $z$ because we will not use variational
arguments in the analysis of this case.
To determine the constants $\Cle$, $\Dle$, $\Cri$, $\Dri$ and
the free-boundary points $\gamma$, $\xi$, we first require
that $\overline{w}$ should be $C^1$ at $\gamma$ and
$\zeta$, which is suggested by the ``principle of smooth fit''.
This requirement yields the expressions
\begin{align}
\Cle & = -\frac{1}{n-m} \bigl[ (m-1) \gamma - mK \bigr]
\gamma^{-n} , \quad
\Cri = -\frac{1}{n-m} \bigl[ (m-1) \zeta - mK \bigr] \zeta^{-n}
, \label{CD32} \\
\Dle & = \frac{1}{n-m} \bigl[ (n-1) \gamma - nK \bigr]
\gamma^{-m} \quad \text{and} \quad
\Dri = \frac{1}{n-m} \bigl[ (n-1) \zeta
- nK \bigr] \zeta^{-m} . \label{CD42} 
\end{align}
We next require that $\overline{w}$ should be 
continuous at $z$ and satisfy the identity 
\be
(1+\beta) w_+' (z) = (1-\beta) w_-' (z)
\ee
that is associated with (\ref{VI2}).
These requirements yield the identities
\begin{align*}
\Cri & = \frac{n(1-\beta) - m(1+\beta)}{(n-m)(1+\beta)}
\Cle -\frac{2m\beta}{(n-m)(1+\beta)} \Dle z^{-(n-m)} \\
\text{and} \quad
\Dri & = \frac{2n\beta}{(n-m)(1+\beta)} \Cle z^{n-m} +
\frac{n(1+\beta) - m(1-\beta)}{(n-m)(1+\beta)} \Dle .
\end{align*}
Using the expressions in (\ref{CD32}), (\ref{CD42})
to substitute for $\Cle$, $\Dle$, $\Cri$, $\Dri$, 
we obtain the system of equations
\begin{align}
\bigl[ (m-1) \zeta - mK \bigr] z^n \zeta^{-n} - \frac{n(1-\beta)
- m(1+\beta)}{(n-m)(1+\beta)} \bigl[ (m-1) \gamma - mK \bigr]
z^n \gamma^{-n} & \nonumber \\
- \frac{2m\beta}{(n-m)(1+\beta)} \bigl[ (n-1) \gamma -
nK \bigr] z^m \gamma^{-m} & = 0 \label{GH2} \\
\text{and} \quad
\bigl[ (n-1) \zeta - nK \bigr] z^m \zeta^{-m} + \frac{2n\beta}
{(n-m)(1+\beta)} \bigl[ (m-1) \gamma - mK \bigr] z^n \gamma^{-n}
& \nonumber \\
- \frac{n(1+\beta) - m(1-\beta)}{(n-m)(1+\beta)} \bigl[ (n-1)
\gamma - nK \bigr] z^m \gamma^{-m} & = 0 . \label{GH1}
\end{align}
By (a) subtracting (\ref{GH2}) from (\ref{GH1}) and (b)
solving (\ref{GH1}) for  $\bigl[ (m-1) \gamma - mK \bigr]
z^n \gamma^{-n}$ and substituting for the resulting
expression in (\ref{GH2}), we obtain the system of equations
\begin{align}
G (\gamma, \zeta; z) = 0 \quad \text{and} \quad 
H (\gamma, \zeta; z) = 0 , \label{GH} 
\end{align}
which is equivalent to (\ref{GH2}) and (\ref{GH1}), where 
\begin{align}
G(x, y; z) = \mbox{} & \bigl[ (n-1) y - nK \bigr] z^m y^{-m}
- \bigl[ (m-1) y - mK \bigr] z^n y^{-n} \nonumber \\
& \mbox{} - \bigl[ (n-1) x - nK \bigr] z^m x^{-m}
+ \bigl[ (m-1) x - mK \bigr] z^n x^{-n} \label{G} \\
\text{and} \quad
H(x, y; z) = \mbox{} & y^{-n} F(y;z) - \frac{1-\beta}{1+\beta}
\bigl[ (n-1) x - nK \bigr] x^{-m} , \label{H}
\end{align}
in which expressions, $F$ is the function defined by (\ref{F}).
   
We prove the final result of the paper in
Section~\ref{sec:proofs}.

\begin{thm} \label{prop:gz}
Consider the optimal stopping problem defined by
(\ref{XGBM}), (\ref{vEx}) and suppose that
Assumption~\ref{assm} holds true.
Also, suppose that $\beta \in \mbox{} ]0,1[$ (Case~(IV)
of  Lemma~\ref{lem:psi}).
The following statements hold true:
\smallskip

\noindent {\rm (a)}
The system of equations (\ref{GH}) has a unique solution
$(\gamma, \zeta)$ such that $\gzo < \gamma < z
< \zeta$ if and only if $z > \zplus$, where $\zplus$ is
as in Lemma~\ref{lem:HJB-ineq}.(iii).
\smallskip

\noindent {\rm (b)}
Given any $z > \zplus$ and the associated solution
$(\gamma, \zeta)$ to the system of equations
(\ref{GH}), the function $\overline{w}$ defined by
(\ref{w-complex1}), for $\Cle, \Dle, \Cri, \Dri > 0$
given by (\ref{CD32})--(\ref{CD42}) 
identifies with the value function $v$ of the
discretionary stopping  problem and the
$(\F_t)$-stopping time $\overline{\tau}^\star$
defined by \eqref{optautil} is optimal.
\end{thm}

In the context of the previous result and the part
of Theorem~\ref{propSol1} addressing the case when
the problem parameters are as in Case~(IV) of 
Lemma~\ref{lem:psi} (see (\ref{a-3}) in particular),
we can see that Figure~10 transforms ``continuously''
into Figure~12 and then into Figure~13 as $z$ increases
from 0 to $\infty$, thanks to the identity
$\overline{w} (\cdot; \zminus) = w(\cdot; \zminus)$,
which follows from (\ref{w-zplus}).
\bigskip

\begin{picture}(160,110)
\put(20,15){\begin{picture}(120,95)

\put(0,0){\line(1,0){120}}
\put(120,0){\line(0,1){95}}
\put(0,0){\line(0,1){95}}
\put(0,95){\line(1,0){120}}

\put(10,10){\vector(1,0){100}}
\put(10,10){\vector(0,1){80}}
\put(6,88){$y$}
\put(108,6){$x$}

\qbezier(30,10)(70,40)(110,70)

\color{blue}
\put(0,-0.2){\qbezier(10,10)(30,11)(50,25)}
\put(0,-0.1){\qbezier(10,10)(30,11)(50,25)}
\put(0,0){\qbezier(10,10)(30,11)(50,25)}
\put(0,0.1){\qbezier(10,10)(30,11)(50,25)}
\put(0,0.1){\qbezier(10,10)(30,11)(50,25)}

\put(0,-0.2){\qbezier(50,25)(58,31)(64,40)}
\put(0,-0.1){\qbezier(50,25)(58,31)(64,40)}
\put(0,0){\qbezier(50,25)(58,31)(64,40)}
\put(0,0.1){\qbezier(50,25)(58,31)(64,40)}
\put(0,0.2){\qbezier(50,25)(58,31)(64,40)}

\put(0,-0.2){\qbezier(64,40)(74,43)(85,51.25)}
\put(0,-0.1){\qbezier(64,40)(74,43)(85,51.25)}
\put(0,0){\qbezier(64,40)(74,43)(85,51.25)}
\put(0,0.1){\qbezier(64,40)(74,43)(85,51.25)}
\put(0,0.2){\qbezier(64,40)(74,43)(85,51.25)}

\color{red}
\put(0,-0.2){\qbezier(85,51.25)(95,58.75)(110,70)}
\put(0,-0.1){\qbezier(85,51.25)(95,58.75)(110,70)}
\put(0,0){\qbezier(85,51.25)(95,58.75)(110,70)}
\put(0,0.1){\qbezier(85,51.25)(95,58.75)(110,70)}
\put(0,0.2){\qbezier(85,51.25)(95,58.75)(110,70)}

\put(50,25){\circle*{1.2}}

\color{black}
\put(0,0){\qbezier(85,51.25)(95,59)(100,75)}

\put(30,9){\line(0,1){2}}
\put(28,5.5){$K$}

\put(85,9){\line(0,1){2}}
\put(0,0){\qbezier[50](85,10)(85,30.6)(85,51.25)}
\put(84,5.5){$a$}

\put(64,9){\line(0,1){2}}
\put(0,0){\qbezier[40](64,10)(64,25)(64,40)}
\put(63,5.5){$\zplus$}

\put(50,9){\line(0,1){2}}
\put(0,0){\qbezier[20](50,10)(50,17.5)(50,25)}
\put(49,5.5){$\gzo$}

\put(39,43){\vector(1,-1){14.5}}
\put(30,44){$v(x)$}

\put(89,68){\vector(1,0){8.1}}
\put(80,67){$\psi(x)$}

\end{picture}}

\put(20,10){\small{{\bf Figure 12.} The value function
$v$ if the problem parameters are as in}}
\put(42,5){\small{Case~(IV) of Lemma~\ref{lem:psi}
and $z = \zplus$.}}
\end{picture}

\begin{picture}(160,110)
\put(20,15){\begin{picture}(120,95)

\put(0,0){\line(1,0){120}}
\put(120,0){\line(0,1){95}}
\put(0,0){\line(0,1){95}}
\put(0,95){\line(1,0){120}}

\put(10,10){\vector(1,0){100}}
\put(10,10){\vector(0,1){80}}
\put(6,88){$y$}
\put(108,6){$x$}

\qbezier(30,10)(70,40)(110,70)

\color{blue}
\put(0,-0.2){\qbezier(10,10)(30,11)(50,25)}
\put(0,-0.1){\qbezier(10,10)(30,11)(50,25)}
\put(0,0){\qbezier(10,10)(30,11)(50,25)}
\put(0,0.1){\qbezier(10,10)(30,11)(50,25)}
\put(0,0.1){\qbezier(10,10)(30,11)(50,25)}

\put(0,-0.2){\qbezier(60,32.5)(65,36)(75,50)}
\put(0,-0.1){\qbezier(60,32.5)(65,36)(75,50)}
\put(0,0){\qbezier(60,32.5)(65,36)(75,50)}
\put(0,0.1){\qbezier(60,32.5)(65,36)(75,50)}
\put(0,0.2){\qbezier(60,32.5)(65,36)(75,50)}

\put(0,-0.2){\qbezier(75,50)(86,53)(100,62.5)}
\put(0,-0.1){\qbezier(75,50)(86,53)(100,62.5)}
\put(0,0){\qbezier(75,50)(86,53)(100,62.5)}
\put(0,0.1){\qbezier(75,50)(86,53)(100,62.5)}
\put(0,0.2){\qbezier(75,50)(86,53)(100,62.5)}

\color{red}
\put(0,-0.2){\qbezier(50,25)(55,28.75)(60,32.5)}
\put(0,-0.1){\qbezier(50,25)(55,28.75)(60,32.5)}
\put(0,0.1){\qbezier(50,25)(55,28.75)(60,32.5)}
\put(0,0.2){\qbezier(50,25)(55,28.75)(60,32.5)}

\put(0,-0.2){\qbezier(100,62.5)(105,66.25)(110,70)}
\put(0,-0.1){\qbezier(100,62.5)(105,66.25)(110,70)}
\put(0,0.1){\qbezier(100,62.5)(105,66.25)(110,70)}
\put(0,0.2){\qbezier(100,62.5)(105,66.25)(110,70)}

\color{black}
\put(30,9){\line(0,1){2}}
\put(28,5.5){$K$}

\put(50,9){\line(0,1){2}}
\put(0,0){\qbezier[20](50,10)(50,17.5)(50,25)}
\put(49,5.5){$\gzo$}

\put(60,9){\line(0,1){2}}
\put(0,0){\qbezier[30](60,10)(60,21.25)(60,32.5)}
\put(59,5.5){$\gamma$}

\put(75,9){\line(0,1){2}}
\put(0,0){\qbezier[45](75,10)(75,30)(75,50)}
\put(74,5.5){$z$}

\put(100,9){\line(0,1){2}}
\put(0,0){\qbezier[55](100,10)(100,36.25)(100,62.5)}
\put(99,5.5){$\zeta$}

\put(39,44){\vector(4,-1){24}}
\put(30,44){$v(x)$}

\end{picture}}

\put(20,10){\small{{\bf Figure 13.} The value function
$v$ if the problem parameters are as in}}
\put(42,5){\small{Case~(IV) of Lemma~\ref{lem:psi}
and $z \in \mbox{} ]\zplus, \infty[$.}}
\end{picture}

\section{Proofs of Theorems~\ref{propSol1}-\ref{prop:gz}}
\label{sec:proofs}

If we denote by $g$ any of the functions $w$,
$\underline{w}$ or $\overline{w}$, defined by
(\ref{w-}), (\ref{wJ}) and (\ref{w-complex1}), respectively,
then
\be
\lim _{y \rightarrow 0} \frac{g(y)}{\varphi (y)} =
\lim _{y \rightarrow 0} \frac{(y - K)^+}{\varphi (y)}
= \lim _{y \rightarrow \infty} \frac{(y - K)^+} {\psi(y)} =
\lim _{y \rightarrow \infty} \frac{g(y)}{\psi(y)} = 0 .
\ee
In view of this observation and
Theorem~\ref{VThm}.(III)-(IV), we can see that we will
prove Theorems~\ref{propSol1}, \ref{prop:xi}
and~\ref{prop:gz} if we establish the claims made
on the solvability of their associated free-boundary
problems as well as show that the corresponding
functions $w$, $\underline{w}$ and $\overline{w}$
satisfy the variational inequality (\ref{VI1})--(\ref{VI2}).
\bigskip

\noindent
{\bf Proof of Theorem~\ref{propSol1}.}
By construction, $w$ is $C^2$ inside $]0, \infty[ \mbox{}
\setminus \{ a, z \}$ and $C^1$ at $a$ if $a \neq z$.
It is straightforward to verify that $w$ satisfies (\ref{VI2}).
In view of its structure, we will verify that $w$
satisfies (\ref{VI1}) if we prove that
\begin{gather}
x-K \leq w(x) \quad \text{for all } x < a , \label{VI-21} \\
\text{and} \quad
\half \sigma^2 x^2 w'' (x) + bx w' (x) - r w(x) \leq 0
\quad \text{for all } x > a . \label{VI-31}
\end{gather}
In view of (\ref{w-}) and (\ref{Gamma}), we can see
that (\ref{VI-21}) is equivalent to
\be
\frac{x-K}{\psi (x)} \leq \frac{a-K}{\psi (a)} \quad
\text{for all } x < a .
\ee
In the context of (\ref{a-1}) with $z < \gzb$ or
(\ref{a-2}) with $z \leq \zminus$ or (\ref{a-3}) with
$z \leq \zplus$, this inequality is equivalent to
(\ref{stop-ineq}), which is true thanks to
Lemma~\ref{lem:HJB-ineq}.
In the context of (\ref{a-1}) with $z \geq \gzb$ or
(\ref{a-2}) with $z \geq \gzc$, this inequality follows
immediately from the fact that the function
$x \mapsto (x-K) x^{-n}$ is strictly increasing in
$]0, \gzo[$.
On the other hand, (\ref{VI-31}) is equivalent to
$bx - r (x-K) \leq 0$ for all $x > a$, which is true
because, in all cases, $a > \gzc = \frac{rK}{r-b}$.
\mbox{}\hfill$\Box$
\bigskip

\noindent
{\bf Proof of Theorem~\ref{prop:xi}.}
Fix any $z \in \mbox{} ]0, \gzc[$.
Using the identity in (\ref{rb-nm}), we calculate
\begin{align}
\frac{\partial J}{\partial x} (x; z) & = - (n-1)(m-1) A z^{n-m}
x^{-n-1} (x - \gzc) \left[ \left( \frac{x}{z} \right) ^{n-m} -1
\right] \nonumber \\
& \begin{cases} < 0 , & \text{if } x \in \mbox{} ]z,
\gzc[ , \\ > 0 , & \text{if } x \in \mbox{} ]\gzc, \infty[ .
\end{cases} \label{J_x}
\end{align}
Combining this result with the observations that
\be
J(z;z) = 0 \quad \text{ and } \quad
\lim _{x \rightarrow \infty} J(x; z) = \infty ,
\ee
we can see that there exists a unique $\upxi (z) \in
\mbox{} ]\gzc, \infty[$ such that 
\ben
J(x; z) \begin{cases} 
< 0 & \text{for all } x \in \mbox{} ]z, \upxi (z)[ , \\
> 0 & \text{for all } x \in \mbox{} ]\upxi(z), \infty[
, \end{cases}
\quad \text{and} \quad
\frac{\partial J}{\partial x} (x;z) > 0 \quad \text{for all } 
x \geq \upxi (z) . \label{Jxi}
\een

We next consider the function $h$ defined by
\be
h(x;z) = C(z) x^n + D(z) x^m - (x-K) , \quad
(x,z) \in \mbox{} ]0, \infty[ \mbox{} \times \mbox{}
]0, \gzc[ ,
\ee
where $C$ and $D$ are given by (\ref{CD}) for
$\xi = \upxi (z)$, and we note that
\ben
h(z;z) = 0 \quad \text{and} \quad h \bigl( \upxi (z)
; z \bigr) = \frac{\partial h}{\partial x} \bigl( \upxi (z)
; z \bigr) = 0 . \label{h-BCs}
\een
The calculation
\be
\half \sigma^2 x^2 \, \frac{\partial^2}{\partial x^2} \left(
x \, \frac{\partial h}{\partial x} (x;z) \right) + bx \,
\frac{\partial}{\partial x} \left( x \, \frac{\partial h}{\partial x}
(x;z) \right) - r x \, \frac{\partial h}{\partial x} (x;z) =
(r-b) x > 0
\ee
and the maximum principle imply that the function
$x \mapsto x \, \frac{\partial h}{\partial x} (x;z)$ has
no positive maximum.
Combining this observation with the fact that
$\frac{\partial h}{\partial x} (\upxi (z); z) = 0$, we
can see that, if we define
\be
\bar{x} (z) = \inf \left\{ x \in [z, \upxi (z)] \ \Big| \ \
\frac{\partial h} {\partial x} (x;z) = 0 \right\} ,
\ee
then $\frac{\partial h}{\partial x} (x;z) \leq 0$ for all
$x \in [\bar{x} (z), \upxi (z)]$.
Since $h(z;z) = h \bigl( \upxi (z) ; z \bigr) = 0$ and
$h (\cdot ;z)$ is not constant, it is not possible that
either $\frac{\partial h}{\partial x} (x;z) \leq 0$ for all
$x \in [z, \upxi (z)]$ or $\frac{\partial h}{\partial x} (x;z)
\geq 0$ for all $x \in [z, \upxi (z)]$.
Therefore,
\ben
\overline{x} \in \mbox{} ]z,\upxi (z)[ \quad \text{and}
\quad \frac{\partial h}{\partial x} (x;z) \begin{cases}
> 0 & \text{for all } x \in [z, \overline{x} (z)[ , \\  \leq 0
& \text{for all } x \in \mbox{} ]\overline{x} (z), \upxi (z)[
. \end{cases}  \label{h'}
\een
Furthermore, the function $\underline{w}$ defined by
(\ref{wJ}) for $\xi = \upxi (z)$ is such that
\ben
\underline{w} (x) > x-K \quad \text{ for all } x \in \mbox{}
]z, \upxi (z)[ . \label{w-ominus1}
\een

To derive the monotonicity of $\upxi$, we first note that
the inequality $z \, \frac{\partial h}{\partial x} (z;z) > 0$,
which follows from (\ref{h'}), the identity $h(z;z) = 0$
and the expression for $D$ given by (\ref{CD}) with
$\xi = \upxi (z)$, imply that
\be
\bigl[ (n-1) \upxi (z) - nK \bigr] \upxi ^{-m} (z) < \bigl[
(n-1) z - nK \bigr] z^{-m} .
\ee
In view of (\ref{Jeqn}) and the first expression in
(\ref{J}), we can see that this inequality
is equivalent to
\be
\bigl[ (m-1) \upxi (z) - mK \bigr] z^{n-m} \upxi ^{-n} (z)
< \bigl[ (m-1) z - mK \bigr]  z^{-m} .
\ee
In view of (\ref{Jeqn}), it follows that
\begin{align}
\frac{\partial J}{\partial z} \bigl( \upxi (z); z \bigr) & = 
- (n-m) A z^{-1} \Bigl( \bigl[  (m-1) \upxi (z) - mK \bigr]
z^{n-m} \upxi ^{-n} (z) - \bigl[ (m-1) z - mK \bigr] z^{-m} 
\Bigr) \nonumber \\
& > 0 . \nonumber
\end{align}
Differentiating the identity $J \bigl( \upxi (z) ; z \bigr) = 0$
with respect to  $z$ and using this inequality, along with
(\ref{Jxi}), we obtain
\be
\upxi' (z) = -  \frac{\frac{\partial J}{\partial z} \bigl(
\upxi (z) ; z \bigr)} {\frac{\partial J}{\partial x} \bigl(
\upxi (z) ; z \bigr) } < 0 ,
\ee
which proves that $\upxi$ is strictly decreasing.
Furthermore, the limits in (\ref{xi-lims}) hold true
thanks to (\ref{J_x}) and the fact that
\be
0 = \lim _{z \downarrow 0} J \bigl( \upxi (z); z \bigr)
= A \lim _{z \downarrow 0} \bigl[ (n-1) \upxi (z) - nK \bigr]
\upxi^{-m} (z) ,
\ee
which follows from the first expression in (\ref{J}) and
the fact that $\upxi (z) > \gzc$ for all $z < \gzc$.

To complete the proof, we fix any $z \in \mbox{}
]\zminus, \gzc[$.
By construction, the function $\underline{w}$
defined by (\ref{wJ}) for $\xi = \upxi (z)$ is continuous,
$C^1$ inside $]0, \infty[ \mbox{} \setminus \{ z \}$
and $C^2$ inside $]0, \infty[ \mbox{} \setminus
\{ z, \xi \}$.
In view of its structure, we will verify that it satisfies
the variational inequality (\ref{VI1})--(\ref{VI2})
if we prove that
\begin{gather}
(1+\beta) \underline{w}_+' (z) \leq (1-\beta)
\underline{w}_-' (z) ,
\label{VI(-)-1} \\
x-K \leq \underline{w} (x) \quad \text{for all }
x < \upxi (z) , \label{VI(-)-2} \\
\text{and} \quad
\half \sigma^2 x^2 \underline{w}'' (x) + bx
\underline{w}' (x) - r \underline{w}(x) \leq 0
\quad \text{for all } x > \upxi (z) . \label{VI(-)-3}
\end{gather}
In view of the definition (\ref{wJ}) of $\underline{w}$
and the identity  $h(z,z) = 0$, we can see
that (\ref{VI(-)-1}) is equivalent to 
\be
(n-m) (1+\beta) C(z) z^n + m (1+\beta) (z-K) \leq
n (1-\beta) (z-K) .
\ee
Furthermore, using the expressions for $A$ and
$C(z)$ given by (\ref{A}) and (\ref{CD}), we can see
that  this inequality is equivalent to
\ben
- \bigl[ (m-1) \upxi (z) - mK \bigr] z^n \upxi^{-n}
(z) - (n-m) A (z-K) \leq 0 \quad
\stackrel{(\ref{Jeqn})}{\Leftrightarrow}
\quad F \bigl( \upxi (z); z \bigr) \geq 0 . \label{Fxi}
\een
If $z \in [\gc \gzc, \gzc[$, then this inequality follows
immediately from (\ref{F>0-ii}) and the fact that
$\upxi (z) > \gzc$.
If $z \in \mbox{} ]\zminus, \gc \gzc[$, then the
conclusions in Remark~\ref{rem:zminus...} and
(\ref{h-BCs}) imply that
\be
C \bigl( \gz (z) \bigr) = A \Gamma (z) , \quad
D \bigl( \gz (z) \bigr) = B(z) \Gamma (z)
\quad \text{and} \quad
\upxi \bigl( \gz (z) \bigr) = \upalpha (z) ,
\ee
where $\upalpha (z)$ and $\gz (z) \in \mbox{}
]z, \gzc[$ are as in Lemma~\ref{lem:a(z)}.(ii) and
Lemma~\ref{lem:HJB-ineq}.(ii), respectively.
The last of these identities and the fact that
$\upxi$ is strictly decreasing imply that
$\upxi (z) > \upalpha (z)$, and (\ref{Fxi}) follows
from (\ref{gotha-ii2}) in Lemma~\ref{lem:a(z)}.(ii).

Recalling that $\gzc < \gzo = \frac{nK}{n-1}$ (see
(\ref{cp-III}) in Lemma~\ref{lem:psi})), we can see
that
\be
\frac{d}{dx} \bigl[ \underline{w} (x) - (x-K) \bigr]
< (n-1) (z-\gzo) z^{-1} < 0 \quad \text{for all }
x < z \in \mbox{} ]\zminus, \gzc[ .
\ee
This result, the identity $\underline{w} (z)
- (z-K) = 0$ and (\ref{w-ominus1}) imply that
(\ref{VI(-)-2}) holds true.
Finally, (\ref{VI(-)-3}) is equivalent to $bx - r (x-K)
\leq 0$ for all $x > \upxi (z)$, which is true
because $\upxi (z) > \gzc = \frac{rK}{r-b}$.
\mbox{}\hfill$\Box$
\bigskip

\noindent
{\bf Proof of Theorem~\ref{prop:gz}.}
In view of the inequality $\gzc < \gzo$ (see (\ref{cp-IV}))
and the definition (\ref{G}) of $G$, we can see that
\begin{align}
& G(x, x; z) = 0 , \qquad
\lim _{y \rightarrow \infty} G (x, y; z) = \infty
\nonumber \\
\text{and} \quad
\frac{\partial G}{\partial y} (x, y; z) & = - (n-1) (m-1)
(y - \gzc) \left[ \left( \frac{y}{z} \right) ^{n-m}
-1 \right] z^n y^{-n-1} \nonumber \\
& \begin{cases} < 0 , & \text{for all } y \in \mbox{}
\bigl] \gzo, z \bigr[ , \\ > 0  , & \text{for all } y
> z . \end{cases} \nonumber
\end{align}
It follows that, given any $z > \gzo$, there exists a unique
function $L(\cdot; z) : [\gzo, z[ \mbox{} \rightarrow \mbox{}
]z, \infty[$, such that 
\ben
G \bigl( x, L(x; z); z \bigr) = 0 \quad \text{and} \quad
z < L(x; z) \quad \text{for all } x \in [\gzo, z[ .
\label{L}
\een
In particular, this function is such that
\ben
G (x, y; z) \begin{cases} < 0 , & \text{for all } \gzo \leq
x < y < L(x; z) , \\ > 0 , & \text{for all } \gzo \leq x <
L(x; z) < y , \end{cases} \quad \text{and} \quad \lim
_{x \rightarrow z} L(x; z) = z . \label{G-signs}
\een
Furthermore, differentiating the identity $G \bigl( x,
L(x; z) ;z \bigr) = 0$ with respect to $x$, we obtain
\ben
\frac{\partial L}{\partial x} (x; z) = - \frac{(x - \gzc) \left[
1 - \left( \frac{x}{z} \right) ^{n-m} \right] x^{-n-1}}
{\bigl( L(x; z) - \gzc \bigr) \left[ \left( \frac{L(x; z)}{z}
\right) ^{n-m} -1 \right] L^{-n-1} (x; z)} < 0
\text{ \ for all } x \in \mbox{} ]\gzo, z[ . \label{L'}
\een

In view of (\ref{F(z,z)}), the definition (\ref{H}) of $H$
and the limit in (\ref{G-signs}), we can see that
\be
\lim _{x \rightarrow z} H \bigl( x, L(x; z); z \bigr) =
H(z,z;z) = - \frac{2\beta}{1+\beta} z^{-m+1} < 0 .
\ee
On the other hand, we use (\ref{Fa()}), (\ref{L}),
(\ref{L'}) and the inequality $\gzc < \gzo$ to
obtain
\begin{align}
\frac{\partial}{\partial x} H & \bigl( x, L(x; z); z \bigr)
\nonumber \\
= \mbox{} & (n-1) (m-1) \frac{1-\beta}{1+\beta}
(x - \gzc) x^{-m-1}
\nonumber \\
& \times \left[ \frac{\bigl[ n(1-\beta) - m(1+\beta) \bigr]
L^{n-m} (x; z) + 2n\beta z^{n-m}} {(n-m) (1-\beta)}
\frac{1 - \left( \frac{x}{z} \right) ^{n-m}}
{\left( \frac{L(x; z)}{z} \right) ^{n-m} - 1} x
^{-n+m} + 1 \right] \nonumber \\
< \mbox{} & 0 \qquad \quad \text{for all } x
\in \mbox{} ]\gzo, z[ . \nonumber
\end{align}
These calculations imply that
\ben
\text{there exists a unique } x^* \in \mbox{} ]\gzo, z[
\text{ such that } H \bigl( x^*, L(x^*; z); z \bigr)
= 0 \label{gamma*}
\een
if and only if
\ben
H \bigl( \gzo, L (\gzo; z); z \bigr) = L^{-n} (\gzo; z)
F \bigl( L(\gzo; z); z \bigr) > 0 . \label{H-ineq}
\een

The analysis leading to (\ref{L})  and
(\ref{gamma*})--(\ref{H-ineq}) imply that the system
of equations (\ref{GH}) has a unique solution
$(\gamma, \zeta) $ such that $\gzo < \gamma < z
< \zeta$ if and only if (\ref{H-ineq}) holds true,
in which case, $(\gamma, \zeta) = \bigl( x^*,
L(x^*; z) \bigr)$.
We can show that (\ref{H-ineq}) is equivalent to
$z > \zplus$, where $\zplus$ is as in
Lemma~\ref{lem:HJB-ineq}.(iii), as follows.

If the problem's parameters are such that
$n > \frac{1+\beta}{1-\beta}$ and $z \geq \gzb$,
then (\ref{cp-IV}) in Lemma~\ref{lem:psi},
(\ref{gotha-iii1})--(\ref{gotha-iii2}) in
Lemma~\ref{lem:a(z)} and (\ref{L}) imply that
(\ref{H-ineq}) holds true.
On the other hand, if the problem's parameters
are such that either $n > \frac{1+\beta}{1-\beta}$
and $z \in \mbox{} ]\gzo, \gzb[$ or $n \leq
\frac{1+\beta}{1-\beta}$ and $z > \gzo$,
then (\ref{gotha-iii1})--(\ref{gotha-iii2}) imply that
(\ref{H-ineq}) holds true if and only if $L (\gzo; z)
> \upalpha (z)$, where $\upalpha (z) > z$ is the
unique solution to the equation  $F(x;z) = 0$
derived in Lemma~\ref{lem:a(z)}.(iii).
In view of (\ref{G-signs}), we can see that the
inequality  $L (\gzo; z) > \upalpha (z)$ is equivalent
to
\begin{align}
G \bigl( \gzo, \upalpha (z); z \bigr) = \mbox{} &
\bigl[ (n-1) \upalpha (z) - nK \bigr] z^m \upalpha^{-m}
(z) \nonumber \\
& - \bigl[ (m-1) \upalpha (z) - mK \bigr] z^n
\upalpha^{-n} (z) - \frac{n-m}{n} z^n \gzo ^{-n+1}
< 0 . \label{G-ineq}
\end{align}
Using the identity $F \bigl( \upalpha (z); z \bigr) = 0$
to eliminate the term $\bigl[ (m-1) a - mK \bigr]$
and the identity $B(z) z^m + Az^n = z^n$ to simplify,
we derive the expression
\be
G \bigl( \gzo, \upalpha (z); z \bigr) B(z) \upalpha^m (z)
= \bigl[ (n-1) \upalpha (z) - nK \bigr] z^n - \frac{n-m}{n}
z^n \gzo ^{-n+1} B(z) \upalpha^m (z) .
\ee
Similarly, we can eliminate the term $\bigl[ (n-1)
\upalpha (z) - nK \bigr]$ to obtain
\be
G \bigl( \gzo, \upalpha (z); z \bigr) A \upalpha^n (z)
= - \bigl[ (m-1) \upalpha (z) - mK \bigr] z^n -
\frac{n-m}{n} z^n \gzo ^{-n+1} A \upalpha^n (z) .
\ee
In the case that we consider here (either
$n > \frac{1+\beta}{1-\beta}$ and $z \in \mbox{}
]\gzo, \gzb[$ or $n \leq \frac{1+\beta}{1-\beta}$
and $z > \gzo$), the fact that $z < \upalpha (z)$
(see (\ref{gotha-iii1}) in Lemma~\ref{lem:a(z)})
implies that
\be
\psi \bigl( \upalpha (z) ; z \bigr) = A \upalpha^n (z)
+ B(z) \upalpha^m (z) ,
\ee
while, the facts that $\gzo = \frac{nK}{n-1} < z$
imply that
\be
\frac{\gzo - K}{\psi (\gzo ; z)} = \frac{\gzo - K}{\gzo^n}
= \frac{1}{n} \gzo ^{-n+1} .
\ee
Therefore, adding up the last two expressions for $G$
yields
\be
G \bigl( \gzo, \upalpha (z); z \bigr) = (n-m) z^n \left[
\frac{\upalpha (z) - K}{\psi \bigl( \upalpha (z) ; z \bigr)}
- \frac{\gzo - K} {\psi (\gzo ; z)} \right] .
\ee
It follows that (\ref{G-ineq}) is equivalent to
$z > \zplus$, as required, thanks to
Lemma~\ref{lem:HJB-ineq}.(iii).

By construction, the function $\overline{w}$ defined by
(\ref{w-complex1}) is continuous, $C^1$ inside $]0, \infty[
\mbox{} \setminus \{ z \}$ and $C^2$ inside $]0, \infty[
\mbox{} \setminus \{ z, \xi \}$.
In view of its structure, we will verify that $\overline{w}$
satisfies the variational inequality (\ref{VI1})--(\ref{VI2})
if we prove that
\begin{gather}
x-K \leq \overline{w} (x) \quad \text{for all } x \in \mbox{}
]0, \gzo[ \mbox{} \cup \mbox{} ]\gamma, \zeta[ ,
\label{VI-23} \\
\text{and} \quad
\half \sigma^2 x^2 \overline{w}'' (x) + bx \overline{w}'
(x)  - r \overline{w} (x) \leq 0 \quad \text{for all } x \in
\mbox{} ]\gzo, \gamma[ \mbox{} \cup \mbox{}
]\zeta,\infty[ . \label{VI-33}
\end{gather}
To establish (\ref{VI-23}), we first note that
\be
\frac{d}{dx} \bigl[ \overline{w} (x) - (x-K) \bigr] = \gzo
^{-n+1} x^{n-1} -1 > 0 \quad \text{for all } x \in \mbox{}
]0, \gzo[ .
\ee
Combining this observation with the fact that
$\overline{w} (\gzo) - (\gzo - K) = 0$, we can
see that (\ref{VI-23}) holds true for all  $x \in \mbox{}
]0, \gzo[$.
On the other hand, the inequalities
\be
\frac{mK}{m-1} < \frac{nK}{n-1} = \gzo< \gamma
< \zeta
\ee 
and the expressions (\ref{CD32}), (\ref{CD42})
of $C_\ell$, $D_\ell$, $C_r$, $D_r$ imply that these
constants are all strictly positive.
Therefore, the restrictions of the function $x \mapsto
\overline{w} (x) - (x-K)$ to the intervals  $]\gamma, z[$
and $]z, \zeta[$ are both strictly convex. 
Combining this observation with the facts that 
\begin{align}
& \overline{w} (\gamma) - (\gamma - K) =
\left. \frac{d}{dx} \bigl[ \overline{w} (x) - (x - K)
\bigl] \right| _{x = \gamma} = 0 \nonumber \\
\text{and} \quad 
& \overline{w} (\zeta) - (\zeta - K) \bigr] =
\left. \frac{d}{dx} \bigl[ \overline{w} (x) - (x-K) \bigr]
\right| _{x = \zeta} = 0 , \nonumber
\end{align}
we can see that (\ref{VI-23}) also holds true for
$x \in  \mbox{} ]\gamma, \zeta[$. 
On the other hand, (\ref{VI-33}) is is equivalent
to $bx - r (x-K) \leq 0$ for all $x \in \mbox{}
]\gzo, \gamma[ \mbox{} \cup \mbox{} ]\zeta,\infty[$,
which is true because $\gzo > \gzc = \frac{rK}{r-b}$.
\mbox{}\hfill$\Box$
\bigskip

{\small 

}


\begin{thebibliography}{20}


\bibitem{AS}
{\sc Alvarez,\,L.H.R., and Salminen,\,P.} (2017),
Timing in the presence of directional predictability:
Optimal stopping of Skew Brownian Motion,
{\em Mathematical Methods of Operations Research\/},
{\bf 86}, 377--400.

\bibitem{AS98}
{\sc Assing,\,S., and Schmidt,\,W.M.} (1998),
{\em Continuous Strong Markov Processes in
Dimension One, A Stochastic Calculus Approach\/},
Springer.

\bibitem{BC}
{\sc Bassan, B. and Ceci, C.} (2002),
Optimal stopping problems with discontinuous reward:
regularity of the value function and viscosity solutions,
{\em Stochastics and Stochastics Reports\/}, {\bf 72}, 55--77.


\bibitem{Benes}
{\sc Bene\v{s},\,V.E.} (1992),
Some combined control and stopping problems,
paper presented at the {\em CRM Workshop on
Stochastic Systems\/}, Montr\'{e}al, November 1992.

\bibitem{BL}
{\sc Bensoussan, A., and Lions,\,J.L.} (1978),
{\em Applications of Variational Inequalities in
Stochastic Control\/}, North-Holland.

\bibitem{BS}
{\sc Borodin,\,A.N., and Salminen,\,P.} (2002),
{\em Handbook of Brownian Motion - Facts and Formulae\/},
Birkh\"{a}user.

\bibitem{CS}
{\sc Corns,\,T.R.A., and Satchell,\,S.E.} (2007),
Skew Brownian motion and pricing European
options,
{\em The European Journal of Finance\/},
{\bf 13}, 523--544.

\bibitem{COT}
{\sc Cox, A.M.G., Obloj, J. and Touzi, N.} (2019),
The Root solution to the multi-marginal embedding
problem: an optimal stopping and time-reversal
approach,
{\em Probability Theory and Related Fields\/},
{\bf 173}, 211--259.

\bibitem{CW}
{\sc Cox, A.M.G., and Wang, J.} (2013),
Root's barrier: construction, optimality and applications
to variance options,
{\em Annals of Applied Probabability\/}, {\bf 23},
859--894.

\bibitem{CM}
{\sc Crocce,\,F., and Mordecki,\,E.} (2014),
Explicit solutions in  one-sided optimal stopping problems
for one-dimensional diffusions,
{\em Stochastics: An International Journal of Probability
and Stochastic  Processes}, {\bf 86}, 491--509.

\bibitem{DZ}
{\sc Davis,\,M.H.A., and Zervos,\,M.} (1994),
A problem of singular stochastic control with discretionary stopping,
{\em Annals of Applied Probability}, {\bf 4}, 226-240.


\bibitem{DGS06a}
{\sc Decamps,\,M., Goovaerts,\,M., and Schoutens,\,W.}
(2006),
Self exciting threshold interest rates models,
{\em International Journal of Theoretical and Applied
Finance\/}, {\bf 9}, 1093--1122.

\bibitem{DGS06b}
{\sc Decamps,\,M., Goovaerts,\,M., and Schoutens,\,W.}
(2006),
Asymmetric skew  Bessel processes and their applications
to finance,
{\em Journal of Computational and Applied Mathematics\/},
{\bf 186}, 130--147.


\bibitem{ES91}
{\sc Engelbert,\,H.J., and Schmidt,\,W.} (1991),
Strong Markov continuous local martingales and solutions
of one-dimensional stochastic differential equations.\ III.,
{\em Mathematische Nachrichten\/},
{\bf 151}, 149--197.


\bibitem{F}
{\sc Friedman, A.} (2006),
{\em Stochastic Differential Equations and Applications\/},
reprint of the 1975 and 1976 original published in two volumes,
Dover Publications.

\bibitem{Ha}
{\sc H\"{a}m\"{a}l\"{a}inen,\,J.} (2015),
Portfolio selection with directional return estimates,
available at
SSRN: http://ssrn.com/abstract=2279823.


\bibitem{HS81}
{\sc Harrison,\,J.M., and Shepp,\,L.A.} (1981),
On skew Brownian motion, {\em Annals of Probability\/},
{\bf 9}, 309--313.

\bibitem{IMcK74}
{\sc It\^{o},\,K., and McKean,\,H.P.} (1974),
{\em Diffusion Processes and their Sample Paths\/},
Springer.


\bibitem{KS}
{\sc Karatzas,\,I., and Sudderth,\,W.D.} (1999),
Control and stopping of a diffusion process on an
interval, {\em Annals of Applied Probability},
{\bf 9}, 188--196. 

\bibitem{K}
{\sc Krylov,\,N.V.} (1980),
{\em Controlled Diffusion Processes},
Springer-Verlag.

\bibitem{Lam}
{\sc Lamberton, D.} (2009),
Optimal stopping with irregular reward functions,
{\em Stochastic Processes and their Applications},
{\bf 119}, 3253--3284.

\bibitem{LZ}
{\sc Lamberton,\,D., and Zervos,\,M.} (2013),
On the optimal stopping of a one-dimensional diffusion,
{\it Electronic Journal of Probability}, {\bfseries 18},
1--49.

\bibitem{Lej06}
{\sc Lejay,\,A.} (2006),
On the constructions of the skew Brownian motion,
{\em Probability Surveys\/}, {\bf 3}, 413--466.

\bibitem{Lon}
{\sc Lon,\,P.C.} (2011),
{\em Two explicitly solvable problems with discretionary
stopping\/},  PhD thesis, The London School of Economics
and Political Science (LSE),
{\tt http://etheses.lse.ac.uk/id/eprint/337}.

\bibitem{LRZ}
{\sc Lon,\,P.C., Rodosthenous,\,N., and Zervos,\,M.}
(2015),
On the optimal stopping of a skew geometric Brownian motion,
in {\em Modern Trends in Controlled Stochastic Processes:
Theory and Applications, Volume~II\/} (A.\,B.\,Piunovskiy, ed.),
Luniver Press, 231--245.

\bibitem{Mu99}
{\sc Murphy,\,J.J.} (1999),
{\em Technical Analysis of the Financial Markets: A
Comprehensive Guide to Trading Methods and
Applications\/}, Penguin.

\bibitem{NS}
{\sc Nilsen,\,W., and Sayit,\,H.} (2011),
No arbitrage in markets with bounces  and sinks,
{\em International Review of Applied Financial Issues
\& Economics}, {\bf 3}, 696--699.

\bibitem{O}
{\sc {\O}ksendal, B.} (2003),
{\em Stochastic Differential Equations. An Introduction with
Applications\/}, 6th edition, Springer.

\bibitem{OR}
{\sc {\O}ksendal, B. and Reikvam, K.} (1998),
Viscosity solutions of optimal stopping problems.
\emph{Stochastics and Stochastics Reports\/},
{\bf 62}, 285--301.

\bibitem{P}
{\sc Peskir,\,G.} (2007),
Principle of smooth fit and diffusions with angles,
{\em Stochastics An International Journal of
Probability and Stochastic Processes\/},
{\bf 79}, 293--302.


\bibitem{PS}
{\sc Peskir,\,G., and Shiryaev,\,A.N.} (2006),
{\em Optimal Stopping and Free-Boundary Problems\/},
Lectures in Mathematics, ETH Z\"{u}rich, Birkh\"{a}user.

\bibitem{RY}
{\sc Revuz,\,D., and Yor,\,M.} (1999),
{\em Continuous Martingales and Brownian Motion\/},
3rd edition, Springer-Verlag.

\bibitem{Ro}
{\sc Rossello,\,D.} (2012),
Arbitrage in skew Brownian motion models,
{\em Insurance: Mathematics and Economics},
{\bf 50}, 50--56.


\end{thebibliography}
\end{document}